\documentclass[reqno]{amsart}
\newtheorem{Pa}{Paper}[section]
\newtheorem{Tm}[Pa]{{\bf Theorem}}
\newtheorem{La}[Pa]{{\bf Lemma}}
\newtheorem{Cy}[Pa]{{\bf Corollary}}
\newtheorem{Rk}[Pa]{{\bf Remark}}
\newtheorem{Pn}[Pa]{{\bf Proposition}}
\newtheorem{Pb}[Pa]{{\bf Problem}}

\newcommand{\be}{{\bf e}}
\newcommand{\tP}{\widetilde{P}}
\newcommand{\tC}{\widetilde{C}}
\newcommand{\tE}{\widetilde{E}}
\newcommand{\te}{\widetilde{e}}
\newcommand{\tc}{\widetilde{c}}
\newcommand{\cE}{{\mathcal {E}}}
\newcommand{\cS}{{\mathcal {S}}}

\begin{document}
\newcommand{\HS}{{\bf HS}}
\newcommand{\gdots}{\mathinner{\mkern1mu\vbox{\kern7pt\hbox{.}}
\mkern2mu\raise2.5pt\hbox{.}\mkern2mu\raise5pt\hbox{.}\mkern1mu}}
\newcommand{\ssucceq}{\buildrel *\over  \succeq }
\newcommand{\C}{{\mathbb C}}
\newcommand{\B}{{\mathbb B}}
\newcommand{\D}{{\mathbb D}}
\newcommand{\T}{{\mathbb T}}
\newcommand{\R}{{\mathbb R}}
\newcommand{\q}{\eta}
\newcommand{\N}{{\mathbb N}}
\newcommand{\PP}{{\mathbb P}}

\title[Boundary Nevanlinna--Pick problem]
{Boundary Nevanlinna--Pick interpolation problems for
generalized Schur functions}

\author{Vladimir Bolotnikov and Alexander Kheifets}
\address{Department of Mathematics,
The College of William and Mary, Williamsburg VA 23187-8795, USA}
\email{vladi@math.wm.edu},
\address{Department of Mathematics,
University of Massachusetts Lowell, Lowell, MA, 01854, USA}
\email{Alexander\_Kheifets@uml.edu}

\date{}
\maketitle
\thispagestyle{empty}

\vskip 1truecm

\baselineskip=12pt
\begin{quote}
Three boundary Nevanlinna-Pick interpolation problems at finitely many
points are formulated for generalized Schur functions. For each problem,
the set of all  solutions is parametrized in terms of a linear fractional
transformation with a Schur class parameter.
\end{quote}

\baselineskip=15pt

\section{Introduction}
\setcounter{equation}{0}

The Schur class $\cS$ of complex-valued
analytic functions mapping the unit disk ${\mathbb D}$ into the closed
unit disk $\overline{\mathbb D}$ can be characterized in terms of positive
kernels as follows: a function $w$ belongs to $\cS$ if and only if the
kernel
\begin{equation}
K_w(z,\zeta):=\frac{1-\overline{w(\zeta)}w(z)}{1-\bar{\zeta}z}
\label{0.1}
\end{equation}
is positive definite on $\D$ (in formulas: $K_w\succeq 0$), i.e., if and
only if the Hermitian matrix
\begin{equation}
\left[K_w(z_j,z_i)\right]_{i,j=1}^n=
\left[\frac{1-\overline{w(z_i)}w(z_j)}{1-\bar{z}_iz_j}\right]_{i,j=1}^n
\label{0.2}
\end{equation}
is positive semidefinite for every choice of an integer $n$ and of $n$
points $z_1,\ldots,z_n\in\D$. The significance of this characterization
for interpolation theory is that it gives the necessity part in the
Nevanlinna-Pick interpolation theorem:  {\em given points $z_{1},
\dots, z_{n} \in\D$ and $w_{1}, \dots, w_{n} \in\C$, there exists
$w\in\cS$ with $w(z_j)=w_{j}$ for $j = 1, \ldots, n$ if and only
if the associated Pick matrix $P=\left[\frac{1- \overline{w}_iw_j}{1 -
\overline{z}_iz_j} \right]$ is positive semidefinite}.

\medskip

There are at least two obstacles to get an immediate boundary analogue
of the latter result just upon sending the points $z_1,\ldots,z_n$ in
\eqref{0.2} to the unit circle $\T$. Firstly, the boundary nontangential
(equivalently, radial) limits
\begin{equation}
w(t):=\lim_{z\to t}w(z)
\label{0.2a}
\end{equation}
exist at almost every (but not every) point $t$ on $\T$. Secondly,
although the nontangential limits
\begin{equation}
d_w(t):=\lim_{z\rightarrow t} \,
\frac{1-|w(z)|^2}{1-|z|^2}\ge 0 \quad (t\in\T)
\label{0.2b}
\end{equation}
exist at every $t\in\T$, they can be infinite. However, if $d_w(t)<\infty$,
then it is readily seen that the limit \eqref{0.2a} exists and is
unimodular. Then we can pass to limits in \eqref{0.2} to get the
necessity part of the following interpolation result:

\medskip

{\em Given points $t_1, \dots, t_n\in\T$ and numbers $w_{1}, \dots,
w_{n}$ and $\gamma_1,\ldots,\gamma_n$ such that
\begin{equation}
|w_i|=1\quad\mbox{and}\quad \gamma_i\ge 0\quad\mbox{for}\quad
i=1,\ldots,n,
\label{0.3}
\end{equation}
there exists $w\in\cS$ with
\begin{equation}
w(t_i)=w_i\quad\mbox{and}\quad d_w(t_i)\le \gamma_i\quad\mbox{for}\quad
i=1,\ldots,n
\label{0.4}
\end{equation}
if and only if the associated Pick matrix
\begin{equation}
P=[P_{ij}]_{i,j=1}^n\quad \mbox{with the entries}\quad
P_{ij}=\left\{\begin{array}{ccc} {\displaystyle
\frac{1-\overline{w}_iw_j}{1-\overline{t}_it_j}}&\mbox{for} & i\neq j\\
\gamma_i&\mbox{for}& i=j\end{array}\right.
\label{0.5}
\end{equation}
is positive semidefinite.}

\medskip

This result in turn, suggests the following well known
boundary Nevanlinna--Pick interpolation problem.
\begin{Pb}
Given points $t_1, \dots, t_n\in\T$ and numbers $w_{1}, \dots,
w_{n}$, $\gamma_1,\ldots,\gamma_n$ as in \eqref{0.3} and such that
the Pick matrix $P$ defined in \eqref{0.5} is positive semidefinite,
find all functions $w\in\cS$ satisfying interpolation conditions
\eqref{0.4}.
\label{P:0.1}
\end{Pb}
Note that assumptions \eqref{0.3} and $P\ge 0$ are not restrictive
since they are necessary for the problem to have a solution.

\medskip

The boundary Nevanlinna--Pick interpolation problem was worked out using
quite different approaches: the method of fundamental matrix
inequalities \cite{kov}, the recursive Schur algorithm \cite{dd84},
the Grassmannian approach \cite{ballhelton1}, via realization theory
\cite{bgr}, and via unitary extensions of partially defined isometries
\cite{b2, kh}. If $P$ is singular, then Problem \ref{P:0.1} has a unique
solution which is a finite Blaschke product
of degree $r={\rm rank} \, P$. If $P$ is positive definite, Problem
\ref{P:0.1} has infinitely many solutions that can
be described in terms of a linear fractional transformation with a free
Schur class parameter.

\medskip

Note that a similar problem with equality sign in the
second series of conditions in (\ref{0.4}) was considered in
\cite{Sarasonnp, Geor, boldym1}:
\begin{Pb}
Given the data as in Problem \ref{P:0.1},
find all functions $w\in\cS$ such that
\begin{equation}
w(t_i)=w_i\quad\mbox{and}\quad d_w(t_i)= \gamma_i\quad\mbox{for}\quad
i=1,\ldots,n
\label{0.5a}
\end{equation}
\label{P:0.2}
\end{Pb}
The solvability criteria for this modified
problem is also given in terms of the Pick matrix \eqref{0.5} but it is
more subtle: condition $P\ge 0$ is necessary (not sufficient, in general)
for the Problem \ref{P:0.2} to have a solution while the condition $P>0$
is
sufficient.

\medskip

The objective of this paper is to study the above problems in the
setting of {\em generalized Schur functions}. A function $w$ is called
{\em a generalized Schur function} if it is of the form
\begin{equation}
w(z)=\frac{S(z)}{B(z)},
\label{0.6}
\end{equation}
for some Schur function $S\in\cS$ and a finite Blaschke product $B$.
Without loss of generality we can (and will) assume that $S$ and $B$
in representation \eqref{0.6} have no common zeroes.
For a fixed integer $\kappa\ge 0$, we denote by $\cS_\kappa$ the
class of generalized Schur functions with $\kappa$ poles inside $\D$,
i.e., the class of functions of the form \eqref{0.6} with a Blaschke
product $B$ of degree $\kappa$. Thus, $\cS_\kappa$ is a class of functions
$w$ such that
\begin{enumerate}
\item $w$ is meromorphic in $\D$ and has $\kappa$ poles inside $\D$
counted with multiplicities.
\item $w$ is bounded on an annulus $\{z: \; \rho<|z|<1\}$ for some
$\rho\in(0, \, 1)$.
\item Boundary nontangential limits $w(t):={\displaystyle\lim_{z\to
t}w(z)}$ exist and satisfy $|w(t)|\leq 1$ for almost all $t\in\T$.
\end{enumerate}
It is clear that the class $\cS_0$ coincides with the classical Schur
class.

\medskip

The class $\cS_\kappa$ can be characterized alternatively (and
sometimes this characterization is taken as the definition of the class)
as the set of functions $w$ meromorphic on $\D$ and such that the kernel
$K_w(z,\zeta)$ defined in \eqref{0.1}
has $\kappa$ negative squares on $\D\cap \rho(w)$ ($\rho(w)$ stands for
the domain of analyticity of $w$); in formulas: ${\rm sq}_-(K_w)=\kappa$.
The last equality means that for every choice of an integer $n$ and of
$n$ points $z_1, \ldots,z_n\in\D\cap \rho(w)$, the Hermitian matrix
\eqref{0.6} has at most $\kappa$ negative eigenvalues:
\begin{equation}
{\rm sq}_-\left[\frac{1-\overline{w(z_i)}w(z_j)}
{1-\bar{z}_iz_j}\right]_{i,j=1}^n\le \kappa,
\label{0.7}
\end{equation}
and for at least one such choice it has exactly $\kappa$ negative
eigenvalues counted with multiplicities. In what follows, we will
say  ``$w$ has $\kappa$ negative squares'' rather than ``the kernel
$K_w$ has $\kappa$ negative squares''.

\medskip

Due to representation \eqref{0.6} and in view of the quite simple
structure of finite Blaschke products, most of the results concerning the
boundary behavior of generalized Schur functions can be derived
from the corresponding classical results for the Schur class functions.
For example, the nontangential boundary limit $d_w(t)$ (defined in
\eqref{0.2b}) exists for every $t\in\T$ and satisfies $d_w(t)>-\infty$
(not necessarily nonnegative, in contrast to the definite case). Indeed,
if $w$ is of the form \eqref{0.6}, then
\begin{equation}
\frac{1-|w(z)|^2}{1-|z|^2}=
\frac{1}{|B(z)|^2} \, \left(\frac{1-|S(z)|^2}{1-|z|^2}-
\frac{1-|B(z)|^2}{1-|z|^2}\right).
\label{0.7a}
\end{equation}
Passing to the limits as $z$ tends to $t\in\T$ in the latter equality and
taking into account that  $|B(t)|=1$, we get
$$
d_w(t)=d_S(t)-d_B(t)>-\infty,
$$
since $d_w(t_0)\ge 0$ and $d_B(t)<\infty$. Furthermore, as in
the definite case, if  $d_w(t)<\infty$, then the nontangential limit
\eqref{0.2a} exists and is unimodular.

\medskip

Now we formulate indefinite analogues of Problems \ref{P:0.1} and
\ref{P:0.2}. The data set for these problems will consist
of $n$ points $t_1,\ldots,t_n$ on $\T$, $n$ unimodular numbers
$w_1,\ldots,w_n$ and $n$ real numbers $\gamma_1,\ldots,\gamma_n$:
\begin{equation}
t_i\in\T,\quad |w_i|=1,\quad \gamma_i\in\R\qquad (i=1,\ldots,n).
\label{0.8}
\end{equation}
As in the definite case, we associate to the interpolation data
\eqref{0.8} the Pick matrix $P$ via the formula \eqref{0.5}
which is still Hermitian (since $\gamma_j\in\R$), but not positive
semidefinite, in general. Let $\kappa$ be the number of
its negative eigenvalues:
\begin{equation}
\kappa:={\rm sq}_-P,
\label{0.9}
\end{equation}
where
\begin{equation}
P=[P_{ij}]_{i,j=1}^n\quad \mbox{and}\quad
P_{ij}=\left\{\begin{array}{ccc} {\displaystyle
\frac{1-\overline{w}_iw_j}{1-\overline{t}_it_j}}&\mbox{for} & i\neq j,\\
\gamma_j&\mbox{for}& i=j.\end{array}\right.
\label{0.10}
\end{equation}
The next problem is an indefinite analogue of Problem \ref{P:0.2}
and it coincides with Problem \ref{P:0.2} if $\kappa=0$.
\begin{Pb}
Given the data set \eqref{0.8}, find all functions
$w\in\cS_\kappa$ (with $\kappa$ defined in \eqref{0.9}) such that
\begin{equation}
d_w(t_i):=\lim_{z\rightarrow t_i} \,
\frac{1-|w(z)|^2}{1-|z|^2}=\gamma_i\qquad (i=1,\ldots,n)
\label{0.11}
\end{equation}
and
\begin{equation}
w(t_i):=\lim_{z\rightarrow t_i} \, w(z)=w_i\qquad(i=1,\ldots,n).
\label{0.12}
\end{equation}
\label{P:0.4}
\end{Pb}
The analogue of Problem \ref{P:0.1} is:
\begin{Pb}
Given the data set \eqref{0.8}, find all functions
$w\in\cS_\kappa$ (with $\kappa$ defined in \eqref{0.9}) such that
\begin{equation}
d_w(t_i)\le \gamma_i\quad\mbox{and}\quad w(t_i)=w_i\qquad (i=1,\ldots,n).
\label{0.13}
\end{equation}
\label{P:0.5}
\end{Pb}
Interpolation conditions for the two above problems are clear:
existence of the nontangential limits $d_w(t_i)$'s implies existence of
the nontangential limits $w(t_i)$'s;  upon prescribing the values of
these limits (or upon prescribing upper bounds for $d_w(t_i)$'s) we come
up with interpolation conditions \eqref{0.11}--\eqref{0.13}.
The choice \eqref{0.9} for the index of $\cS_\kappa$ should be explained
in some more detail.
\begin{Rk}
If a generalized Schur function $w$ satisfies interpolation conditions
\eqref{0.13}, then it has at least $\kappa={\rm sq}_-P$ negative squares.
\label{R:0.6}
\end{Rk}
Indeed, if $w$ is a generalized Schur function of the class
$\cS_{\widetilde{\kappa}}$ and $t_1,\ldots,t_n$ are distinct points on
$\T$ such that
$$
d_w(t_i)<\infty\quad\mbox{for}\quad i=1,\ldots,n,
$$
then the nontangential boundary limits $w(t_i)$'s exist (and are
unimodular)
and one can pass to the limit in \eqref{0.7} (as $t_i\to z_i$ for
$i=1,\ldots,n$) to conclude that the Hermitian matrix
\begin{equation}
P^w(t_1,\ldots,t_n)=\left[P^w_{ij}\right]_{i,j=1}^n\quad\mbox{with}\quad
P^w_{ij}=\left\{\begin{array}{ccc} {\displaystyle
\frac{1-\overline{w(t_i)}w(t_j)}{1-\overline{t}_it_j}}&\mbox{for} & i\neq
j \\ d_w(t_i)&\mbox{for}& i=j\end{array}\right.
\label{0.14}
\end{equation}
satisfies
\begin{equation}
{\rm sq}_-P^w(t_1,\ldots,t_n)\le \widetilde{\kappa}.
\label{0.15}
\end{equation}
If $w$ meets conditions \eqref{0.12}, then the nondiagonal entries in
the matrices $P^w(t_1,\ldots,t_n)$ and $P$ coincide which clearly
follows from the definitions \eqref{0.10} and \eqref{0.14}. It follows
from
the same definitions that
$$
P-P^w(t_1,\ldots,t_n)=\left[\begin{array}{ccc}\gamma_1-d_w(t_1) &&0 \\ &
\ddots & \\ 0 && \gamma_n-d_w(t_n)\end{array}\right]
$$
and thus, conditions \eqref{0.11} and the first series of conditions
in \eqref{0.13} can be written equivalently in the matrix form as
\begin{equation}
P^w(t_1,\ldots,t_n)=P\quad\mbox{and}\quad P^w(t_1,\ldots,t_n)\le P,
\label{0.15a}
\end{equation}
respectively. Each one of the two last relations implies, in view of
\eqref{0.15} that
$$
{\rm sq}_-P\le \widetilde{\kappa}.
$$
Thus, the latter condition is necessary for existence of
a function $w$ of the class $\cS_{\widetilde{\kappa}}$ satisfying
interpolation conditions \eqref{0.13} (or \eqref{0.11} and \eqref{0.12}).
The choice \eqref{0.9} means that we are concerned about generalized
Schur functions with the minimally possible negative index.

\medskip

Problems \ref{P:0.4} and \ref{P:0.5} are indefinite analogues of Problems
\ref{P:0.2} and \ref{P:0.1}, respectively. Now we introduce another
boundary interpolation problem that does not appear in the context
of classical Schur functions.
\begin{Pb}
Given the data set \eqref{0.8}, find all functions
$w\in\cS_{\kappa^\prime}$ for some $\kappa^\prime\le \kappa={\rm sq}_-P$
such that conditions \eqref{0.13} are satisfied at all but
$\kappa-\kappa^\prime$ points $t_1,\ldots,t_n$.
\label{P:0.6}
\end{Pb}
In other words, a solution $w$ to the last problem is allowed
to have less then $\kappa$ negative squares and to omit some of
interpolation conditions (but not too many of them). The significance of
Problem \ref{P:0.6} will be explained in the next section.

\medskip

\section{Main results}
\setcounter{equation}{0}

The purpose of the paper is to obtain parametrizations of solution sets
$\; {\mathbb S}_{13}$, ${\mathbb S}_{14}$ and ${\mathbb S}_{16}$ for
Problems \ref{P:0.4}, \ref{P:0.5} and \ref{P:0.6}, respectively.
First we note that
\begin{equation}
{\mathbb S}_{13}\subseteq {\mathbb S}_{14}\subseteq {\mathbb S}_{16}
\quad\mbox{and}\quad
{\mathbb S}_{14}={\mathbb S}_{16}\cap \cS_\kappa.
\label{0.16}
\end{equation}
Inclusions in \eqref{0.16} are self-evident. If $w$ is a solution
of Problems \ref{P:0.6} with $\kappa^\prime=\kappa$, then
$\kappa-\kappa^\prime=0$ which means that conditions \eqref{0.13} are
satisfied at all points $t_1,\ldots,t_n$ and thus, $w\in{\mathbb
S}_{14}$. Thus,  ${\mathbb S}_{14}\subseteq{\mathbb S}_{16}\cap
\cS_\kappa$. The reverse inclusion is evident, since
${\mathbb S}_{14}\subseteq\cS_\kappa$.
Note also that if $\kappa=0$, then Problems \ref{P:0.5}
and \ref{P:0.6} are equivalent: ${\mathbb S}_{14}={\mathbb S}_{16}$.

\medskip

It turns out that in the indefinite setting (i.e., when $\kappa>0$),
Problem \ref{P:0.6} plays the same role as Problem \ref{P:0.5} does in
the classical setting: it always has a solution and, in the indetereminate
case, the solution set ${\mathbb S}_{16}$ admits a linear fractional
parametrization with the free Schur class parameter. The case when $P$ is
singular, is relatively simple:
\begin{Tm}
Let $P$ be singular. Then Problem \ref{P:0.6} has a unique solution
$w$ which is the ratio of two finite Blaschke products
$$
w(z)=\frac{B_1(z)}{B_2(z)}
$$
with no common zeroes and such that
$$
{\rm deg} \, B_1+{\rm deg} \, B_2={\rm rank} \, P.
$$
Furthermore, if ${\rm deg} \, B_2=\kappa$, then $w$ is also a solution
of Problem \ref{P:0.5}.
\label{T:0.8}
\end{Tm}
The proof will be given in Section 7. Now we turn to a more interesting
case when $P$ is not singular. In this case, we pick an arbitrary point
$\mu\in\T\setminus\{t_1,\ldots,t_n\}$
and introduce the $2\times 2$ matrix valued function
\begin{eqnarray}
\Theta(z)&=&\left[\begin{array}{cc}\Theta_{11}(z) & \Theta_{12}(z)\\
\Theta_{21}(z)&\Theta_{22}(z)\end{array}\right]\label{0.17} \\
&=&I_{2}+(z-\mu)\left[\begin{array}{c}C \\ E\end{array}\right]
(zI_n-T)^{-1}P^{-1}(I_n-\mu T^*)^{-1}\left[\begin{array}{cc}
C^* & -E^*\end{array}\right]\nonumber
\end{eqnarray}
where
\begin{equation}
T=\left[\begin{array}{ccc}t_1 &&
\\ &\ddots & \\ && t_n\end{array}\right],\quad
E=\begin{bmatrix} 1 & \ldots & 1\end{bmatrix},\quad
C=\begin{bmatrix} w_1 & \ldots & w_n\end{bmatrix}.
\label{0.18}
\end{equation}
Note that the Pick matrix $P$ defined in \eqref{0.10} satisfies the
following
identity
\begin{equation}
 P-T^*PT=E^*E-C^*C.
\label{0.19}
\end{equation}
Indeed, equality of nondiagonal entries in (\ref{0.19}) follows from
the definition \eqref{0.14} of $P$, whereas diagonal entries in both sides
of \eqref{0.19} are zeroes. Identity \eqref{0.19} and all its ingredients
will play an important role in the subsequent analysis.

\medskip

The function $\Theta$ defined in \eqref{0.17} is rational and has
simple poles at $t_1,\ldots,t_n$. Note some extra properties of $\Theta$.
Let $J$ be a signature matrix defined as
\begin{equation}
J=\left[\begin{array}{cr} 1 & 0 \\ 0 & -1\end{array}\right].
\label{0.20}
\end{equation}
It turns out that $\Theta$ is $J$-unitary on the unit circle, i.e., that
\begin{equation}
\Theta(t)J\Theta(t)^*=J\quad\mbox{for every} \; \; t\in\T\cap \rho(\Theta)
\label{0.21}
\end{equation}
and the kernel
\begin{equation}
K_{\Theta,J}(z,\zeta):=\frac{J-\Theta(z)J\Theta(\zeta)^*}{1-z\bar{\zeta}}
\label{0.21a}
\end{equation}
has $\kappa={\rm sq}_-P$ negative squares on $\D$:
\begin{equation}
{\rm sq}_-K_{\Theta,J}=\kappa.
\label{0.22}
\end{equation}
We shall use the symbol ${\mathcal W}_\kappa$ for the class of $2\times 2$
meromorphic functions satisfying conditions \eqref{0.21} and \eqref{0.22}.
It is well known that for {\em every} function  $\Theta\in{\mathcal
W}_\kappa$, the linear fractional transformation
\begin{equation}
{\bf T}_\Theta: \; \cE \longrightarrow\frac{\Theta_{11}\cE+\Theta_{12}}
{\Theta_{21}\cE+\Theta_{22}}
\label{0.23}
\end{equation}
is well defined for every Schur class function $\cE$ and maps $\cS_0$ into
$\bigcup_{\kappa^\prime\le\kappa}\cS_{\kappa^\prime}$. This map is not
onto and the question about its range is of certain interest. If $\Theta$
is of the form \eqref{0.17}, the range of the transformation \eqref{0.23}
is ${\mathbb S}_{16}$:
\begin{Tm}
Let $P$, $T$, $E$ and $C$ be defined as in \eqref{0.10} and \eqref{0.18}
and let $w$ be a function meromorphic on $\D$. If $P$ is invertible,
then $w$ is a solution of Problem \ref{P:0.6} if and only if it is
of the form
\begin{equation}
w(z)={\bf T}_\Theta[\cE](z):=\frac{\Theta_{11}(z)\cE(z)+\Theta_{12}(z)}
{\Theta_{21}(z)\cE(z)+\Theta_{22}(z)},
\label{0.24}
\end{equation}
for some Schur function $\cE\in\cS_0$.
\label{T:0.9}
\end{Tm}
It is not difficult to show that every rational function $\Theta$ from
the class ${\mathcal W}_\kappa$ with simple poles at $t_1,\ldots,t_n\in\T$
and normalized to $I_2$ at $\mu\in\T$, is necessarily of the form
\eqref{0.17} for some row vector $C\in\C^{1\times n}$ with unimodular
entries, with $E$ as in \eqref{0.18} and with a Hermitian invertible
matrix $P$ having $\kappa$ negative squares and being subject to
the Stein identity \eqref{0.19}. Thus, Theorem \ref{T:0.9}
clarifies the interpolation meaning of the range of a linear
fractional transformation based on a rational function $\Theta$ of the
class ${\mathcal W}_\kappa$ with simple poles on the boundary of the unit
disk.

\medskip

The necessity part in Theorem \ref{T:0.9} will be obtained in Section 3
using an appropriate adaptation of the V. P. Potapov's method of
the Fundamental Matrix Inequality (FMI) to the context of generalized
Schur functions. The proof of the sufficiency part rests on Theorems
\ref{T:0.20} and  \ref{T:0.21} which are of certain independent
interest. To formulate these theorems, let us introduce
the numbers $\tc_1,\ldots,\tc_n$ and $\te_1,\ldots,\te_n$ by
\begin{equation}
\tc_i^*:=-\lim_{z\to t_i}(z-t_i)\Theta_{21}(z)\quad\mbox{and}
\quad \te_i^*:=\lim_{z\to
t_i}(z-t_i)\Theta_{22}(z)\quad(i=1,\ldots,n)
\label{0.25}
\end{equation}
(for notational convenience we will write sometimes $a^*$ rather than
$\overline{a}$ for $a\in\C$). It turns out $|\tc_i|=|\te_i|\neq 0$ (see
Lemma \ref{L:2.1} below for the proof) and therefore the following
numbers
\begin{equation}
\q_i:=\frac{\tc_i}{\te_i}=\frac{\te_i^*}{\tc_i^*}=
-\lim_{z\to t_i}\frac{\Theta_{22}(z)}{\Theta_{21}(z)} \quad
(i=1,\ldots,n)
\label{0.25a}
\end{equation}
are unimodular:
\begin{equation}
|\q_i|=1\quad (i=1,\ldots,n).
\label{0.25b}
\end{equation}
Furthermore let $\widetilde{p}_{ii}$ stand for the $i$-th diagonal entry
of the matrix  $P^{-1}$, the inverse of the Pick matrix. It is
self-evident that for a
fixed  $i$, any function $\cE\in\cS_0$ satisfies exactly one of the
following six conditions:
\begin{eqnarray}
&&{\bf C}_1: \quad\mbox{The function $\cE$ fails to have a nontangential
boundary limit $\q_i$ at $t_i$.}\nonumber\\
&& {\bf C}_2:\quad
\cE(t_i):=\lim_{z\to
t_i}\cE(z)=\q_i\quad\mbox{and}\quad
d_{\cE}(t_i):=\frac{1-|\cE(z)|^2}{1-|z|^2}=\infty.\label{0.26}\\
&& {\bf C}_3:\quad \cE(t_i)=\q_i\quad\mbox{and}\quad
-\frac{\widetilde{p}_{ii}}{|\te_i|^2}<d_{\cE}(t_i)<\infty.
\label{0.27}  \\
&& {\bf C}_4:\quad \cE(t_i)=\q_i\quad\mbox{and}\quad
0\le d_{\cE}(t_i)<-\frac{\widetilde{p}_{ii}}{|\te_i|^2}.
\label{0.28}\\
&& {\bf C}_5:\quad\cE(t_i)=\q_i\quad\mbox{and}\quad
d_{\cE}(t_i)=-\frac{\widetilde{p}_{ii}}{|\te_i|^2}>0.
\label{0.29}\\
&& {\bf C}_6:\quad  \cE(t_i)=\q_i\quad\mbox{and}\quad
d_{\cE}(t_i)=\widetilde{p}_{ii}=0.
\label{0.30}
\end{eqnarray}
Note that condition ${\bf C}_1$ means that either the nontangential
boundary limit $\cE(t_i):={\displaystyle\lim_{z\to t_i}\cE(z)}$ fails to
exist or it exists and is not equal to $\q_i$.
Let us denote by ${\bf C}_{4-6}$ the disjunction of conditions
${\bf C}_4$, ${\bf C}_5$ and ${\bf C}_6$:
\begin{equation}
{\bf C}_{4-6}: \quad \cE(t_i)=\q_i\quad\mbox{and}\quad
d_{\cE}(t_i)\le -\frac{\widetilde{p}_{ii}}{|\te_i|^2}.
\label{0.31}
\end{equation}
The next theorem gives a classification of interpolation
conditions that are or are not satisfied by a function
$w$ of the form \eqref{0.24} in terms of the corresponding parameter
$\cE$.
\begin{Tm}
Let the Pick matrix $P$ be invertible, let $\cE$ be a Schur class
function, let $\Theta$ be given by \eqref{0.17}, let $w={\bf
T}_\Theta[\cE]$ and let $t_i$ be an interpolation node.
\begin{enumerate}
\item The nontangential boundary limits $d_w(t_i)$ and $w(t_i)$ exist
and are subject to
$$
d_w(t_i)=\gamma_i\quad\mbox{and}\quad
w(t_i)=w_i
$$
if and only if the parameter $\cE$ meets either condition
${\bf C}_1$ or ${\bf C}_2$.
\item The nontangential boundary limits $d_w(t_i)$ and $w(t_i)$ exist
and are subject to
$$
d_w(t_i)<\gamma_i\quad\mbox{and}\quad w(t_i)=w_i
$$
if and only if the parameter $\cE$ meets condition ${\bf C}_3$.
\item The nontangential boundary limits $d_w(t_i)$ and $w(t_i)$ exist
and are subject to
$$
\gamma_i<d_w(t_i)<\infty\quad\mbox{and}\quad w(t_i)=w_i.
$$
if and only if the parameter $\cE$ meets condition ${\bf C}_4$.
\item If  $\cE$ meets ${\bf C}_5$, then
$w$ is subject to one of the following:
\begin{enumerate}
\item[(a)] The limit $w(t_i)$ fails to exist.
\item[(b)] The limit $w(t_i)$ exists and $w(t_i)\neq w_i$.
\item[(c)] $w(t_i)=w_i$ and $d_w(t_i)=\infty$.
\end{enumerate}
\item If  $\cE$ meets ${\bf C}_6$, then $w$ is the ratio
of two finite Blaschke products,
$$
d_w(t_i)<\infty\quad\mbox{and}\quad w(t_i)\neq w_i.
$$
\end{enumerate}
\label{T:0.20}
\end{Tm}
We note an immediate consequence of the last theorem.
\begin{Cy}
A function $w={\bf T}_\Theta[\cE]$ meets the $i$-th interpolation
conditions for Problem \ref{P:0.5}:
$$
d_w(t_i)\le \gamma_i\quad\mbox{and}\quad w(t_i)=w_i
$$
if and only if the corresponding parameter $\cE\in\cS_0$ meets
the condition ${\bf C}_{1-3}:={\bf C}_{1}\vee {\bf
C}_{2}\vee{\bf C}_{3}$ at $t_i$.
\label{C:0.22}
\end{Cy}
Note that Problem \ref{P:0.4} was considered in \cite{bgr} for rational
generalized Schur functions. It was shown (\cite[Theorem 21.1.2]{bgr})
that all rational solutions
of Problem \ref{P:0.4} are parametrized by the formula \eqref{0.24} when
$\cE$ varies over the set of all rational Schur functions such that
(in the current terminology)
$$
\cE(t_i)\neq \q_i\quad\mbox{for} \; \; i=1,\ldots,n.
$$
Note that if $\cE$ is a {\em rational} Schur function admitting a
unimodular value $\cE(t_0)$ at a boundary point $t_0\in\T$, then the
limit $d_w(t_0)$ always exists and equals
$t_0\cE^\prime(t_0)\cE(t_0)^*$. The latter follows from the converse
Carath\'eodory-Julia theorem (see e.g., \cite{sarasonsubh, shapiro}):
\begin{eqnarray*}
d_w(t_0):=\lim_{z\to t_0}\frac{1-|\cE(z)|^2}{1-|z|^2}&=&
\lim_{z\to t_0}\frac{1-\cE(z)\cE(t_0)^*}{1-z\bar{t}_0}\\
&=&\lim_{z\to t_0}\frac{\cE(t_0)-\cE(z)}{t_0-z}\cdot
\frac{\cE(t_0)^*}{\bar{t}_0}\\
&=&t_0\cE^\prime(t_0)\cE(t_0)^*<\infty.
\end{eqnarray*}
Thus, a Schur function $\cE$ cannot satisfy condition ${\bf C}_2$
at a boundary point $t_i$ therefore, Statement (1) in Theorem \ref{T:0.20}
recovers Theorem 21.1.2 in \cite{bgr}. The same conclusion can be done
when $\cE$ is not rational but still analytic at $t_i$.
In the case when $\cE$ is not rational and admits the nontangential
boundary limit $\cE(t_i)=\q_i$, the situation is more subtle: Statement
(1) shows that even in this case (if the convergence of $\cE(z)$ to
$\cE(t_i)$ is not too fast), the function $w={\bf T}[\cE]$ may
satisfy interpolation conditions \eqref{0.11}, \eqref{0.12}.

\medskip

The next theorem concerns the number of negative squares of the
function $w={\bf T}_\Theta[\cE]$.
\begin{Tm}
If the Pick matrix $P$ is invertible and has $\kappa$ negative
eigenvalues, then a Schur function $\cE\in\cS_0$ may satisfy
conditions ${\bf C}_{4-6}$ at at most $\kappa$ interpolation nodes.
Furthermore, if $\cE$ meets conditions ${\bf C}_{4-6}$ at exactly
$\ell$ ($\le \kappa$) interpolation nodes, then the function $w={\bf
T}_\Theta[\cE]$ belongs to the class $\cS_{\kappa-\ell}$.
\label{T:0.21}
\end{Tm}
Corollary \ref{C:0.22} and Theorem \ref{T:0.21} imply the sufficiency part
in Theorem \ref{T:0.9}. Indeed, any Schur function $\cE$ satisfies either
conditions ${\bf C}_{4-6}$ or ${\bf C}_{1-3}$ at every interpolation
node $t_i$  ($i=1,\ldots,n$). Let $\cE$ meet conditions ${\bf C}_{4-6}$
at $t_{i_1},\ldots, t_{i_{\ell}}$ and ${\bf C}_{1-3}$ at
other $n-\ell$ interpolation nodes $t_{j_1},\ldots,
t_{j_{n-\ell}}$. Then, by Corollary \ref{C:0.22}, the function
$w={\bf T}_\Theta[\cE]$ satisfies interpolation conditions \eqref{0.13}
for $i\in\{j_1,\ldots, j_{n-\ell}\}$ and fails to satisfy
at least one of these conditions at the remaining $\ell$
interpolation nodes. On the other hand, $w$ has exactly
$\kappa-\ell$ negative squares, by Theorem \ref{T:0.21}.
Thus, for every $\cE\in\cS_0$, the function $w={\bf T}_\Theta[\cE]$ solves
Problem \ref{P:0.6}.

\medskip

Note also that Theorems \ref{T:0.9} and \ref{T:0.21} lead to
parametrizations of solution sets for Problems \ref{P:0.4} and
\ref{P:0.5}. Indeed,  by inclusions \eqref{0.16}, every solution
$w$ to Problem \ref{P:0.4} (or to Problem \ref{P:0.5}) is also of the form
\eqref{0.24} for some $\cE\in\cS_0$.  Thus, there is a
chance to describe the solution sets ${\mathbb S}_{13}$ and ${\mathbb
S}_{14}$ by appropriate selections of the parameter $\cE$ in
\eqref{0.24}. Theorem \ref{T:0.21} indicates how these selections have to
be made.
\begin{Tm}
A function $w$ of the form \eqref{0.24} is a solution to Problem
\ref{P:0.4} if and only if the corresponding parameter $\cE\in\cS_0$
satisfies either condition ${\bf C}_1$ or  ${\bf C}_2$ for every
$i\in\{1,\ldots,n\}$.
\label{T:0.10}
\end{Tm}
\begin{Tm}
A function $w$ of the form \eqref{0.24} is a solution to Problem
\ref{P:0.5} if and only if the corresponding parameter $\cE\in\cS_0$
either fails to have a nontangential boundary limit $\q_i$ at $t_i$ or
$$
\cE(t_i)=\q_i\quad\mbox{and}\quad
d_{\cE}(t_i)>-\frac{\widetilde{p}_{ii}}{|\te_i|^2}
$$
for every $i=1,\ldots,n$ (in other words, $\cE$ meets one of conditions
${\bf C}_1$, ${\bf C}_2$, ${\bf C}_3$ at each interpolation node $t_i$).
\label{T:0.11}
\end{Tm}
As a consequence of Theorems \ref{T:0.9} and \ref{T:0.11} we
get curious necessary and sufficient conditions (in terms of the
interpolation data \eqref{0.8}) for Problems \ref{P:0.5} and \ref{P:0.6}
to be equivalent (that is, to have the same solution sets).
\begin{Cy}
Problems \ref{P:0.5} and \ref{P:0.6} are equivalent if and only if
all the diagonal entries of the inverse $P^{-1}$ of the Pick matrix are
positive.
\label{C:0.12}
\end{Cy}
Indeed, in this case, all the conditions in Theorem \ref{T:0.11}
are fulfilled for every $\cE\in\cS_0$ and every $i\in\{1,\ldots,n\}$ and
formula \eqref{0.21} gives a free Schur class parameter
description of all solutions $w$ of Problem \ref{P:0.5}.

\bigskip

In the course of the proof of Theorem \ref{T:0.21} we will discuss the
following related question: given indices
$i_1,\ldots,i_\ell\in\{1,\ldots,n\}$,
does  there exist a parameter $\cE\in\cS_0$ satisfying conditions
${\bf C}_{4-6}$ at $t_{i_1},\ldots,t_{i_\ell}$? Due to Theorems
\ref{T:0.9} and \ref{T:0.20}, this question can be posed equivalently:
does there exist a solution $w$ to Problem \ref{P:0.6} that misses
interpolation conditions at $t_{i_1},\ldots,t_{i_\ell}$ (Theorem
\ref{T:0.21} claims that if such a function exists, it belongs to the
class $\cS_{\kappa-\ell}$). The question admits a simple answer in terms
of a certain submatrix of
$P^{-1}=\left[\widetilde{p}_{ij}\right]_{i,j=1}^n$,
the inverse of the Pick matrix.
\begin{Tm}
There exists a parameter $\cE$ satisfying conditions
${\bf C}_{4-6}$ at $t_{i_1},\ldots,t_{i_\ell}$ if and only if the
$\ell\times \ell$ matrix
$$
{\mathcal P}:=\left[\widetilde{p}_{i_\alpha, i_\beta}
\right]_{\alpha,\beta=1}^\ell
$$
is negative semidefinite. Moreover, if ${\mathcal P}$ is
negative definite, then there are infinitely many such
parameters. If ${\mathcal P}$ is negative semidefinite
(singular), then there is only one such parameter, which is a
Blaschke product of degree $r={\rm rank} \, {\mathcal P}$.
\label{T:0.15}
\end{Tm}
Note that all the results announced above have their counterparts in the
context of the regular Nevanlinna-Pick problem with all the interpolation
nodes inside the unit disk \cite{bol}

\medskip

The paper is organized as follows: Section 3 contains some
needed auxiliary results which can be found (probably in a different
form) in many sources  and are included for the sake of completeness.
In Section 4 we prove the necessity part in  Theorem \ref{T:0.9} (see
Remark \ref{R:2.4}).  In Section 5 we prove Theorem \ref{T:0.20}. In
Section 6 we present the proofs of Theorems \ref{T:0.15} and \ref{T:0.21}
and complete the proof of Theorem \ref{T:0.9} (see Remark \ref{R:6:2}).
The proof of Theorem \ref{T:0.8} is contained in Section 7; some
illustrative numerical examples are presented in Section 8.

\medskip

\section{Some preliminaries}
\setcounter{equation}{0}

In this section we present some auxiliary results needed in the sequel.
We have already mentioned the Stein identity
\begin{equation}
 P-T^*PT=E^*E-C^*C
\label{1.1u}
\end{equation}
satisfied  by the Pick matrix $P$ constructed in \eqref{0.10} from the
interpolation data. Most of the facts recalled in this section rely on
this identity rather than on the special form \eqref{0.18} of matrices
$T$, $E$ and $C$.
\begin{La}
Let $T$, $E$ and $C$ be defined as in \eqref{0.18}, let $P$ defined in
\eqref{0.10} be invertible and let $\mu$ be a point on
$\T\setminus\{t_1,\ldots,t_n\}$. Then
\begin{enumerate}
\item The row vectors
\begin{equation}
\tE=\begin{bmatrix}\te_1 & \ldots &
\te_n\end{bmatrix}\quad\mbox{and}\quad
\tC=\begin{bmatrix}\tc_1 & \ldots & \tc_n\end{bmatrix}
\label{2.1}
\end{equation}
defined by
\begin{equation}
\left[\begin{array}{c}\tC \\ \tE\end{array}\right]
=\left[\begin{array}{c}C \\ E\end{array}\right]
(\mu I-T)^{-1}P^{-1}(I-\mu T^*)
\label{2.2}
\end{equation}
satisfy the Stein identity
\begin{equation}
P^{-1}-TP^{-1}T^*=\tE^*\tE-\tC^*\tC.
\label{2.3}
\end{equation}
\item The numbers $\tc_i$ and $\te_i$ are subject to
\begin{equation}
|\te_i|=|\tc_i|\neq 0 \quad\mbox{for} \; \; i=1,\ldots,n.
\label{2.4}
\end{equation}
\item The nondiagonal entries $\widetilde{p}_{ij}$ of $P^{-1}$
are given by
\begin{equation}
\widetilde{p}_{ij}=\frac{\te_i^*\te_j-\tc_i^*\tc_j}{1-t_i\bar{t}_j}
\qquad (i\neq j).
\label{2.4a}
\end{equation}
\end{enumerate}
\label{L:2.1}
\end{La}
{\bf Proof:} Under the assumption that $P$ is invertible, identity
(\ref{2.3}) turns out to be equivalent to (\ref{1.1u}).
Indeed, by (\ref{2.2}) and (\ref{1.1u}),
\begin{eqnarray*}
&&\tE^*\tE-\tC^*\tC\\
&&=(I-\bar{\mu}T)P^{-1}(\bar{\mu}I-T^*)^{-1}
\left[E^*E-C^*C\right] (\mu I-T)^{-1}P^{-1}(I-\mu T^*)\\
&&=(I-\bar{\mu}T)P^{-1}(\bar{\mu}I-T^*)^{-1}
\left[P-T^*PT\right] (\mu I-T)^{-1}P^{-1}(I-\mu T^*)\\
&&=(I-\bar{\mu}T)P^{-1}\left[(I-\mu
T^*)^{-1}P+PT(\mu I-T)^{-1}\right]P^{-1}(I-\mu T^*)\\
&&=(I-\bar{\mu}T)P^{-1}+\bar{\mu}TP^{-1}(I-\mu T^*)\\
&&=P^{-1}-TP^{-1}T^*.
\end{eqnarray*}
Let $P^{-1}=\left[\widetilde{p}_{ij}\right]_{i,j=1}^n$.
Due to \eqref{2.1} and \eqref{0.18}, equality of the $ij$-th entries
in \eqref{2.3} can be displayed as
\begin{equation}
\widetilde{p}_{ij}-t_i\bar{t}_j\widetilde{p}_{ij}=
\te_i^*\te_j-\tc_i^*\tc_j
\label{2.4k}
\end{equation}
and implies \eqref{2.4a} if $i\neq j$. Letting $i=j$ in \eqref{2.4k}
and taking into account that $|t_i|=1$, we get $|\te_i|=|\tc_i|$ for
$i=1,\ldots,n$. It remains to show that $\te_i$ and $\tc_i$ do not vanish.
To this end  let us assume that
\begin{equation}
\te_i=\tc_i=0.
\label{2.5}
\end{equation}
Let ${\bf e}_i$ be the $i$-th column of the identity matrix $I_n$.
Multiplying \eqref{2.3} by ${\bf e}_i$ on the right we get
$$
P^{-1}{\bf e}_i-TP^{-1}T^*{\bf e}_i=\tE^*\te_i-\tC^*\tc_i=0
$$
or equivalently, since $T^*{\bf e}_i=\bar{t}_i{\bf e}_i$,
$$
(I-\bar{t}_iT)P^{-1}{\bf e}_i=0.
$$
Since the points $t_1,\ldots,t_n$ are distinct, all the
diagonal entries but the $i$-th in the diagonal matrix $I-\bar{t}_iT$
are not zeroes; therefore, it follows from the last equality that
all the entries in the vector $P^{-1}{\bf e}_i$ but the $i$-th entry are
zeroes. Thus,
\begin{equation}
P^{-1}{\bf e}_i=\alpha {\bf e}_i
\label{2.6}
\end{equation}
for some $\alpha\in\C$ and, since $P$ is not singular, it follows that
$\alpha\neq 0$. Now we compare the $i$-th columns in the equality
\eqref{2.2} (i.e., we multiply both parts in \eqref{2.2} by  ${\bf e}_i$
on the right). For the left hand side we have, due to assumption
\eqref{2.5},
$$
\left[\begin{array}{c}\tC \\ \tE\end{array}\right]{\bf
e}_i=\left[\begin{array}{c}\tc_i \\ \te_i\end{array}\right]=
\left[\begin{array}{c}0 \\ 0\end{array}\right].
$$
For the right hand side, we have, due to \eqref{2.6} and
\eqref{0.18},
$$
\left[\begin{array}{c}C \\ E\end{array}\right]
(\mu I-T)^{-1}P^{-1}(I-\mu T^*){\bf e}_i=
\alpha\frac{1-\mu\overline{t}_i}{\mu-t_i}\left[\begin{array}{c}C \\
E\end{array}\right]{\bf e}_i=-\alpha \overline{t}_i\left[\begin{array}{c}
w_i \\ 1 \end{array}\right].
$$
By \eqref{2.2}, the right hand side expressions in the two last
equalities must be the same, which is not the case. The obtained
contradiction  completes the proof of \eqref{2.4}.\qed
\begin{Rk}
The numbers $\te_i$ and $\tc_i$ introduced in \eqref{2.1}, \eqref{2.2}
coincide with those in \eqref{0.25}.
\label{R:2.1}
\end{Rk}
For the proof we first note that the formula \eqref{0.17} for $\Theta$
can be written, on account of \eqref{2.2}, as
\begin{equation}
\Theta(z)=I_{2}+(z-\mu)\left[\begin{array}{c}C \\ E\end{array}\right]
(zI_n-T)^{-1}(\mu I_n-T)^{-1}\left[\begin{array}{cc}
\tC^* & -\tE^*\end{array}\right]
\label{2.6p}
\end{equation}
and then, since
$$
\lim_{z\to t_i}(z-t_i)(zI-T)^{-1}={\bf e}_i{\bf e}_i^*\quad
\mbox{and}\quad
{\bf e}_i^*(\mu I-T)^{-1}=(\mu -t_i)^{-1}{\bf e}_i^*
$$
(recall that ${\bf e}_i$ is the $i$-th column of the identity matrix
$I_n$), we have
\begin{eqnarray}
\lim_{z\to t_i}(z-t_i)\Theta(z)&=&
\lim_{z\to t_i}(z-\mu)\left[\begin{array}{c}C \\
E\end{array}\right]{\bf e}_i{\bf e}_i^*(\mu I-T)^{-1}
\left[\begin{array}{cc}\tC^* & -\tE^*\end{array}\right]\nonumber\\
&=&-\left[\begin{array}{c}C \\
E\end{array}\right]{\bf e}_i{\bf e}_i^*
\left[\begin{array}{cc}\tC^* & -\tE^*\end{array}\right]\nonumber \\
&=&-\left[\begin{array}{c}w_i \\ 1\end{array}\right]
\left[\begin{array}{cc}\tc_i^* & -\te_i^*\end{array}\right].\label{uh}
\end{eqnarray}
Comparing the bottom entries in the latter equality we get
\eqref{0.25}.\qed

\medskip

In the rest of the section we recall some needed results
concerning the function $\Theta$ introduced in \eqref{0.17}. These
results are well known in a more general situation when $T$, $C$ and
$E$ are matrices such that the pair $(\begin{bmatrix} C\\ E\end{bmatrix},
\; T)$ is observable:
\begin{equation}
\bigcap_{j\ge 0}{\rm Ker} \, \begin{bmatrix} C\\ E\end{bmatrix}T^j=\{0\},
\label{1.0}
\end{equation}
and $P$ is an invertible Hermitian matrix satisfying the Stein identity
\eqref{1.1u} (see e.g., \cite{bgr}). Note that the matrices defined in
\eqref{0.18} satisfy a stronger condition:
\begin{equation}
\bigcap_{j\ge 0} {\rm Ker} \, CT^j=\bigcap_{j\ge 0}{\rm Ker} \,
ET^j=\{0\}.
\label{1.0a}
\end{equation}
\begin{Rk}
Under the above assumptions, the function $\Theta$ defined via formula
\eqref{0.17} belongs to the class ${\mathcal W}_\kappa$ with
$\kappa={\rm sq}_-P$.
\label{R:2.0}
\end{Rk}
{\bf Proof:} The desired membership follows from the formula
\begin{equation}
K_{\Theta,J}(z,\zeta)=\left[\begin{array}{c}C \\ E\end{array}\right]
(zI-T)^{-1}P^{-1}(\bar{\zeta}I-T^*)^{-1}\left[\begin{array}{cc}
C^* & E^*\end{array}\right]
\label{1.25}
\end{equation}
for the kernel $K_{\Theta,J}$  defined in \eqref{0.21a}. The calculation
is straightforward and relies on the Stein identity \eqref{1.1u} only
(see e.g., \cite{bgr}). It follows from \eqref{1.25} that $\Theta$ is $J$
-unitary on $\T$ (that is, satisfies condition \eqref{0.21}) and that
$$
{\rm sq}_-K_{\Theta,J}\le {\rm sq}_-P=\kappa.
$$
Condition \eqref{1.0} guarantees that in fact ${\rm sq}_-K_{\Theta,J}
=\kappa$ (see \cite{bgr}).\qed
\begin{Rk}
{\rm Since $\Theta$ is $J$-unitary on $\T$ it holds, by the symmetry
principle, that $\Theta(z)^{-1}=J\Theta(1/\bar{z})^*J$,  which together
with formula \eqref{0.17} leads us to
\begin{equation}
\Theta(z)^{-1}=I_2-(z-\mu)\left[\begin{array}{c}C \\
E\end{array}\right](\mu I-T)^{-1}P^{-1}
(I-zT^*)^{-1}\begin{bmatrix}C^* & - E^*\end{bmatrix}.
\label{2.9a}
\end{equation}}
\label{R:2.2}
\end{Rk}
Besides \eqref{1.25} we will need realization formulas for two
related kernels. Verification of these fomulas \eqref{1.299} and
\eqref{2.8} is also straightforward and is based on the Stein identities
\eqref{1.1u} and \eqref{2.3}, respectively.
\begin{Rk}
Let $\Theta$ be defined as in \eqref{0.17}. The following identities hold
for every choice of $z, \, \zeta\not\in\{t_1,\ldots,t_n\}$:
\begin{eqnarray}
\frac{\Theta(\zeta)^{-*}J\Theta(z)^{-1}-J}{1-z\bar{\zeta}}&=&
\left[\begin{array}{r}C \\ -E\end{array}\right]
(I-\bar{\zeta}T)^{-1}
P^{-1}(I-zT^*)^{-1}\begin{bmatrix}C^* & -E^*\end{bmatrix},\label{1.299}\\
\frac{J-\Theta(\zeta)^*J\Theta(z)}{1-z\bar{\zeta}}&=&
\left[\begin{array}{r}\tC \\ -\tE\end{array}\right]
(\bar{\zeta}I-T^*)^{-1}P(zI-T)^{-1}\left[\begin{array}{cc}
\tC^* & -\tE^*\end{array}\right].\label{2.8}
\end{eqnarray}
\label{R:2.3}
\end{Rk}
Let us consider conformal partitionings
\begin{equation}
P=\begin{bmatrix}P_{11} & P_{12} \\ P_{21} & P_{22}\end{bmatrix},\quad
P^{-1}=\begin{bmatrix}\tP_{11} & \tP_{12} \\ \tP_{21} &
\tP_{22}\end{bmatrix},\quad
T=\begin{bmatrix}T_1 & 0 \\ 0 & T_2\end{bmatrix},
\label{2.11a}
\end{equation}
\begin{equation}
E=\begin{bmatrix}E_1 & E_2\end{bmatrix},\quad
C=\begin{bmatrix}C_1 & C_2\end{bmatrix},\quad
\tE=\begin{bmatrix}\tE_1 & \tE_2\end{bmatrix},\quad
\tC=\begin{bmatrix}\tC_1 & \tC_2\end{bmatrix}
\label{2.11d}
\end{equation}
where $P_{22}, \, \tP_{22}, \, T_2\in\C^{\ell\times \ell}$ and $E_2, \,
C_2, \, \tE_2, \, \tC_2\in\C^{1\times \ell}$. Note that these
decompositions contain one restrictive assumption: it is assumed that
the matrix $T$ is block diagonal.
\begin{La}
Let us assume that $P_{11}$ is invertible and let
${\rm sq}_-P_{11}=\kappa_1\le \kappa$. Then $\tP_{22}$ is invertible,
${\rm sq}_-\tP_{22}=\kappa-\kappa_1$ and the functions
\begin{equation}
\Theta^{(1)}(z)=I_{2}+(z-\mu)\left[\begin{array}{c}C_1 \\
E_1\end{array}\right] (zI-T_1)^{-1}P_{11}^{-1}(I-\mu
T_1^*)^{-1}\left[\begin{array}{cc}C_1^* & -E_1^*\end{array}\right]
\label{2.11b}
\end{equation}
and
\begin{equation}
\widetilde{\Theta}^{(2)}(z)=
I_{2}+(z-\mu)\left[\begin{array}{c}\tC_2 \\ \tE_2\end{array}\right]
(I-\mu T_2^*)^{-1}\tP_{22}^{-1}(zI-T_2)^{-1}
\left[\begin{array}{cc}\tC_2^* & -\tE_2^*\end{array}\right]
\label{2.11c}
\end{equation}
belong to ${\mathcal W}_{\kappa_1}$ and ${\mathcal W}_{\kappa-\kappa_1}$,
respectively. Furthermore, the function $\Theta$ defined in \eqref{0.17}
admits a factorization
\begin{equation}
\Theta(z)=\Theta^{(1)}(z)\widetilde{\Theta}^{(2)}(z).
\label{2.12}
\end{equation}
\label{L:2.5}
\end{La}
{\bf Proof:} The first statement follows by standard Schur complement
arguments: since $P$ and $P_{11}$ are invertible, the matrix
$P_{22}-P_{21}P_{11}^{-1}P_{12}$ (the Schur complement of $P_{11}$ in
$P$) is invertible and has $\kappa-\kappa_1$ negative eigenvalues. Since
the block $\tP_{22}$ in $P^{-1}$ equals
$(P_{22}-P_{21}P_{11}^{-1}P_{12})^{-1}$,
it also has $\kappa-\kappa_1$ negative eigenvalues. Realization formulas
\begin{equation}
K_{\Theta^{(1)},J}(z,\zeta)=R(z)P_{11}^{-1}R(\zeta)^*\quad
\mbox{and}\quad
K_{\widetilde{\Theta}^{(2)},J}(z,\zeta)=\widetilde{R}(z)
\widetilde{P}_{22}\widetilde{R}(\zeta)^*,
\label{2.12k}
\end{equation}
where we have set for short
$$
R(z)=\left[\begin{array}{c}C_1 \\
E_1\end{array}\right] (zI-T_1)^{-1},\quad
\widetilde{R}(z)=\left[\begin{array}{c}\tC_2 \\
\tE_2\end{array}\right](I-\mu
T_2^*)^{-1}\widetilde{P}_{22}^{-1}(zI-T_2)^{-1},
$$
are established exactly as in Remark \ref{R:2.0} and rely on the Stein
identities
\begin{equation}
P_{11}-T_{1}^*P_{11}T_1=E_1^*E_1-C_1^*C_1\quad\mbox{and}\quad
\widetilde{P}_{22}^{-1}-T_2\widetilde{P}_{22}T_2^*=\tE_2^*\tE_2-\tC_2^*\tC_2
\label{2.12m}
\end{equation}
which hold true, being parts of identities \eqref{1.1u} and \eqref{2.3}.
Formulas \eqref{2.12k} guarantee that the rational functions
$\Theta^{(1)}$ and $\widetilde{\Theta}^{(2)}$ are $J$-unitary on $\T$ and
moreover, that
\begin{equation}
{\rm sq}_-K_{\Theta^{(1)},J}\le {\rm sq}_-
P_{11}=\kappa_1\quad\mbox{and}\quad
{\rm sq}_-K_{\widetilde{\Theta}^{(2)},J}\le {\rm
sq}_-\widetilde{P}_{22}=\kappa-\kappa_1.
\label{2.12v}
\end{equation}
Assuming that the factorization formula \eqref{2.12} is already
proved, we have
$$
K_{\Theta,J}(z,\zeta)
=K_{\Theta^{(1)},J}(z,\zeta)+\Theta^{(1)}(z)
K_{\widetilde{\Theta}^{(2)},J}(z,\zeta)\Theta^{(1)}(\zeta)^*
$$
and thus,
$$
\kappa={\rm sq}_-K_{\Theta,J}\le {\rm sq}_-K_{\Theta^{(1)},J}+{\rm
sq}_-K_{\widetilde{\Theta}^{(2)},J}
$$
which together with inequalities \eqref{2.12v} imply
$$
{\rm sq}_-K_{\Theta^{(1)},J}=\kappa_1\quad\mbox{and}\quad
{\rm sq}_-K_{\widetilde{\Theta}^{(2)},J}=\kappa-\kappa_1.
$$
It remains to prove \eqref{2.12}. Making use of the well known
equality
\begin{equation}
P^{-1}=\left[\begin{array}{cc}P_{11}^{-1} & 0 \\ 0 &
0\end{array}\right]+\left[\begin{array}{c}-P_{11}^{-1}P_{12} \\
1\end{array}\right]\widetilde{P}_{22}\left[\begin{array}{cc}-
P_{21}P_{11}^{-1} & 1\end{array}\right]
\label{2.15}
\end{equation}
we conclude from \eqref{2.2} that
\begin{eqnarray}
\left[\begin{array}{c}\tC_2 \\ \tE_2\end{array}\right]
&=&\left[\begin{array}{c}C \\ E\end{array}\right]
(\mu I_n-T)^{-1}P^{-1}(I_n-\mu T^*)\begin{bmatrix}
0 \\ I_\ell\end{bmatrix}\nonumber\\
&=&\left[\begin{array}{c}C \\ E\end{array}\right]
(\mu I_n-T)^{-1}\left[\begin{array}{c}-P_{11}^{-1}P_{12} \\
1\end{array}\right]\tP_{22}(I_\ell-\mu T_2^*).
\label{2.15a}
\end{eqnarray}
This last relation allows us to rewrite \eqref{2.11c} as
\begin{equation}
\widetilde{\Theta}^{(2)}(z)=
I_{2}+(z-\mu)\left[\begin{array}{c}C \\ E\end{array}\right]
(\mu I-T)^{-1}\left[\begin{array}{c}-P_{11}^{-1}P_{12} \\
1\end{array}\right](zI-T_2)^{-1}
\left[\begin{array}{cc}\tC_2^* & -\tE_2^*\end{array}\right].
\label{2.15c}
\end{equation}
Now we substitute \eqref{2.15} into the formula \eqref{0.17} defining
$\Theta$ and take into account \eqref{2.11b} and \eqref{2.15a} to get
\begin{eqnarray*}
\Theta(z)&=&\Theta^{(1)}(z)+(z-\mu)\left[\begin{array}{c}C \\
E\end{array}\right] (zI_n-T)^{-1}
\left[\begin{array}{c}-P_{11}^{-1}P_{12} \\
1\end{array}\right]\widetilde{P}_{22}\\
&&\qquad\times \left[\begin{array}{cc}-
P_{21}P_{11}^{-1} & 1\end{array}\right]
(I_n-\mu T^*)^{-1}\left[\begin{array}{cc}
C^* & -E^*\end{array}\right]\\
&=&\Theta^{(1)}(z)+(z-\mu)\left[\begin{array}{c}C \\
E\end{array}\right] (zI_n-T)^{-1}
\left[\begin{array}{c}-P_{11}^{-1}P_{12} \\
1\end{array}\right]\\
&& \qquad \times(\mu I-T_2)^{-1}
\left[\begin{array}{cc}\tC_2^* & -\tE_2^*\end{array}\right].
\end{eqnarray*}
Thus, \eqref{2.12} is equivalent to
\begin{eqnarray*}
\widetilde{\Theta}^{(2)}(z)&=&
I_2+(z-\mu)\Theta^{(1)}(z)^{-1}\left[\begin{array}{c}C \\
E\end{array}\right] (zI_n-T)^{-1}
\left[\begin{array}{c}-P_{11}^{-1}P_{12} \\
1\end{array}\right]\\
&&\qquad\times (\mu I -T_2)^{-1}
\left[\begin{array}{cc}\tC_2^* & -\tE_2^*\end{array}\right].
\end{eqnarray*}
Comparing the last relation with \eqref{2.15c} we conclude that to
complete the proof it suffices to show that
\begin{eqnarray}
&&\Theta^{(1)}(z)^{-1}\left[\begin{array}{c}C \\
E\end{array}\right] (zI_n-T)^{-1}
\left[\begin{array}{c}-P_{11}^{-1}P_{12} \\
1\end{array}\right](\mu I -T_2)^{-1}\nonumber\\
&&=\left[\begin{array}{c}C \\ E\end{array}\right]
(\mu I-T)^{-1}\left[\begin{array}{c}-P_{11}^{-1}P_{12} \\
1\end{array}\right](zI-T_2)^{-1}.
\label{2.16}
\end{eqnarray}
The explicit formula for $\Theta^{(1)}(z)^{-1}$ can be obtained similarly
to \eqref{2.9a}:
\begin{equation}
\Theta^{(1)}(z)^{-1}=I_2-(z-\mu)\left[\begin{array}{c}C_1 \\
E_1\end{array}\right](\mu I-T_1)^{-1}P_{11}^{-1}
(I-zT_1^*)^{-1}\begin{bmatrix}C_1^* & - E_1^*\end{bmatrix}.
\label{2.13}
\end{equation}
Next, comparing the top block entries in the Stein identity
\eqref{1.1u} we get, due to decompositions \eqref{2.11a} and
\eqref{2.11d},
$$
\begin{bmatrix}P_{11} & P_{12}\end{bmatrix}
-T_1^*\begin{bmatrix}P_{11} & P_{12}\end{bmatrix}T=E_1^*E-C_1^*C
$$
which, being multiplied by $(I-zT_1^*)^{-1}$ on the left and by
$(zI-T)^{-1}$ on the right, leads us to
\begin{eqnarray}
&&(I-zT_1^*)^{-1}\left(E_1^*E-C_1^*C\right)(zI-T)^{-1}\nonumber\\
&&=(I-zT_1^*)^{-1}T_1^*\begin{bmatrix}P_{11} & P_{12}\end{bmatrix}
+\begin{bmatrix}P_{11} & P_{12}\end{bmatrix}(zI-T)^{-1}.
\label{2.17}
\end{eqnarray}
Upon making use of \eqref{2.16} and \eqref{2.17} we have
\begin{eqnarray*}
&&\Theta^{(1)}(z)^{-1}\left[\begin{array}{c}C \\
E\end{array}\right] (zI_n-T)^{-1}
\left[\begin{array}{c}-P_{11}^{-1}P_{12} \\
1\end{array}\right]\nonumber\\
&&=\left[\begin{array}{c}C \\ E\end{array}\right] (zI-T)^{-1}
\left[\begin{array}{c}-P_{11}^{-1}P_{12} \\
1\end{array}\right]\\
&&\quad +(z-\mu)\left[\begin{array}{c}C_1 \\
E_1\end{array}\right](\mu I-T_1)^{-1}
\begin{bmatrix}I & P_{11}^{-1}P_{12}\end{bmatrix}(zI-T)^{-1}
\left[\begin{array}{c}-P_{11}^{-1}P_{12} \\
1\end{array}\right]\nonumber\\
&&=-\left[\begin{array}{c}C_1 \\ E_1\end{array}\right]
(zI-T_1)^{-1}P_{11}^{-1}P_{12}+\left[\begin{array}{c}C_2 \\
E_2\end{array}\right](zI-T_2)^{-1}\\
&&\quad +(z-\mu)\left[\begin{array}{c}C_1 \\
E_1\end{array}\right](\mu I-T_1)^{-1}\left(
P_{11}^{-1}P_{12}(zI-T_2)^{-1}
-(zI-T_1)^{-1}P_{11}^{-1}P_{12}\right)\\
&&=-\left[\begin{array}{c}C_1 \\ E_1\end{array}\right]
(\mu I-T_1)^{-1}P_{11}^{-1}P_{12}(\mu I-T_2)+\left[\begin{array}{c}C_2 \\
E_2\end{array}\right](zI-T_2)^{-1}\\
&&=\left[\begin{array}{c}C \\ E\end{array}\right]
(\mu I-T)^{-1}\left[\begin{array}{c}-P_{11}^{-1}P_{12} \\
1\end{array}\right](zI-T_2)^{-1}
\end{eqnarray*}
which proves \eqref{2.16} and therefore, completes the proof of the
lemma.\qed
\begin{Rk}
{\rm The case when $\ell=1$ in Lemma \ref{L:2.5} will be of special
interest.
In this case,
$$P_{22}=\gamma_n,\quad \tP_{22}=\widetilde{p}_{nn}, \quad
T_2=t_n,\quad C_2=w_n,\quad E_2=1,\quad\tC_2=\tc_n,\quad\tE_2=\te_n.
$$
Then the formula \eqref{2.11c} for $\widetilde{\Theta}^{(2)}$ simplifies
to
\begin{equation}
\widetilde{\Theta}^{(2)}(z)=
I_{2}+\frac{z-\mu}{(1-\mu\bar{t}_n)(z-t_n)}\left[\begin{array}{c}\tc_n \\
\te_n\end{array}\right]\widetilde{p}_{nn}^{-1}
\left[\begin{array}{cc}\tc_n^* & -\te_n^*\end{array}\right].
\label{2.20}
\end{equation}}
\label{R:2.6}
\end{Rk}

\section{Fundamental Matrix Inequality}
\setcounter{equation}{0}

In this section we characterize the solution set ${\mathbb S}_{16}$ of
Problem \ref{P:0.6} in terms of certain Hermitian kernel. We start
with some  simple observations.
\begin{Pn}
Let $K(z,\zeta)$ be a Hermitian kernel defined on $\Omega\subseteq\C$
and with ${\rm sq}_-K=\kappa$. Then
\begin{enumerate}
\item For every choice of an integer $p$, of a Hermitian $p\times
p$ matrix $A$ and of a $p\times 1$ vector valued function $B$,
$$
{\rm sq}_-\begin{bmatrix} A & B(z) \\ B(\zeta)^* & K(z,\zeta)\end{bmatrix}
\le \kappa +p.
$$
\item If $\lambda_1,\ldots,\lambda_p$ are points in $\Omega$ and if
\begin{equation}
A=\left[K(\lambda_j,\lambda_i)\right]_{i,j=1}^p\quad\mbox{and}\quad
B(z)=\begin{bmatrix} K(z,\lambda_1) \\ \vdots \\
K(z,\lambda_p)\end{bmatrix},
\label{3.1a}
\end{equation}
then
\begin{equation}
{\rm sq}_-\begin{bmatrix} A & B(z) \\ B(\zeta)^* & K(z,\zeta)\end{bmatrix}
=\kappa.
\label{3.1b}
\end{equation}
\end{enumerate}
\label{P:2.1}
\end{Pn}
{\bf Proof:} For the proof of the first
statement we have to show that for every integer $m$ and every choice of
points $z_1,\ldots,z_m\in\Omega$, the block matrix
\begin{equation}
M=\left[\begin{bmatrix} A & B(z_j) \\
B(z_i)^* & K(z_j,z_i)\end{bmatrix}\right]_{i,j=1}^m
\label{3.1c}
\end{equation}
has at most $\kappa+p$ negative eigenvalues. It is easily seen that
$M$ contains $m$ block identical rows of the form
$$
\begin{bmatrix} A & B(z_1) & A & B(z_2) & \ldots & A &
B(z_n)\end{bmatrix}.
$$
Deleting all these rows but one and deleting also the corresponding
columns, we come up with the $(m+p)\times(m+p)$ matrix
$$
\widetilde{M}=\begin{bmatrix} A & B(z_1)& \ldots & B(z_m)\\
B(z_1)^* & K(z_1,z_1) & \ldots & K(z_1,z_m)\\
\vdots & \vdots &&\vdots \\
B(z_m)^* & K(z_m,z_1) & \ldots & K(z_m,z_m)\end{bmatrix}
$$
having the same number of positive and negative eigenvalues as $M$.
The bottom $m\times m$ principal submatrix of $\widetilde{M}$ has at most
$\kappa$ negative eigenvalues since ${\rm sq}_-K=\kappa$. Since
$\widetilde{M}$ is Hermitian, we have by the Cauchy's interlacing theorem
(see e.g., \cite[p. 59]{Bhatia}), that ${\rm sq}_-\widetilde{M}\le
\kappa+p$. Thus, ${\rm sq}_-M\le \kappa+p$ which completes the proof of
Statement 1.

\smallskip

If $A$ and $B$ are of the form \eqref{3.1a}, then the matrix
$M$ in \eqref{3.1c} is of the form
$\left[K(\zeta_j,\zeta_i)\right]_{i,j=1}^{m+pm}$
where all the points $\zeta_i$ live in $\Omega$. Since ${\rm
sq}_-K=\kappa$, it follows that ${\rm sq}_-M\le \kappa$ for every choice
of $z_1,\ldots,z_m$ in $\Omega$ which means that the kernel
$\begin{bmatrix} A & B(z) \\ B(\zeta)^* & K(z,\zeta)\end{bmatrix}$
has at most $\kappa$ negative squares on $\Omega$. But it has at least
$\kappa$ negative squares since it contains the kernel $K(z,\zeta)$ as
a principal block. Thus, \eqref{3.1b} follows.\qed
\begin{Tm}
Let $P$, $T$, $E$ and $C$ be defined as in \eqref{0.10} and \eqref{0.18},
let $w$ be a function meromorphic on $\D$ and let the kernel $K_w$
be defined as in \eqref{0.1}. Then
$w$ is a solution of Problem \ref{P:0.6} if and only if the kernel
\begin{equation}
{\bf K}_w(z,\zeta):=\begin{bmatrix}P & (I-zT^*)^{-1}(E^*-C^*w(z)) \\
(E-w(\zeta)^*C)(I-\bar{\zeta}T)^{-1} & K_w(z,\zeta)\end{bmatrix}
\label{2.1u}
\end{equation}
has $\kappa$ negative squares on $\D\cap\rho(w)$:
\begin{equation}
{\rm sq}_-{\bf K}_w(z,\zeta)=\kappa.
\label{2.2u}
\end{equation}
\label{T:2.2}
\end{Tm}
{\bf Proof of the necessity part:} Let $w$ be a solution of Problem
\ref{P:0.6}, i.e., let $w$ belong to the class $\cS_{\kappa'}$ for some
$\kappa'\le \kappa$ and satisfy conditions \eqref{0.13} at all but
$\kappa-\kappa'$ interpolation nodes.

\smallskip

First we consider the case when $w\in\cS_\kappa$. Then $w$ satisfies all
the conditions \eqref{0.13} (i.e., $w$ is also a solution to Problem
\ref{P:0.5}). Furthermore, ${\rm sq}_-K_w=\kappa$ and by the second
statement in Proposition \ref{P:2.1}, the kernel
\begin{equation}
{\bf K}^{(1)}(z,\zeta):=\begin{bmatrix}
K_w(z_{1},z_{1})& \ldots & K_w(z_n,z_1) &
K_w(z,z_{1})\\ \vdots & & \vdots&\vdots \\
K_w(z_1,z_n) &  \ldots & K_w(z_n,z_n) & K_w(z,z_n)\\
K_w(z_1,\zeta) & \ldots & K_w(z_n,\zeta) & K_w(z,\zeta)\end{bmatrix}
\label{2.3u}
\end{equation}
has $\kappa$ negative squares on $\D\cap\rho(w)$ for every
choice of points $z_{1},\ldots,z_n\in \D\cap\rho(w)$.
Since the limits $d_w(t_i)$ and $w(t_i)=w_i$ exist for
$i=1,\ldots,n$, it follows that
\begin{equation}
\left[K_w(z_j,z_i)\right]_{i,j=1}^n
=\left[\frac{1-w(z_i)^*w(z_j)}{1-\bar{z}_iz_j}
\right]_{i,j=1}^n\longrightarrow P^w(t_1,\ldots,t_n)
\label{2.4u}
\end{equation}
(by definition \eqref{0.14} of the matrix $P^w(t_1,\ldots,t_n)$)
and also
$$
K_w(z_i,\zeta)=\frac{1-w(\zeta)^*w(z_i)}{1-\bar{\zeta}z_i}
\longrightarrow
\frac{1-w(\zeta)^*w_i}{1-\bar{\zeta}t_i}\qquad (i=1,\ldots,n).
$$
Note that by the structure \eqref{0.18} of the matrices $T$, $E$ and $C$,
$$
(E-w(\zeta)^*C)(I-\bar{\zeta}T)^{-1}=
\begin{bmatrix}{\displaystyle\frac{1-w(\zeta)^*w_1}{1-\bar{\zeta}t_1}}
& \ldots & {\displaystyle\frac{1-w(\zeta)^*w_n}
{1-\bar{\zeta}t_n}}\end{bmatrix}
$$
which, being combined with the previous relation, gives
\begin{equation}
\begin{bmatrix}K_w(z_1,\zeta) & \ldots & K_w(z_n,\zeta)\end{bmatrix}
\longrightarrow (E-w(\zeta)^*C)(I-\bar{\zeta}T)^{-1}.
\label{2.5u}
\end{equation}
Now we take the limit in \eqref{2.3u} as $z_i\to t_i$ for $i=1,\ldots,n$;
on account of \eqref{2.4u} and \eqref{2.5u}, the limit kernel has the form
$$
{\bf K}^{(2)}(z,\zeta):=\begin{bmatrix}P^w(t_1,\ldots,t_n) &
(I-zT^*)^{-1}(E^*-C^*w(z)) \\
(E-w(\zeta)^*C)(I-\bar{\zeta}T)^{-1} & K_w(z,\zeta)\end{bmatrix}.
$$
Since ${\bf K}^{(2)}$ is the limit of a family of kernels each of which
has $\kappa$ negative squares, ${\rm sq}_-{\bf K}^{(2)}\le \kappa$.
It remains to note that the kernel ${\bf K}_w$ defined in \eqref{2.1u}
is expressed in terms of ${\bf K}^{(2)}$ as
$$
{\bf K}_w(z,\zeta)={\bf
K}^{(2)}(z,\zeta)+\begin{bmatrix}P-P^w(t_1,\ldots,t_n) & 0 \\ 0 & 0
\end{bmatrix}
$$
and since the second term on the right hand side is positive semidefinite
(due to the first series of conditions in \eqref{0.13}; see also
\eqref{0.15a}),
$$
{\rm sq}_-{\bf K}_w\le {\rm sq}_-{\bf K}^{(2)}\le \kappa.
$$
On the other hand, since ${\bf K}_w$ contains the kernel
$K_w$ as a principal submatrix, ${\rm sq}_-{\bf K}_w\ge {\rm
sq}_-K_w=\kappa$ which eventually leads us to \eqref{2.2u}.
Note that in this part of the proof we have not used the fact that
${\rm sq}_-P=\kappa$.

\smallskip

Now we turn to the general case: let $w\in\cS_{\kappa^\prime}$ for some
$\kappa^\prime\le \kappa$ and let conditions \eqref{0.13} be fulfilled at
all but $\ell:=\kappa-\kappa^\prime$ interpolation nodes $t_i$'s.
We may assume without loss of generality that conditions \eqref{0.13} are
satisfied at $t_i$ for $i=1,\ldots,n-\ell$:
\begin{equation}
d_w(t_i)\le \gamma_i\quad\mbox{and}\quad w(t_i)=w_i\qquad
(i=1,\ldots,n-\ell).
\label{2.6u}
\end{equation}
Let us consider conformal partitionings \eqref{2.11a}, \eqref{2.11d} for
matrices $P$, $T$, $C$ and $E$ and let us set for short
\begin{equation}
F_i(z)=(I-zT_i^*)^{-1}\left(E_i^*-C_i^*w(z)\right) \qquad (i=1,2)
\label{2.6m}
\end{equation}
so that
\begin{equation}
\begin{bmatrix} F_1(z)\\
F_2(z)\end{bmatrix}=(I-zT^*)^{-1}\left(E^*-C^*w(z)\right).
\label{2.6t}
\end{equation}
The matrix $P_{11}$ is the Pick matrix of the truncated interpolation
problem with the data $t_i, \, w_i, \,  \gamma_i$ ($i=1,\ldots,n-\ell$)
and with interpolation conditions \eqref{2.6u}. By the first part of the
proof, the kernel
\begin{equation}
\widetilde{\bf K}_w(z,\zeta):=\begin{bmatrix}P_{11} & F_1(z)\\
F_1(\zeta)^* & K_w(z,\zeta)\end{bmatrix}
\label{2.7k}
\end{equation}
has $\kappa^\prime$ negative squares on $\D\cap\rho(w)$.
Now we apply the first statement in Proposition \ref{P:2.1} to
\begin{equation}
K(z,\zeta)=\widetilde{\bf K}_w(z,\zeta), \quad
B(z)=\begin{bmatrix}P_{21} & F_2(z)\end{bmatrix}\quad\mbox{and}\quad
A=P_{22}
\label{2.7u}
\end{equation}
to conclude that
\begin{equation}
{\rm sq}_-\begin{bmatrix} P_{22} & B(z) \\ B(\zeta)^* &
\widetilde{\bf K}_w(z,\zeta)\end{bmatrix}\le {\rm sq}_-
\widetilde{\bf K}_w+\ell=
\kappa^\prime+(\kappa-\kappa^\prime)=\kappa.
\label{2.8u}
\end{equation}
By \eqref{2.7u} and \eqref{2.7k}, the latter kernel equals
$$
\begin{bmatrix} P_{22} & B(z) \\ B(\zeta)^* &
\widetilde{\bf K}_w(z,\zeta)\end{bmatrix}=
\begin{bmatrix} P_{22} & P_{21} & F_2(z) \\ P_{12} & P_{11} & F_1(z)\\
F_2(\zeta)^* & F_1(\zeta)^* & K_w(z,\zeta)\end{bmatrix}.
$$
Now it follows from \eqref{2.1u} and \eqref{2.7k} that
$$
{\bf K}_w(z,\zeta)=U
\begin{bmatrix} P_{22} & B(z) \\ B(\zeta)^* &
\widetilde{\bf K}_w(z,\zeta)\end{bmatrix}U^*,\quad\mbox{where}\quad
U=\left[\begin{array}{llc}0 & I_{n-\ell} & 0 \\ I_{\ell} & 0 & 0 \\
0 & 0 & 1\end{array}\right]
$$
which, on account of \eqref{2.8u}, implies that
${\rm sq}_-{\bf K}_w\le\kappa$. Finally, since ${\bf K}_w$
contains $P$ as a principal submatrix, ${\rm sq}_-{\bf K}_w\ge {\rm
sq}_-P=\kappa$ which now implies \eqref{2.2u} and completes the proof of
the necessity part of the theorem.
The proof of the sufficiency part will be given in Sections 6 and 7
(see Remarks \ref{R:6.3} and \ref{R:7.2} there).\qed

\medskip

In the case when $P$ is invertible,
all the functions satisfying \eqref{2.2u} can be described
in terms of a linear fractional transformation.
\begin{Tm}
Let the Pick matrix $P$ be invertible and let $\Theta=
[\Theta_{ij}]$ be the $2\times 2$ matrix valued function defined in
\eqref{0.17}. A function $w$ meromorphic on $\D$ is subject to FMI
\eqref{2.2u} if and only if it is of the form
\begin{equation}
w(z)={\bf T}_\Theta[\cE]:=\frac{\Theta_{11}(z)\cE(z)+\Theta_{12}(z)}
{\Theta_{21}(z)\cE(z)+\Theta_{22}(z)}
\label{2.9u}
\end{equation}
for some Schur function $\cE\in\cS_0$.
\label{T:2.3}
\end{Tm}
{\bf Proof:}  The proof is about the same as in the definite case.
Let ${\bf S}$ be the Schur complement of $P$ in the kernel ${\bf K}_w$
defined in \eqref{2.1u}:
$$
{\bf S}(z,\zeta):=K_w(z,\zeta)-
(E-w(\zeta)^*C)(I-\bar{\zeta}T)^{-1}P^{-1}(I-zT^*)^{-1}(E^*-C^*w(z)).
$$
Obvious equalities
$$
K_w(z,\zeta):=\frac{1-w(\zeta)^*w(z)}{1-\bar{\zeta}z}=
-\begin{bmatrix}w(\zeta)^*&
1\end{bmatrix}J\begin{bmatrix}w(z) \\ 1\end{bmatrix}
$$
where $J$ is the matrix introduced in \eqref{0.20}, and
$$
E-w(\zeta)^*C=-\begin{bmatrix}w(\zeta)^*&
1\end{bmatrix}J\left[\begin{array}{r}C \\ -E\end{array}\right]
$$
allows us to represent ${\bf S}$ in the form
\begin{eqnarray*}
{\bf S}(z,\zeta)&=&-\begin{bmatrix}w(\zeta)^*&
1\end{bmatrix}
\left\{\frac{J}{1-z\bar{\zeta}}+\left[\begin{array}{r}C \\ -E
\end{array}\right](I-\bar{\zeta}T)^{-1}P^{-1}\right.\\
&& \quad \times\left.(I-zT^*)^{-1}\left[\begin{array}{cc}
C^* & -E^*\end{array}\right]
\right\}\begin{bmatrix}w(z) \\ 1\end{bmatrix}
\end{eqnarray*}
or, on account of identity \eqref{1.299}, as
$$
{\bf S}(z,\zeta)=-\begin{bmatrix}w(\zeta)^*&
1\end{bmatrix}\frac{\Theta(\zeta)^{-*}J\Theta(z)^{-1}}{1-z\bar{\zeta}}
\begin{bmatrix}w(z) \\ 1\end{bmatrix}.
$$
By the standard Schur complement argument,
$$
{\rm sq}_-{\bf K}_w={\rm sq}_-P+{\rm sq}_-{\bf S}
$$
which implies, since ${\rm sq}_-P=\kappa$, that \eqref{2.2u} holds if and
only if the kernel ${\bf S}$ is positive definite on $\rho(w)\cap\D$:
\begin{equation}
-\begin{bmatrix}w(\zeta)^*&
1\end{bmatrix}\frac{\Theta(\zeta)^{-*}J\Theta(z)^{-1}}{1-z\bar{\zeta}}
\begin{bmatrix}w(z) \\ 1\end{bmatrix}\succeq 0.
\label{2.10u}
\end{equation}
It remains to show that \eqref{2.10u} holds if and only if $w$ is of the
form \eqref{2.9u}. To show the ``only if'' part, let us consider
meromorphic functions $u$ and $v$ defined by
\begin{equation}
\begin{bmatrix}u(z) \\ v(z)\end{bmatrix}:=\Theta(z)^{-1}
\begin{bmatrix}w(z) \\ 1\end{bmatrix}.
\label{2.11u}
\end{equation}
Then inequality \eqref{2.10u} can be written in terms of these functions as
\begin{equation}
-\begin{bmatrix}u(\zeta)^*& v(\zeta)^*\end{bmatrix}
\frac{J}{1-\bar{\zeta}z}\begin{bmatrix}u(z)
\\ v(z)\end{bmatrix}=\frac{v(\zeta)^*v(z)-
u(\zeta)^*u(z)}{1-\bar{\zeta}z}\succeq 0.
\label{2.12u}
\end{equation}
As it follows from definition \eqref{2.11u}, $u$ and $v$ are analytic
on $\rho(w)\cap\D$. Moreover,
\begin{equation}
v(z)\neq 0 \quad\mbox{for every} \; \;  z\in\rho(w)\cap\D.
\label{2.13u}
\end{equation}
Indeed, assuming that $v(\xi)=0$ at some point $\xi\in\D$,
we conclude from \eqref{2.12u} that
$u(\xi)= 0$ and then \eqref{2.11u} implies that $\det \,
\Theta(\xi)^{-1}=0$ which is a contradiction. Due to \eqref{2.13u},
we can introduce the meromorphic function
\begin{equation}
\cE(z)=\frac{u(z)}{v(z)}
\label{2.14u}
\end{equation}
which is analytic on $\rho(w)\cap\D$. Writing \eqref{2.12u} in terms of
$\cE$ as
$$
v(\zeta)^* \cdot \frac{1-\cE(\zeta)^*\cE(z)}{1-\bar{\zeta}z}\cdot v(z)
\succeq 0 \quad (z,\zeta\in\rho(w)\cap\D),
$$
we then take advantage of \eqref{2.13u} to conclude that
$$
\frac{1-\cE(\zeta)^*\cE(z)}
{1-\bar{\zeta}z}\succeq 0 \quad (z,\zeta\in\rho(w)\cap\D).
$$
The latter means that $\cE$ is (after an analytic continuation to the all
of $\D$) a Schur function. Finally, it follows from \eqref{2.11u} that
$$
\begin{bmatrix}w \\ 1\end{bmatrix}=\Theta
\begin{bmatrix}u \\ v\end{bmatrix}=
\begin{bmatrix}\Theta_{11}u+\Theta_{12}v\\
\Theta_{21}u+\Theta_{22}v\end{bmatrix}
$$
which in turn implies
$$
w=\frac{\Theta_{11}u+\Theta_{12}v}{\Theta_{21}u+\Theta_{22}v}
=\frac{\Theta_{11}\cE+\Theta_{12}}{\Theta_{21}\cE+\Theta_{22}}={\bf
T}_{\Theta}[\cE].
$$
Now let $\cE$ be a Schur function. Then the function
$$
V(z)=\Theta_{21}(z)\cE(z)+\Theta_{22}(z)
$$
does not vanish identically. Indeed, since  $\Theta$ is rational and
$\Theta(\mu)=I_2$, it follows that $\Theta_{22}(z)\approx 1$ and
$\Theta_{21}(z)\approx 0$ if $z$ is close enough to $\mu$. Since
$|\cE(z)|\le 1$ everywhere in $\D$, the function $V$ does not vanish
on ${\mathcal U}_\delta=\left\{z\in\D: \; \; |z-\mu|<\delta\right\}$
if $\delta$ is small enough. Thus,
formula \eqref{2.9u} makes sense and can be written equivalently as
$$
\begin{bmatrix}w(z) \\ 1\end{bmatrix}=\Theta(z)\begin{bmatrix}\cE(z) \\
1\end{bmatrix}\cdot\frac{1}{V(z)}
$$
Then it is readily seen that
\begin{eqnarray*}
\frac{1-\cE(\zeta)^*\cE(z)}{1-\bar{\zeta}z}&=&-
\begin{bmatrix}\cE(\zeta)^*&
1\end{bmatrix}\frac{J}{1-\bar{\zeta}z}\begin{bmatrix}\cE(z)
\\ 1\end{bmatrix}\\
&=&-\frac{1}{V(\zeta)^*V(z)}\cdot\begin{bmatrix}w(\zeta)^*&
1\end{bmatrix}\frac{\Theta(\zeta)^{-*}J\Theta(z)^{-1}}{1-z\bar{\zeta}}
\begin{bmatrix}w(z) \\ 1\end{bmatrix}
\end{eqnarray*}
for $z,\zeta\in\rho(w)\cap\D$. Since $\cE$ is a Schur function, the latter
kernel is positive on $\rho(w)\cap\D$ and since $V\not\equiv 0$,
\eqref{2.10u} follows.\qed
\begin{Rk}
Combining Theorems \ref{T:2.2} and \ref{T:2.3} we get the necessity
part in Theorem \ref{T:0.9}.
\label{R:2.4}
\end{Rk}
Indeed, by the necessity part in Theorem
\ref{T:2.2}, any solution $w$ of Problem \ref{P:0.6} satisfies
\eqref{2.2u}; then by Theorem \ref{T:2.3}, $w={\bf T}_\Theta[\cE]$ for
some $\cE\in\cS_0$.

\medskip

In the case when $\kappa=0$, Theorem \ref{T:2.2} was established in
\cite{kov}.
\begin{Tm}
Let the Pick matrix $P$ be positive semidefinite. Then a function $w$
defined on $\D$ is a solution to Problem \ref{P:0.1} (i.e., belongs to the
Schur class $\cS_0$ and meets conditions \eqref{0.4}) if and only if
\begin{equation}
{\bf K}_w(z,\zeta)\succeq 0 \quad (z,w\in\D)
\label{2.15u}
\end{equation}
where ${\bf K}_w(z,\zeta)$ is the kernel defined in \eqref{2.1u}.
\label{T:2.5}
\end{Tm}
Under the a priori assumption that $w$ is a Schur function, condition
\eqref{2.15u} can be replaced by a seemingly weaker matrix inequality
$$
{\bf K}_w(z,z)\ge 0 \quad \mbox{for every} \; \; z\in\D
$$
which is known in interpolation theory
as a Fundamental Matrix Inequality (FMI) of V. P. Potapov. We will follow
this terminology and will consider relation \eqref{2.2u} as an indefinite
analogue of V. P. Potapov's FMI. It is appropriate to note that a
variation of the  Potapov's method was first applied to the
Nevanlinna-Pick problem (with finitely many interpolation nodes inside the
unit disk) for generalized Schur functions in \cite{Gol}. We conclude this
section with another theorem concerning the classical case which will be
useful for the subsequent analysis.
\begin{Tm} (1) If the Pick matrix $P$ is positive definite then all the
solutions $w$ to Problem \ref{P:0.1} are parametrized by the formula
\eqref{0.24} with the coefficient matrix $\Theta$ defined as in
\eqref{0.17} with $\cE$ being a free Schur class parameter.

(2) If $P$ is positive semidefinite and singular, then Problem
\ref{P:0.1} has a unique solution $w$ which is a Blaschke product of
degree $r={\rm rank} \, P$. Furthermore, this unique solution
can be represented as
\begin{equation}
w(z)=\frac{x^*(I-zT_2^*)^{-1}E^*}{x^*(I-zT_2^*)^{-1}C^*}
\label{5.14p}
\end{equation}
where $T$, $C$ and $E$ are defined as in \eqref{0.18} and where
$x$ is any nonzero vector such that $Px=0$.
\label{T:2.6}
\end{Tm}
These results are well known and has been established using different
methods in \cite{b2, kov, ballhelton1, bgr, kh}. In regard to methods
used in the present paper, note that the first statement follows
immediately from Theorems \ref{T:2.5} and \ref{T:2.3}.
This demonstrates how the Potapov's method works in the definite case
(and this is exactly how the result was established in \cite{kov}).
The second statement also can be derived from Theorem \ref{T:2.5}:
if $w$ solves Problem \ref{P:0.1}, then the kernel ${\bf K}_w(z,\zeta)$
defined in \eqref{2.1u} is positive definite. Multiplying it by the
vector $\begin{bmatrix}x \\1 \end{bmatrix}$ on the right and by its
adjoint on the left we come to the positive definite kernel
$$
\begin{bmatrix}x^*Px & x^*(I-zT^*)^{-1}(E^*-C^*w(z)) \\
(E-w(\zeta)^*C)(I-\bar{\zeta}T)^{-1}x & K_w(z,\zeta)\end{bmatrix}\succeq
0.
$$
Thus, for every $x\neq 0$ such that $Px=0$, we also have
$$
x^*(I-zT^*)^{-1}(E^*-C^*w(z))\equiv 0.
$$
Solving the latter identity for $w$ we arrive at formula \eqref{5.14p}.
The numerator and the denominator in \eqref{5.14p} do not vanish
identically due to conditions \eqref{1.0a}. Since $x$ can be chosen
so that $n-{\rm rank} \, P-1$ its coordinates are zeros, the rational
function $w$ is of McMillan degree $r={\rm rank} \, P$. Due to the Stein
identity  \eqref{1.1u}, $w$ is inner and therefore, it is a finite
Blachke product of degree $r$.

\medskip

\section{Parameters and interpolation conditions}
\setcounter{equation}{0}

In this section we prove Theorem \ref{T:0.20}. It will be
done in several steps formulated as separate theorems. In what follows,
$U_{\cE}$ and $V_{\cE}$ will stand for the functions
\begin{equation}
U_{\cE}(z)=\Theta_{11}(z)\cE(z)+\Theta_{12}(z),\qquad
V_{\cE}(z)=\Theta_{21}(z)\cE(z)+
\Theta_{22}(z)
\label{4.1}
\end{equation}
for a fixed Schur function $\cE$, so that
\begin{equation}
\begin{bmatrix}U_{\cE}(z)\\ V_{\cE}(z)\end{bmatrix}=\Theta(z)
\begin{bmatrix}\cE(z)\\ 1\end{bmatrix}
\label{4.2a}
\end{equation}
and (\ref{0.24}) takes the form
\begin{equation}
w(z):={\bf T}_{\Theta}[\cE]=\frac{U_{\cE}(z)}{V_{\cE}(z)}.
\label{4.2}
\end{equation}
Substituting \eqref{2.6p} into \eqref{4.2a} and setting
\begin{equation}
\Psi(z)=(zI-T)^{-1}\left(\tE^*-\tC^*\cE(z)\right)
\label{4.3}
\end{equation}
for short, we get
\begin{eqnarray}
U_{\cE}(z)&=&\cE(z)-(z-\mu)C(\mu I-T)^{-1}\Psi(z),\label{4.4}\\
V_{\cE}(z)&=&1-(z-\mu)E(\mu I-T)^{-1}\Psi(z).
\label{4.5}
\end{eqnarray}
Furthermore, for $w$ of the form \eqref{4.2}, we have
\begin{equation}
\frac{1-w(\zeta)^*w(z)}{1-\bar{\zeta}z}=
\frac{1}{V_{\cE}(\zeta)^*V_{\cE}(z)}\cdot
\frac{V_{\cE}(\zeta)^*V_{\cE}(z)-
U_{\cE}(\zeta)^*U_{\cE}(z)}{1-\overline{\zeta}z}.
\label{4.6}
\end{equation}
Note that
\begin{eqnarray*}
V_{\cE}(\zeta)^*V_{\cE}(z)-U_{\cE}(\zeta)^*U_{\cE}(z)&=&
-\left[\begin{array}{cc} U_{\cE}(\zeta)^* &
V_{\cE}(\zeta)^*\end{array}\right]
J\left[\begin{array}{r} U_{\cE}(z) \\
V_{\cE}(z)\end{array}\right]\nonumber\\
&=&\left[\begin{array}{cc} \cE(\zeta)^* & 1
\end{array}\right]\Theta(\zeta)^{*}J\Theta(z)
\left[\begin{array}{c} \cE(z) \\ 1\end{array}\right]\nonumber\\
&=&1-\cE(\zeta)^*\cE(z)+(1-\bar{\zeta}z)\Psi(\zeta)^*P\Psi(z),
\end{eqnarray*}
where the second equality follows from \eqref{4.2a}, and the third
equality is a consequence of \eqref{2.8} and definition \eqref{4.3}
of $\Psi$.
Now \eqref{4.6} takes the form
\begin{equation}
\frac{1-w(\zeta)^*w(z)}{1-\bar{\zeta}z}=\frac{1}{V_{\cE}(\zeta)^*V_{\cE}(z)}
\left(\frac{1-\cE(\zeta)^*\cE(z)}{1-\bar{\zeta}z}
+\Psi(\zeta)^*P\Psi(z)\right).
\label{4.7a}
\end{equation}
\begin{Rk}
Equality \eqref{4.7a} implies that for every $\cE\in\cS_0$ and
$\Theta\in{\mathcal W}_\kappa$, the function $w={\bf T}_{\Theta}[\cE]$
belongs to the generalized Schur class $\cS_{\kappa'}$ for some
$\kappa'\le \kappa$.

\smallskip

{\rm Indeed, it follows from \eqref{4.7a} that
${\rm sq}_-K_w\le {\rm sq}_-K_{\cE}+{\rm sq}_-P=0+\kappa$.}
\label{R:dop}
\end{Rk}
Upon evaluating \eqref{4.7a}  at $\zeta=z$ we get
\begin{equation}
\frac{1-|w(z)|^2}{1-|z|^2}=\frac{1}{|V_{\cE}(z)|^2}
\left(\frac{1-|\cE(z)|^2}{1-|z|^2}
+\Psi(z)^*P\Psi(z)\right)
\label{4.7}
\end{equation}
and realize that boundary values of $w(t_i)$ and $d_w(t_i)$ can be
calculated from asymptotic formulas for $\Psi$, $U_{\cE}$, $V_{\cE}$
and $\cE$ as $z$ tends to one of the interpolation nodes $t_i$.
These asymptotic relations are presented in the next lemma.
\begin{La}
Let $\cE$ be a Schur function, let $\Psi$, $U_{\cE}$ and $V_{\cE}$
be  defined as in \eqref{4.3}, \eqref{4.4} and \eqref{4.5},
respectively , and let $t_i$ be an interpolation node.
Then the following asymptotic relations hold as $z$ tends to $t_i$
nontangentially:
\begin{eqnarray}
(z-t_i)\Psi(z)&=&\be_i\left(\te_i^*-\tc_i^*\cE(z)\right)+O(|z-t_i|),
\label{4.9}\\
(z-t_i)U_{\cE}(z)&=&w_i\left(\te_i^*-\tc_i^*\cE(z)\right)+O(|z-t_i|),
\label{4.11}\\
(z-t_i)V_{\cE}(z)&=&\left(\te_i^*-\tc_i^*\cE(z)\right)+O(|z-t_i|).
\label{4.12}
\end{eqnarray}
\label{L:dop1}
\end{La}
{\bf Proof:} Recall that $\be_i$ be the $i$-th column in the identity
matrix $I_n$. Since
$$
(z-t_i)(zI-T)^{-1}=\be_i\be_i^*+O(|z-t_i|)\quad \mbox{as} \; \;
z\to t_i,
$$
and since $\cE(z)$ is uniformly bounded on $\D$,
we have by \eqref{4.3},
\begin{eqnarray*}
(z-t_i)\Psi(z)&=&(z-t_i)(zI-T)^{-1}\left(\tE^*-\tC^*\cE(z)\right)\\
&=&\be_i\be_i^*\left(\tE^*-\tC^*\cE(z)\right)+O(|z-t_i|)
\end{eqnarray*}
which proves \eqref{4.9}, since ${\bf e}_i^*\tC^*=\tc_i^*$ and
${\bf e}_i^*\tE^*=\te_i^*$ by \eqref{2.1}.

\smallskip

Now we plug in the asymptotic relation \eqref{4.9} into the formulas
\eqref{4.4} and \eqref{4.9} for $U_{\cE}$ and $V_{\cE}$ and make use of
evident equalities
\begin{equation}
C(\mu I-T)^{-1}{\bf e}_i=\frac{w_i}{\mu-{t}_i}\quad\mbox{and}\quad
E(\mu I-T)^{-1}{\bf e}_i=\frac{1}{\mu-{t}_i}
\label{4.10}
\end{equation}
to get \eqref{4.11} and  \eqref{4.12}:
\begin{eqnarray*}
(z-t_i)U_{\cE}(z)&=&(z-t_i)\cE(z)-(z-t_i)(z-\mu)C(\mu I-T)^{-1}\Psi(z)\\
&=&(\mu -z)C(\mu I-T)^{-1}\be_i\left(\te_i^*-\tc_i^*\cE(z)\right)
+O(|z-t_i|)\\
&=& \frac{\mu -z}{\mu-t_i}w_i\left(\te_i^*-\tc_i^*\cE(z)
\right)+O(|z-t_i|)\\
&=&w_i\left(\te_i^*-\tc_i^*\cE(z)\right)+O(|z-t_i|),\\
(z-t_i)V_{\cE}(z)&=&(z-t_i)-(z-t_i)(z-\mu)E(\mu I-T)^{-1}\Psi(z)\\
&=&(\mu -z)E(\mu I-T)^{-1}\be_i\left(\te_i^*-\tc_i^*\cE(z)\right)
+O(|z-t_i|)\\
&=&\left(\te_i^*-\tc_i^*\cE(z)\right)+O(|z-t_i|).
\end{eqnarray*}
\begin{La}
Let $w\in\cS_\kappa$, let $t_0\in\T$, and let us assume that the limit
\begin{equation}
d:=\lim_{j\to \infty}\frac{1-|w(r_jt_0)|^2}{1-r_j^2}<\infty
\label{4.7d}
\end{equation}
exists and is finite for some sequence of numbers $r_j\in(0, \, 1)$
such that $\lim_{j\to\infty} r_j=1$. Then the nontangential limits
$d_w(t_0)$ and $w(t_0)$ (defined as in \eqref{0.2a} and \eqref{0.2b})
exist and moreover
\begin{equation}
d_w(t_0)=d\quad\mbox{and}\quad |w(t_0)|=1.
\label{4.7e}
\end{equation}
\label{L:dop}
\end{La}
{\bf Proof:} Since $w$ is a generalized Schur function, it admits the
Krein-Langer representation \eqref{0.6} and identity \eqref{0.7a} holds at
every point $z\in\D$. In particular,
\begin{equation}
\frac{1-|w(r_jt_0)|^2}{1-r_j^2}=
\frac{1}{|B(r_jt_0)|^2} \, \left(\frac{1-|S(r_jt_0)|^2}{1-r_j^2}-
\frac{1-|B(r_jt_0)|^2}{1-r_j^2}\right).
\label{4.7b}
\end{equation}
Since $B$ is a finite Blaschke product, it is analytic at $t_0$ and the
limit $d_B(t_0):={\displaystyle\lim_{z\to t_0}\frac{1-|B(z)|^2}{1-|z|^2}}$
exists and is finite. Assumption \eqref{4.7d} implies therefore that the
limit
$$
\lim_{j\to \infty}\frac{1-|S(r_jt_0)|^2}{1-r_j^2}=d+d_B(t_0)
$$
exists and is finite. Since $S\in\cS_0$, we then conclude by the
Carath\'eodory-Julia theorem (see e.g.,
\cite{sarasonangder, sarasonsubh, shapiro}) that
the nontangential limits $d_S(t_0)$ and $S(t_0)$ exist and moreover,
\begin{equation}
d_S(t_0)=d+d_B(t_0)\quad\mbox{and}\quad |S(t_0)|=1.
\label{4.7c}
\end{equation}
Now we pass to limits in \eqref{0.6} and \eqref{0.7a} as $z$ tends to
$t_0$ nontangentially to get
$$
w(t_0):=\lim_{z\to t_0}w(z)=\frac{S(t_0)}{B(t_0)}
\quad\mbox{and}\quad
d_w(t_0):=\lim_{z\to t_0}\frac{1-|w(z)|^2}{1-|z|^2}=d_S(t_0)-d_B(t_0)
$$
and relations \eqref{4.7c} imply now \eqref{4.7e} and complete the
proof.\qed
\begin{Tm}
If $\cE\in\cS_0$ meets condition ${\bf C}_1$ at $t_i$ (i.e.,
the nontangential boundary limit ${\displaystyle\lim_{z\to
t_i}\cE(z)}$ is not equal to $\q_i={\displaystyle\frac{\te_i^*}{\tc_i^*}}$
or fails to exist), then the function $w={\bf T}_{\Theta}[\cE]$ is
subject to
\begin{equation}
\lim_{z\rightarrow t_i} \, w(z)=w_i\quad\mbox{and}\quad
\lim_{z\rightarrow t_i} \,
\frac{1-|w(z)|^2}{1-|z|^2}=\gamma_i.
\label{4.8}
\end{equation}
\label{T:4.1}
\end{Tm}
{\bf Proof:} By the assumption of the theorem, there exists
$\varepsilon>0$ and a sequence of points $\{r_\alpha
t_i\}_{\alpha=1}^\infty$  tending to $t_i$ radially
($0<r_\alpha<1$ and $r_\alpha\to 1$) such that
\begin{equation}
|\te_i^*-\tc_i^*\cE(r_\alpha t_i)|\ge \varepsilon \quad \mbox{for
every} \; \; \alpha.
\label{4.13}
\end{equation}
Since $\be_i^*P\be_i=\gamma_i$ by the definition \eqref{0.10} of $P$, it
follows from \eqref{4.9} that
$$
|z-t_i|^2\Psi(z)^*P\Psi(z)=|\te_i^*-\tc_i^*\cE(z)|^2 \gamma_i+O(|z-t_i|).
$$
Furthermore, relation
$$
|z-t_i|^2  \cdot |V_{\cE}(z)|^2=|\te_i^*-\tc_i^*\cE(z)|^2+O(|z-t_i|)
$$
is a consequence of \eqref{4.12} and, since $\cE$ is uniformly bounded on
$\D$, it is clear that
$$
\lim_{z\to t_i}|z-t_i|^2 \cdot \frac{1-|\cE(z)|^2}{1-|z|^2}=0.
$$
Now we substitute the three last relations into \eqref{4.7} and let
$z=r_\alpha t_i\to t_i$; due to \eqref{4.13} we have
\begin{eqnarray*}
\lim_{z=r_\alpha t_i\to t_i}\frac{1-|w(z)|^2}{1-|z|^2}&=&
\lim_{z=r_\alpha t_i\to t_i}\frac{|z-t_i|^2 \cdot
{\displaystyle\frac{1-|\cE(z)|^2}{1-|z|^2}}+|z-t_i|^2\Psi(z)^*P\Psi(z)}
{|z-t_i|^2  \cdot |V_{\cE}(z)|^2}\nonumber\\
&=&\frac{0+\gamma_i}{1}=\gamma_i.
\end{eqnarray*}
Since $w$ is a generalized Schur function (by Remark \ref{R:dop}), we can
apply Lemma \ref{L:dop} to conclude that the nontangential limit
$d_w(t_i)$ exists and equals $\gamma_i$. This proves the second relation
in \eqref{4.8}. Furthermore, by \eqref{4.11} and \eqref{4.12} and in view
of \eqref{4.13},
\begin{equation}
\lim_{z=r_\alpha t_i\to t_i}w(z)=\lim_{z\to
t_i}\frac{(z-t_i)U_{\cE}(z)}{(z-t_i)V_{\cE}(z)}=w_i.
\label{4.15b}
\end{equation}
Again by Lemma \ref{L:dop}, the nontangential limit $w(t_i)$ exists;
therefore, it is equal to the subsequential limit \eqref{4.15b}, that is,
to $w_i$. This proves the first relation  in \eqref{4.8}
and completes the proof of the theorem.\qed

\medskip

The next step will will be to handle condition ${\bf C}_2$ (see
\eqref{0.26}). We need an auxiliary result.
\begin{La}
Let $t_0\in\T$ and let $\cE$ be a Schur function such that
\begin{equation}
\lim_{z\to t_0}\cE(z)=\cE_0 \; \; (|\cE_0|=1) \quad\mbox{and}\quad
\lim_{z\to t_0}\frac{1-|\cE(z)|^2}{1-|z|^2}=\infty.
\label{3.6a}
\end{equation}
Then
\begin{equation}
\lim_{z\to t_0}\frac{1-|\cE(z)|^2}{1-|z|^2}\cdot \left|
\frac{z-t_0}{\cE(z)-\cE_0}\right|^2=0\quad \mbox{and}\quad
\lim_{z\to t_0}\frac{z-t_0}{\cE(z)-\cE_0}=0.
\label{3.7}
\end{equation}
\label{L:3.1}
\end{La}
{\bf Proof:} Since $|\cE_0|=1$, we have
\begin{eqnarray*}
2{\rm Re} \, (1-\cE(z)\overline{\cE}_0)&=&(1-\cE(z)\overline{\cE}_0)
+(1-\cE_0\overline{\cE(z)})\\
&=&|1-\cE(z)\overline{\cE}_0|^2+1-|\cE_0|^2\cdot |\cE(z)|^2\\
&\ge& 1-|\cE(z)|^2
\end{eqnarray*}
and thus,
\begin{equation}
|\cE(z)-\cE_0|=|1-\cE(z)\overline{\cE_0}|\ge {\rm Re} \,
(1-\cE(z)\overline{\cE}_0)
\ge \frac{1}{2}\left(1-|\cE(z)|^2\right).
\label{aga}
\end{equation}
Furthermore, for every $z$ in the Stoltz domain
$$
\Gamma_a(t_0)=\{z\in\D: \; |z-t_0|<a(1-|z|^2)\}, \quad a>1,
$$
it holds that
$$
\frac{1-|z|^2}{|z-t_0|}\ge \frac{1-|z|}{|z-t_0|}>\frac{1}{a},
$$
which together with \eqref{aga} leads us to
$$
\left|\frac{\cE(z)-\cE_0}{z-t_0}\right|\ge \frac{1}{2} \cdot
\frac{1-|\cE(z)|^2}{|z-t_0|}=\frac{1}{2} \cdot
\frac{1-|\cE(z)|^2}{1-|z|^2}\cdot
\frac{1-|z|^2}{|z-t_0|}>\frac{1}{2a}\cdot
\frac{1-|\cE(z)|^2}{1-|z|^2}
$$
which is equivalent to
\begin{equation}
\frac{1-|\cE(z)|^2}{1-|z|^2}\cdot \left|
\frac{z-t_0}{\cE(z)-\cE_0}\right|\le 2a.
\label{aga1}
\end{equation}
Note that the denominator $\cE(z)-\cE_0$ in the latter inequality does not
vanish: assuming that $\cE(z_0)=\cE_0$ at
some point  $z_0\in\D$, we would have by the maximum modulus principle
(since $|\cE_0|=1$) that $\cE(z)\equiv \cE_0$ which would contradict the
second assumption in \eqref{3.6a}. Finally, by this latter assumption,
$d_{\cE}(t_0)=\infty$ and relations \eqref{3.7} follow immediately from
\eqref{aga1}.\qed
\begin{Tm}
Let $\cE\in\cS_0$ meet condition ${\bf C}_2$ at $t_i$:
\begin{equation}
\lim_{z\to t_i}\cE(z)=\q_i=\frac{\te_i^*}{\tc_i^*}\quad\mbox{and}\quad
\lim_{z\to t_i}\frac{1-|\cE(z)|^2}{1-|z|^2}=\infty.
\label{4.15}
\end{equation}
Then the function $w={\bf T}_{\Theta}[\cE]$ is subject to
relations \eqref{4.8}.
\label{T:4.2}
\end{Tm}
{\bf Proof:} Let for short
$$
\Delta_i(z):=\frac{\te_i^*-\tc_i^*\cE(z)}{t_i-z}
$$
and note that
\begin{equation}
\Delta_i(z)\neq 0 \quad (z\in\D).
\label{4.16}
\end{equation}
To see this we argue as in the proof of the previous lemma:
assuming that $\cE(z_0)=\q_i$ at some point  $z_0\in\D$, we would have by
the maximum modulus principle (since $|\q_i|=1$) that $\cE(z)\equiv \q_i$
which would contradict the second assumption in \eqref{4.15}.
Furthermore, since $|\q_i|=1$ and due to assumptions \eqref{4.15}, we can
apply Lemma \ref{L:3.1} (with  $\cE_0=\q_i$ and $t_0=t_i$) to conclude
that
\begin{equation}
\lim_{z\to t_0}\frac{1-|\cE(z)|^2}{1-|z|^2}\cdot
\frac{1}{|\Delta_i(z)|^2}=0
\label{4.17}
\end{equation}
and
\begin{equation}
\lim_{z\to t_0}\Delta_i(z)^{-1}=0.
\label{4.18}
\end{equation}
Now we divide both parts in asymptotic relations
\eqref{4.9}--\eqref{4.12} by $(\te_i^*-\tc_i^*\cE(z))$ and write
the obtained equalities in terms of $\Delta_i$ as
\begin{eqnarray*}
\Delta_i(z)^{-1}\Psi(z)&=&{\bf e}_i+\Delta_i(z)^{-1}\cdot O(1),\\
\Delta_i(z)^{-1}U_{\cE}(z)&=&w_i+\Delta_i(z)^{-1}\cdot O(1),\\
\Delta_i(z)^{-1}V_{\cE}(z)&=&1+\Delta_i(z)^{-1}\cdot O(1).
\end{eqnarray*}
By \eqref{4.18}, the following nontangential limits exist
$$
\lim_{z\to t_i}\Delta_i(z)^{-1}\Psi(z)={\bf e}_i, \quad
\lim_{z\to t_i}\Delta_i(z)^{-1}U_{\cE}(z)=w_i, \quad
\lim_{z\to t_i}\Delta_i(z)^{-1}V_{\cE}(z)=1
$$
and we use these limits along with \eqref{4.17} to pass to limits
in \eqref{4.7}:
\begin{eqnarray*}
\lim_{z\to t_i}\frac{1-|w(z)|^2}{1-|z|^2}
&=&\lim_{z\to t_i}\frac{|\Delta_i(z)|^{-2}
{\displaystyle\frac{1-|\cE(z)|^2}{1-|z|^2}}
+|\Delta_i(z)|^{-2}\Psi(z)^*P\Psi(z)}
{|\Delta_i(z)|^{-2}|V_{\cE}(z)|^2}\\
&=&\frac{0+{\bf e}_i^*P{\bf e}_i}{1}=\gamma_i.
\end{eqnarray*}
Finally,
$$
\lim_{z\rightarrow t_i} \, w(z)=\lim_{z\rightarrow t_i} \,
\frac{\Delta_i(z)^{-1}U_{\cE}(z)}{\Delta_i(z)^{-1}V_{\cE}(z)}
=\frac{w_i}{1}=w_i,
$$
which completes the proof.\qed
\begin{Tm}
Let $\widetilde{p}_{ii}$ be the $i$-th diagonal entry of
$P^{-1}=\left[\widetilde{p}_{ij}\right]_{i,j=1}^n$, let
$\cE\in\cS_0$ be subject to
\begin{equation}
\lim_{z\to t_i}\cE(z)=\q_i\quad\mbox{and}\quad
\lim_{z\to t_i}\frac{1-|\cE(z)|^2}{1-|z|^2}=d_\cE(t_i)<\infty.
\label{4.20}
\end{equation}
Let us assume that
\begin{equation}
d_\cE(t_i) \neq \frac{\widetilde{p}_{ii}}{|\te_i|^2}.
\label{4.21}
\end{equation}
Then the function $w:={\bf T}_{\Theta}[\cE]$ satisfies
\begin{equation}
\lim_{z\to t_i} \, w(z)=w_i
\label{4.22}
\end{equation}
and the nontangential limit $d_w(t_i):={\displaystyle\lim_{z\to
t_i}\frac{1-|w(z)|^2}{1-|z|^2}}$ is finite. Moreover,
\begin{equation}
d_w(t_i)<\gamma_i\quad\mbox{if} \quad
d_\cE(t_i)>-\frac{\widetilde{p}_{ii}}{|\te_i|^2}
\label{4.23}
\end{equation}
and
\begin{equation}
d_w(t_i)>\gamma_i\quad\mbox{if} \quad
d_\cE(t_i)<-\frac{\widetilde{p}_{ii}}{|\te_i|^2}.
\label{4.24}
\end{equation}
In other words, $d_w(t_i)<\gamma_i$ if $\cE$ meets condition ${\bf C}_3$
and $d_w(t_i)>\gamma_i$ if $\cE$ meets condition ${\bf C}_4$
at $t_i$.
\label{T:4.3}
\end{Tm}
{\bf Proof:} By the Carath\'eodory-Julia theorem (for Schur functions),
conditions \eqref{4.20} imply that the following nontangential limits
exist
$$
\lim_{z\to t_i} \, \cE^\prime(z)=
\lim_{z\to t_i} \, \frac{\cE(z)-\q_i}{z-t_i}=\bar{t}_i\q_id_\cE(t_i)
$$
and the following asymptotic equality holds
\begin{equation}
\cE(z)=\q_i+(z-t_i)\overline{t}_i\q_id_\cE(t_i)
+o(|z-t_i|)\quad \mbox{as} \; \; z\to t_i.
\label{4.25}
\end{equation}
We shall show that the functions $\Psi$, $U_{\cE}$ and $V_{\cE}$
defined in \eqref{4.3}, \eqref{4.4}, \eqref{4.5} admit the nontangential
boundary limits at every interpolation node $t_i$:
\begin{equation}
\Psi(t_i)=\frac{\bar{t}_i}{\te_i}\left(P^{-1}\be_i-
\be_i(\widetilde{p}_{ii}+|\te_i|^2d_\cE(t_i))\right),
\label{4.26}
\end{equation}
\begin{equation}
U_{\cE}(t_i)=-\frac{\bar{t}_iw_i}{\te_i}
(\widetilde{p}_{ii}+|\te_i|^2d_\cE(t_i))\quad\mbox{and}\quad
V_{\cE}(t_i)=-\frac{\bar{t}_i}{\te_i}
(\widetilde{p}_{ii}+|\te_i|^2d_\cE(t_i)).
\label{4.27}
\end{equation}
To prove \eqref{4.26} we first multiply both parts in the Stein identity
\eqref{2.3}, by $\be_i$ on the right and obtain
$$
P^{-1}\be_i-TP^{-1}T^*\be_i=\tE^*{\te}_i-\tC^*{\tc}_i
$$
which can be written equivalently, since $T^*\be_i=\bar{t}_i\be_i$ and
${\tc}_i={\te}_i\q_i$, as
\begin{equation}
\tE^*-\tC^*\q_i=\frac{\bar{t}_i}{\te_i}(t_i I-T)P^{-1}\be_i.
\label{4.28}
\end{equation}
Substituting \eqref{4.25} into \eqref{4.3} and making use of
\eqref{4.28} we get
\begin{eqnarray}
\Psi(z)&=&(zI-T)^{-1}\left(\tE^*-\tC^*\q_i\right)-
(z-t_i)(zI-T)^{-1}\tC^*\q_id_\cE(t_i)\bar{t}_i+o(1)\nonumber\\
&=&\frac{\bar{t}_i}{\te_i}(zI-T)^{-1}(t_i I-T)P^{-1}\be_i\nonumber\\
&&-(z-t_i)(zI-T)^{-1}\tC^*\q_id_\cE(t_i)\bar{t}_i+o(1).\label{4.29}
\end{eqnarray}
Since the following limits exist
$$
\lim_{z\to t_i}(zI-T)^{-1}(t_i I-T)=I-\be_i\be_i^*,\quad
\lim_{z\to t_i}(z-t_i)(zI-T)^{-1}=\be_i\be_i^*,
$$
we can pass to the limit in \eqref{4.29} as $z\to t_i$ nontangentially
to get
\begin{equation}
\Psi(t_i)=\frac{\bar{t}_i}{\te_i}(I-\be_i\be_i^*)P^{-1}\be_i-
\be_i\be_i^*\tC^*\q_id_\cE(t_i)\bar{t}_i.
\label{4.30}
\end{equation}
Since $\be_i^*P^{-1}\be_i=\widetilde{p}_{ii}$ and
$\be_i^*\tC^*\q_i=\tc_i^*\q_i=\te_i^*$, the right hand side expression
in \eqref{4.30} coincides with that in \eqref{4.26}.

\smallskip

Making use of \eqref{4.25} and \eqref{4.26} we pass to the limits
in \eqref{4.4} and \eqref{4.5} as $z\to t_i$ nontangentially:
\begin{eqnarray}
U_{\cE}(t_i)&=&\cE(t_i)-(t_i-\mu)C(\mu I-T)^{-1}\Psi(t_i)\nonumber\\
&=&\q_i-\frac{1-\mu\bar{t}_i}{\te_i}C(\mu I-T)^{-1}
\left(P^{-1}\be_i-\be_i(\widetilde{p}_{ii}+|\te_i|^2d_\cE(t_i))\right),
\label{4.30a}\\
V_{\cE}(t_i)&=&1-(t_i-\mu)E(\mu I-T)^{-1}\Psi(t_i)\nonumber\\
&=& 1-\frac{1-\mu\bar{t}_i}{\te_i}E(\mu I-T)^{-1}
\left(P^{-1}\be_i-\be_i(\widetilde{p}_{ii}+|\te_i|^2d_\cE(t_i))\right).
\label{4.30b}
\end{eqnarray}
Note that by \eqref{2.1},
\begin{eqnarray}
\frac{1-\mu\bar{t_i}}{{\te}_i}C(\mu I-T)^{-1}P^{-1}\be_i&=&
\frac{1-\mu\bar{t_i}}{{\te}_i}\tC(I-\mu T^*)^{-1}\be_i=
\frac{\tc_i}{\te_i}=\q_i,\label{4.31}\\
\frac{1-\mu\bar{t_i}}{{\te}_i}E(\mu I-T)^{-1}P^{-1}\be_i&=&
\frac{1-\mu\bar{t_i}}{{\te}_i}\tE(I-\mu T^*)^{-1}\be_i=
\frac{\te_i}{\te_i}=1.\label{4.32}
\end{eqnarray}
Making use of these two equalities we simplify \eqref{4.30a} and
\eqref{4.30b} to
$$
U_{\cE}(t_i)=\frac{1-\mu\bar{t}_i}{\te_i}C(\mu
I-T)^{-1}\be_i(\widetilde{p}_{ii}+|\te_i|^2d_\cE(t_i))
$$
and
$$
V_{\cE}(t_i)=\frac{1-\mu\bar{t}_i}{\te_i}E(\mu
I-T)^{-1}\be_i(\widetilde{p}_{ii}+|\te_i|^2d_\cE(t_i)),
$$
respectively, and it is readily seen from \eqref{4.10} that the two
latter equalities coincide with those in \eqref{4.27}.

\smallskip

Now we conclude from \eqref{4.2} and \eqref{4.27} that the nontangential
boundary limits $w(t_i)$ exist for $i=1,\ldots,n$ and
$$
w(t_i)=\lim_{z\to t_i}w(z)=\lim_{z\to t_i}\frac{U_{\cE}(z)}{V_{\cE}(z)}
=\frac{U_{\cE}(t_i)}{V_{\cE}(t_i)}=w_i
$$
which proves \eqref{4.22}. Furthermore, since the nontangential
boundary limits $d_\cE(t_i)$ and
\begin{equation}
|V_{\cE}(t_i)|^2=\frac{(\widetilde{p}_{ii}+|\te_i|^2d_\cE(t_i))^2}{|\te_i|^2}
\label{4.30c}
\end{equation}
exist (by the second assumption in \eqref{4.20} and the second relation
in \eqref{4.27}), we can pass to the limit in \eqref{4.7} as $z$ tends to
$t_i$ nontangentially:
$$
d_w(t_i)=\frac{d_\cE(t_i)+\Psi(t_i)^*P\Psi(t_i)}{|V_{\cE}(t_i)|^2}.
$$
By \eqref{4.30c} and \eqref{4.26} we have
$$
d_w(t_i)=\frac{|\te_i|^2d_\cE(t_i)+
\left(\be_i^*P^{-1}-(\widetilde{p}_{ii}+|\te_i|^2d_\cE(t_i))\be_i^*\right)
P\left(P^{-1}\be_i-\be_i(\widetilde{p}_{ii}+|\te_i|^2d_\cE(t_i))\right)}
{(\widetilde{p}_{ii}+|\te_i|^2d_\cE(t_i))^2}
$$
and elementary algebraic transformations based on equalities
$\be_i^*P^{-1}\be_i=\widetilde{p}_{ii}$,
$\be_i^*P\be_i=\gamma_i$ and $\be_i^*\be_i=1$ lead us to
\begin{equation}
d_w(t_i)=\gamma_i-\frac{1}{\widetilde{p}_{ii}+|\te_i|^2d_\cE(t_i)}.
\label{4.34}
\end{equation}
Statements \eqref{4.23} and \eqref{4.24} follow immediately from
\eqref{4.34}.\qed

\medskip

As we have already mentioned in Introduction, Theorem \ref{T:0.8}
is known for the case $\kappa=0$ (see \cite{Sarasonnp})
At this point we already can recover this result.
\begin{Tm}
Let the Pick matrix $P$ be positive definite and let $T$, $E$, $C$,
$\Theta(z)$ and $\q_i$ be defined as in \eqref{0.18}, \eqref{0.17} and
\eqref{0.25a}. Then all solutions $w$ of Problem \ref{P:0.2} are
parametrized by the formula \eqref{0.24}
when the parameter $\cE$ belongs to the Schur class $\cS_0$ and satisfies
condition ${\bf C}_1\vee{\bf C}_2$ at each interpolation
node: either $\cE$ fails to admit the nontangential boundary limit $\q_i$
at $t_i$ or
$$
\cE(t_i)=\q_i\quad\mbox{and}\quad d_\cE(t_i)=\infty.
$$
\label{T:dop3}
\end{Tm}
{\bf Proof:} Any solution $w$ of Problem \ref{P:0.2} is a solution of
Problem \ref{P:0.1} and then by Statement 1 in Theorem \ref{T:2.6},
it is of the form $w={\bf T}_{\Theta}[\cE]$ for some Schur class function
$\cE$. Since $P>0$, the diagonal entries $\widetilde{p}_{ii}$ of $P^{-1}$
are positive. Therefore, the cases specified in \eqref{0.28}--\eqref{0.30}
(conditions ${\bf C}_4-{\bf C}_6$ cannot occur in this situation, whereas
condition ${\bf C}_3$ simplifies to
$$
{\bf C}_3: \quad \cE(t_i)=\q_i\quad\mbox{and}\quad
d_{\cE}(t_i)<\infty.
$$
In other words, any function $\cE\in\cS_0$ satisfies exactly one of the
conditions ${\bf C}_1$, ${\bf C}_2$ or ${\bf C}_3$ at each one of
interpolation nodes. Therefore, once $\cE$ does not meet condition ${\bf
C}_1$ or condition ${\bf C}_2$ at at least one interpolation node $t_i$,
it meets condition ${\bf C}_3$ at $t_i$. Therefore, it holds for the
function $w={\bf T}_{\Theta}[\cE]$ that
$d_w(t_i)<\gamma_i$ (by Theorem \ref{T:4.3}) and therefore $w$ is not a
solution of Problem \ref{P:0.2}. On the other hand, if $\cE$ meets
condition ${\bf C}_1\vee {\bf C}_2$ at every interpolation node, then
$w={\bf T}_{\Theta}[\cE]$ satisfies interpolation conditions \eqref{4.8}
(by Theorems \ref{T:4.1} and \ref{T:4.2}) that means that $w$ is a
solution of Problem \ref{P:0.2}.\qed
\begin{Rk}
{\rm It is useful to note that for the one-point interpolation problem
(i.e.,
when $n=1$), definition \eqref{2.2} takes the form
$$
\left[\begin{array}{c}\tc_1 \\ \te_1\end{array}\right]
=\left[\begin{array}{c}w_1 \\ 1\end{array}\right]
(\mu-t_1)^{-1}\gamma_1^{-1}(I-\mu
\bar{t}_1)=-\bar{t}_1\left[\begin{array}{c}w_1 \\
1\end{array}\right]\gamma_1^{-1}
$$
and therefore the number $\q_1:=\frac{\tc_1}{\te_1}$ in this case is equal
to $w_1$.}
\label{R:dop4}
\end{Rk}
Now we turn back to the indefinite case. Theorems \ref{T:4.4} and
\ref{T:4.5} below treat the case when condition \eqref{4.21} is dropped.
For notational convenience we let $i=n$ and
$$
T_1=\left[\begin{array}{ccc}t_1 &&
\\ &\ddots & \\ && t_{n-1}\end{array}\right],\quad
E_1=\begin{bmatrix} 1 & \ldots & 1\end{bmatrix},\quad
C_1=\begin{bmatrix} w_1 & \ldots & w_{n-1}\end{bmatrix}
$$
so that decompositions
\begin{equation}
T=\left[\begin{array}{cc} T_1 & 0 \\ 0 & t_n \end{array}\right],
\quad E=\begin{bmatrix} E_1 & 1\end{bmatrix},\quad C=\begin{bmatrix} C_1 &
w_n\end{bmatrix}
\label{4.35}
\end{equation}
are conformal with partitionings
\begin{equation}
P=\begin{bmatrix}P_{11} & P_{12} \\ P_{21} & \gamma_n\end{bmatrix}
\quad\mbox{and}\quad P^{-1}=\begin{bmatrix}\tP_{11} & \tP_{12} \\ \tP_{21}
& \widetilde{p}_{nn}\end{bmatrix}.
\label{4.36}
\end{equation}
\begin{Tm}
Let $\widetilde{p}_{nn}<0$ and let $\cE$ be a Schur function such that
\begin{equation}
\lim_{z\to t_n}\cE(z)=\q_n\quad\mbox{and}\quad
d_\cE(t_n)=-\frac{\widetilde{p}_{nn}}{|\te_n|^2}.
\label{4.37}
\end{equation}
Then the function
\begin{equation}
w:={\bf T}_{\Theta}[\cE]
\label{4.38}
\end{equation}
is subject to one of the following:
\begin{enumerate}
\item The nontangential boundary limit $w(t_n)$ does not exist.
\item The latter limit exists and $w(t_n)\neq w_n$.
\item The latter limit exists, is equal to $w_n$ and $d_w(t_n)=\infty$.
\end{enumerate}
\label{T:4.4}
\end{Tm}
{\bf Proof:} Since $\cE$ is a Schur function,
conditions \eqref{4.37} form a well posed one-point
interpolation problem (similar to Problem \ref{P:0.2}).
By Theorem \ref{T:dop3}, $\cE$ admits a representation
\begin{equation}
\cE={\bf T}_{\widehat{\Theta}}[\widehat{\cE}]
\label{4.39}
\end{equation}
with the coefficient matrix $\widehat{\Theta}$ defined via formula
\eqref{0.17}, but with $P$, $T$, $E$ and $C$ replaced by
$-\frac{\widetilde{p}_{nn}}{|\te_n|^2}$, $t_n$, $1$ and $\q_n$,
respectively:
\begin{equation}
\widehat{\Theta}(z)=I_2-\frac{z-\mu}{(z-t_n)(1-\mu
\overline{t}_n)}\left[\begin{array}{c}\q_n \\
1\end{array}\right]\frac{|\te_n|^2}{\widetilde{p}_{nn}}
\left[\begin{array}{cc}\q_n^* & -1\end{array}\right]
\label{4.40}
\end{equation}
and a parameter $\widehat{\cE}\in\cS_0$ satisfying one of the following
three conditions:
\begin{enumerate}
\item[(a)] The limit $\widehat{\cE}(t_n)$ does not exist.
\item[(b)] The limit $\widehat{\cE}(t_n)$ exists and is not equal to
$\q_n$.
\item[(c)] It holds that
\begin{equation}
\widehat{\cE}(t_n)=\q_n\quad\mbox{and}\quad
d_{\widehat{\cE}}(t_n)=\infty.
\label{4.41}
\end{equation}
\end{enumerate}
We shall show that conditions (a), (b) and (c) for the parameter
$\widehat{\cE}$
are equivalent to statements (1), (2) and (3), respectively, in the
formulation of the theorem. This will complete the proof.

\smallskip

Note that $\q_n$ appearing in (a) and (b) is the same as in \eqref{4.37},
due to Remark \ref{R:dop4}. Since
$\q_n={\displaystyle\frac{\tc_n}{\te_n}}$, we can write \eqref{4.40} as
$$
\widehat{\Theta}(z)=I_2-\frac{z-\mu}{(z-t_n)(1-\mu
\overline{t}_n)}\left[\begin{array}{c}\tc_n \\
\te_n\end{array}\right]\frac{1}{\widetilde{p}_{nn}}
\left[\begin{array}{cc}\tc_n^* & -\te_n^*\end{array}\right]
$$
The inverse of $\widehat{\Theta}$ equals
\begin{equation}
\widehat{\Theta}(z)^{-1}=I_2+\frac{z-\mu}{(z-t_n)(1-\mu
\overline{t}_n)}\left[\begin{array}{c}\tc_n \\
\te_n\end{array}\right]\frac{1}{\widetilde{p}_{nn}}
\left[\begin{array}{cc}\tc_n^* & -\te_n^*\end{array}\right]
\label{4.42u}
\end{equation}
and coincides with the function $\widehat{\Theta}^{(2)}$ in \eqref{2.20}.
Therefore, by Lemma \ref{L:2.5} and by Remark \ref{R:2.6},
\begin{equation}
\Theta(z)=\Theta^{(1)}(z)\widehat{\Theta}(z)^{-1}
\label{4.42}
\end{equation}
where $\Theta^{(1)}$ is given in \eqref{2.11b}. Substituting \eqref{4.40}
into \eqref{4.38} (that is, representing $w$ as a result of composition
of two linear fractional transformations) and taking into account
\eqref{4.42} we get
$$
w:={\bf T}_{\Theta}[\cE]=
{\bf T}_{\Theta}[{\bf T}_{\widehat{\Theta}}[\widehat{\cE}]]=
{\bf T}_{\Theta\widehat{\Theta}}[\widehat{\cE}]=
{\bf T}_{\Theta^{(1)}}[\widehat{\cE}].
$$
Thus, upon setting
\begin{equation}
U_{\widehat{\cE}}(z)=\Theta^{(1)}_{11}(z)\widehat{\cE}(z)+
\Theta^{(1)}_{12}(z),\qquad
V_{\widehat{\cE}}(z)=\Theta^{(1)}_{21}(z)\widehat{\cE}(z)+
\Theta^{(1)}_{22}(z),
\label{4.43}
\end{equation}
we have
\begin{equation}
w={\bf T}_{\Theta^{(1)}}[\widehat{\cE}]=\frac{\Theta^{(1)}_{11}
\widehat{\cE}+\Theta^{(1)}_{12}}{\Theta^{(1)}_{21}\widehat{\cE}+
\Theta^{(1)}_{22}}=\frac{U_{\widehat{\cE}}}{V_{\widehat{\cE}}}.
\label{4.44}
\end{equation}
Note that $\Theta^{(1)}$ is a rational function analytic and invertible
at $t_n$. It follows immediately from \eqref{4.44} that if the boundary
limit $\widehat{\cE}(t_n)$ does not exist, then the boundary $w(t_n)$
does not exist either. Thus, $(a)\Rightarrow (1)$. The rest is broken into
two steps.

\medskip

{\bf Step 1:} {\em Let the nontangential boundary limit
$\widehat{\cE}(t_n)$ exists. Then so do the limits
$U_{\widehat{\cE}}(t_n)$, $V_{\widehat{\cE}}(t_n)$ and $w(t_n)$, and
moreover,}
\begin{equation}
V_{\widehat{\cE}}(t_n):=\lim_{z\to t_n}V_{\widehat{\cE}}(z)\neq 0
\label{4.46}
\end{equation}
and
\begin{equation}
w(t_n)=w_n \quad\mbox{if and only if}\quad
\widehat{\cE}(t_n)=\q_n.
\label{4.49}
\end{equation}
{\bf Proof of Step 1:} Existence of the limits $U_{\widehat{\cE}}(t_n)$
and $V_{\widehat{\cE}}(t_n)$ is clear since $\Theta^{(1)}$ is analytic at
$t_n$. Assume that $V_{\widehat{\cE}}(t_n)=0$. Then
$U_{\widehat{\cE}}(t_n)=0$, since otherwise, the function $w$ of the form
\eqref{4.44} would not be bounded in a neighborhood of $t_n\in\T$ which
cannot occur since $w$ is a generalized Schur function. If
$V_{\widehat{\cE}}(t_n)=U_{\widehat{\cE}}(t_n)=0$, then it follows
from \eqref{4.43} that
$$
\Theta^{(1)}(t_n)\begin{bmatrix}\widehat{\cE}(t_n) \\ 1\end{bmatrix}
=\begin{bmatrix}U_{\widehat{\cE}}(t_n) \\
V_{\widehat{\cE}}(t_n)\end{bmatrix}=0
$$
and thus, the matrix $\Theta^{(1)}(t_n)$ is singular which is a
contradiction. Now it follows from \eqref{4.44} and \eqref{4.46} that
the limit $w(t_n)$ exists. This completes the proof of
$(a)\Leftrightarrow
(1)$. The proof of \eqref{4.49} rests on the
equality
\begin{equation}
\left[\begin{array}{cc}w_n^* & -1\end{array}\right]
\Theta^{(1)}(t_n)=\frac{\bar{t}_n}{\widetilde{p}_{nn}}
\left[\begin{array}{cc}\tc_n^*& -\te_n^*\end{array}\right].
\label{4.45}
\end{equation}
Indeed, it follows from \eqref{4.44} and \eqref{4.45} that
\begin{eqnarray*}
w(t_n)-w_n&=&\frac{U_{\widehat{\cE}}(t_n)-
w_nV_{\widehat{\cE}}(t_n)}{V_{\widehat{\cE}}(t_n)}\\
&=&\frac{w_n}{V_{\widehat{\cE}}(t_n)}\cdot
\begin{bmatrix}{w}_n^* & -1\end{bmatrix}\Theta^{(1)}(t_n)
\begin{bmatrix}\widehat{\cE}(t_n) \\ 1\end{bmatrix}\\
&=& \frac{\bar{t}_nw_n}{\widetilde{p}_{nn}V_{\widehat{\cE}}(t_n)}
\begin{bmatrix}\tc_n^* & -\te_n^*\end{bmatrix}
\begin{bmatrix}\widehat{\cE}(t_n) \\ 1\end{bmatrix}\\
&=& \frac{\bar{t}_n w_n}{\widetilde{p}_{nn}\tc_n^*
V_{\widehat{\cE}}(t_n)}\left(\widehat{\cE}(t_n)-\q_n\right)
\end{eqnarray*}
which clearly implies \eqref{4.49}. It remains to prove \eqref{4.45}.
To this end, note  that by \eqref{uh},
$$
{\rm Res}_{z=t_n}\Theta(z)=-\left[\begin{array}{c}w_n \\ 1\end{array}
\right]\left[\begin{array}{cc}\tc_n^* & -\te_n^*\end{array}\right]
$$
and it is readily seen from \eqref{4.42u} that
$$
{\rm Res}_{z=t_n}\widehat{\Theta}(z)^{-1}=t_n\left[\begin{array}{c}\tc_n
\\ \te_n\end{array}\right]\frac{1}{\widetilde{p}_{nn}}
\left[\begin{array}{cc}\tc_n^* & -\te_n^*\end{array}\right].
$$
Taking into account that $\Theta^{(1)}$ is analytic at $t_n$ and
that $\Theta$ and $\widehat{\Theta}^{-1}$ have simple poles at $t_n$,
we compare the residues of both parts in \eqref{4.42} at $t_n$ to arrive
at
$$
-\left[\begin{array}{c}w_n \\ 1\end{array}
\right]\left[\begin{array}{cc}\tc_n^* & -\te_n^*\end{array}\right]=
\frac{t_n}{\widetilde{p}_{nn}}\Theta^{(1)}(t_n)\left[\begin{array}{c}\tc_n
\\ \te_n\end{array}\right]\left[\begin{array}{cc}\tc_n^* &
-\te_n^*\end{array}\right],
$$
which implies (since $\te_n\neq 0$)
$$
\left[\begin{array}{c}w_n \\ 1\end{array}\right]=
\Theta^{(1)}(t_n)\left[\begin{array}{c}\tc_n
\\ \te_n\end{array}\right]\frac{t_n}{\widetilde{p}_{nn}}.
$$
Equality of adjoints in the latter equality gives
$$
\begin{bmatrix}w_n^* & -1\end{bmatrix}=
\begin{bmatrix}w_n^* & 1\end{bmatrix}J=
\frac{\bar{t}_n}{\widetilde{p}_{nn}}
\begin{bmatrix}\tc_n^* & -\te_n^*\end{bmatrix}J\Theta^{(1)}(t_n)^*J
$$
which is equivalent to \eqref{4.45}, since $\Theta^{(1)}(t_n)$ is
$J$-unitary and thus, $J\Theta^{(1)}(t_n)^*J=\Theta^{(1)}(t_n)^{-1}$.
This completes the proof of \eqref{4.49} which implies in particular,
that $(b)\Leftrightarrow (2)$.

\medskip

{\bf Step 2:} $(c)\Leftrightarrow (3)$.

\medskip

{\bf Proof of Step 2:} Equality $w(t_n)=w_n$ is equivalent to the
first condition in \eqref{4.41} by \eqref{4.49}.
To complete the proof, it suffices to show that if
$\widehat{\cE}(t_n)=\q_n$, then
\begin{equation}
d_w(t_n)=\infty\quad\mbox{if and only if}\quad
d_{\widehat{\cE}}(t_n)=\infty.
\label{4.580}
\end{equation}
To this end, we
write a virtue of relation \eqref{4.7} in terms of the parameter
$\widehat{\cE}$:
\begin{equation}
\frac{1-|w(z)|^2}{1-|z|^2}=\frac{1}{|V_{\widehat{\cE}}(z)|^2}
\left(\frac{1-|\widehat{\cE}(z)|^2}{1-|z|^2}
+\widehat{\Psi}(z)^*P\widehat{\Psi}(z)\right)
\label{4.50}
\end{equation}
where
\begin{equation}
\widehat{\Psi}(z)=(zI-T_1)^{-1}(\mu I-T_1)P_{11}^{-1}(I-\mu T_1^*)^{-1}
\left(E_1^*-C_1^*\widehat{\cE}(z)\right).
\label{4.51}
\end{equation}
Note that to get \eqref{4.51} we represent the right hand side
expression in \eqref{4.3} in terms of $C$ and $E$ (rather than
$\tC$ and $\tE$; this can be achieved with help of \eqref{2.2})
and then replace $P$, $T$, $E$, $C$ and $\cE$ in the obtained formula
by $P_{11}$, $T_1$, $E_1$, $C_1$ and $\widehat{\cE}$, respectively.
Since the nontangential boundary limit
$$
\widehat{\Psi}(t_n)=(t_nI-T_1)^{-1}(\mu I-T_1)P_{11}^{-1}(I-\mu
T_1^*)^{-1} \left(E_1^*-C_1^*\q_n\right)
$$
exists and is finite, equivalence \eqref{4.580} follows from
\eqref{4.50}.\qed
\begin{Tm}
Let $\widetilde{p}_{nn}=0$ and let $\cE$ be a Schur function such that
\begin{equation}
\cE(t_n)=\q_n\quad\mbox{and}\quad d_\cE(t_n)=0.
\label{4.52}
\end{equation}
Then the function $w:={\bf T}_{\Theta}[\cE]$ admits finite
nontangential boundary limits $d_w(t_n)$ and $w(t_n)\neq w_n$.
\label{T:4.5}
\end{Tm}
{\bf Proof:} Conditions \eqref{4.52} state a one-point boundary
interpolation problem for Schur functions $\cE$ with the Pick matrix
equals $d_\cE(t_n)=0$. Then by Statement 2 in Theorem \ref{T:2.6},
the only function $\cE$ satisfying conditions \eqref{4.52} is the
constant function $\cE(z)\equiv \q_n$ (the Blaschke product of degree zero).
Since $|\q_n|=1$, the function $w={\bf T}_{\Theta}[\cE]$ is rational and
unimodular on $\T$. Therefore, it is equal to the ratio of two finite
Blaschke products and therefore, the limits $w(t)$ and $d_w(t)$ exist at
every point $t\in\T$. We shall use decompositions \eqref{4.35} and
\eqref{4.36} with understanding that $\widetilde{p}_{nn}=0$, so that
\begin{equation}
\tP_{21}P_{12}=1\quad\mbox{and}\quad P^{-1}{\bf
e}_n=\begin{bmatrix}\tP_{12} \\ 0\end{bmatrix}.
\label{4.54}
\end{equation}
We shall also make use the formula
\begin{equation}
P_{21}(I-\bar{t}_n T_1)^{-1}=(E_1-w_n^*C_1)
\label{4.60}
\end{equation}
that follows from the Stein identity \eqref{1.1u} upon substituting
partitionings \eqref{4.35}, \eqref{4.36} and comparison the $(1,2)$
block entries.

\smallskip

In the current context, the formula \eqref{4.3} for $\Psi$ simplifies, on
account of \eqref{4.28}, to
\begin{eqnarray*}
\Psi(z)&=&(zI-T)^{-1}\left(\tE^*-\tC^*\q_n\right)\\
&=&\frac{\bar{t}_n}{\te_n}(zI-T)^{-1}(t_nI-T)P^{-1}{\bf e}_n
\end{eqnarray*}
Now we substitute the latter equality into \eqref{4.4} and
\eqref{4.5} and use formulas \eqref{4.31} and \eqref{4.32} (for $i=n$)
to get
$$
U_{\cE}(z)=\frac{1-z\bar{t}_n}{\te_n}C(zI-T)^{-1}P^{-1}{\bf e}_n,
\quad V_{\cE}(z)=\frac{1-z\bar{t}_n}{\te_n}E(zI-T)^{-1}P^{-1}{\bf e}_n.
$$
Taking into account the second equality in \eqref{4.54}, rewrite the
latter two formulas in terms of partitionings \eqref{4.35} and
\eqref{4.36} as
\begin{equation}
U_{\cE}(z)=\frac{1-z\bar{t}_n}{\te_n}C_1(zI-T_1)^{-1}\tP_{12},
\quad V_{\cE}(z)=\frac{1-z\bar{t}_n}{\te_n}E_1(zI-T_1)^{-1}\tP_{12}.
\label{4.59}
\end{equation}
Thus,
$$
w(z):=\frac{U_{\cE}(z)}{V_{\cE}(z)}=\frac{C_1(zI-T_1)^{-1}\tP_{12}}
{E_1(zI-T_1)^{-1}\tP_{12}}.
$$
We shall show that the denominator on the right hand side in the latter
formula does not vanish at $z=t_n$, so that
\begin{equation}
w(t_n):=\lim_{z\to t_n}\frac{C_1(zI-T_1)^{-1}\tP_{12}}
{E_1(zI-T_1)^{-1}\tP_{12}}=\frac{C_1(t_n I-T_1)^{-1}\tP_{12}}
{E_1(t_n I-T_1)^{-1}\tP_{12}}
\label{4.61}
\end{equation}
Then we will have, on account of \eqref{4.60} and the first
equality in \eqref{4.54},
\begin{eqnarray}
w_n-w(t_n)&=&w_n-\frac{C_1(t_n I-T_1)^{-1}\tP_{12}}
{E_1(t_n I-T_1)^{-1}\tP_{12}}\nonumber\\
&=&\frac{(w_n E_1-C_1)(t_n I-T_1)^{-1}\tP_{12}}
{E_1(t_n I-T_1)^{-1}\tP_{12}}\nonumber\\
&=&\frac{w_nt_n(E_1-w_n^*C_1)(I-\bar{t}_nT_1)^{-1}\tP_{12}}
{E_1(t_n I-T_1)^{-1}\tP_{12}}\nonumber\\
&=&\frac{w_nt_n\tP_{21}P_{12}}{E_1(t_n I-T_1)^{-1}\tP_{12}}
=\frac{w_nt_n}{E_1(t_n I-T_1)^{-1}\tP_{12}}\neq 0\label{4.62}
\end{eqnarray}
and thus $w(t_n)\neq w_n$. Thus, it remains to show that the
denominator in \eqref{4.61} is not zero. Assume that
$E_1(t_n I-T_1)^{-1}\tP_{12}=0$.
Since the limit in \eqref{4.61} exists (recall that $w$ is the ratio of
two finite Blaschke products), the latter assumption forces
$C_1(t_n I-T_1)^{-1}\tP_{12}=0$ and therefore, equality
$$
(w_nE_1-C_1)(t_n I-T_1)^{-1}\tP_{12}=0.
$$
But it was already shown in calculation \eqref{4.62} that
$$
(w_nE_1-C_1)(t_n I-T_1)^{-1}\tP_{12}=w_nt_n\neq 0
$$
and the obtained contradiction completes the proof.\qed

\medskip

Recall that the interpolation node $t_n$ in Theorems \ref{T:4.4} and
\ref{T:4.5} was chosen just for notational convenience and can be replaced
by any interpolation node $t_i$. It means that Theorems \ref{T:4.4} and
\ref{T:4.5} prove Statements $(4)$ and $(5)$ in Theorem \ref{T:0.20}.
Furthermore, Theorem \ref{T:4.3} proves the ``if'' parts in Statements
$(4)$ and $(5)$ in Theorem \ref{T:0.20}, whereas Theorems \ref{T:4.1} and
\ref{T:4.2} prove the ``if'' part in Statement $(1)$ in Theorem
\ref{T:0.20}. Finally since conditions ${\bf C}_1$-${\bf C}_6$ are
disjoint, the ``only if'' parts in Statements $(1)$, $(2)$ and  $(3)$ are
obvious. This completes the proof of Theorem \ref{T:0.20}.

\section{Negative squares of the function $w={\bf T}_\Theta[\cE]$.}
\setcounter{equation}{0}

In this section we prove Theorems \ref{T:0.15} and \ref{T:0.21}.
We assume without loss of generality that (maybe after an appropriate
rearrangement of the interpolation nodes) a fixed parameter
$\cE\in\cS_0$ satisfies condition ${\bf C}_{1-3}$ at interpolation nodes
$t_1,\ldots,t_{n-\ell}$ and conditions ${\bf C}_{4-6}$ at the remaining
$\ell$ points. Thus, we assume that
\begin{equation}
\lim_{z\to t_i} \cE(z)=\q_i\quad \mbox{and}\quad
\lim_{z\to t_i} \frac{1-|\cE(z)|^2}{1-|z|^2}\le
-\frac{\widetilde{p}_{ii}}{|\te_i|^2} \quad (i=n-\ell+1,\ldots,n).
\label{5.1}
\end{equation}
Let
\begin{equation}
P^{-1}=\left[\begin{array}{cc} \tP_{11} & \tP_{12} \\ \tP_{21} &
\tP_{22}\end{array}\right]\quad\mbox{with}\quad \tP_{22}\in\C^{\ell\times
\ell}.
\label{5.2}
\end{equation}
Note that under the above assumption, the matrix ${\mathcal P}$ in the
formulation of Theorem \ref{T:0.15} coincides with $\tP_{22}$ in the
decomposition \eqref{5.2}. Thus, to prove Theorem \ref{T:0.15}, it
suffices to show that there exists a Schur function $\cE$ satisfying
conditions \eqref{5.1} if and only if the matrix $\tP_{22}$ is negative
semidefinite.

\medskip

{\bf Proof of Theorem \ref{T:0.15}:} Since $|\q_i|=1$, conditions
\eqref{5.1} form a well posed boundary Nevanlinna-Pick  problem
(similar to Problem \ref{P:0.1}) in the Schur class $\cS_0$. This
problem has a solution $\cE$ if and only if the corresponding Pick
matrix
\begin{equation}
{\mathbb P}=[{\mathbb P}_{ij}]_{i,j=n-\ell+1}^{n}\quad \mbox{with the
entries}\quad {\mathbb P}_{ij}=\left\{\begin{array}{ccc} {\displaystyle
\frac{1-\q_i^*\q_j}{1-\bar{t}_{i}t_{j}}}
&\mbox{for} & i\neq j,\\
-{\displaystyle\frac{\widetilde{p}_{ii}}{|\te_{i}|^2}}
&\mbox{for}&  i=j,\end{array}\right.
\label{5.3}
\end{equation}
is positive semidefinite. Furthermore, there exist infinitely many
functions $\cE$ satisfying \eqref{5.1} if ${\mathbb P}$ is positive
definite and there is a unique such function (which is a Blaschke product
of degree equals ${\rm rank} \, {\mathbb P}$) if ${\mathbb P}$ is
singular. Thus, to complete the proof, it suffices to show that
\begin{equation}
{\mathbb P}>0\Longleftrightarrow \tP_{22}<0,\quad
{\mathbb P}\ge 0\Longleftrightarrow \tP_{22}\le 0\quad\mbox{and}\quad
{\rm rank} \, {\mathbb P}={\rm rank} \,  \tP_{22}.
\label{5.4}
\end{equation}
To this end, note that
\begin{equation}
\bar{t}_i\te_{i}^*\cdot {\mathbb P}_{ij} \cdot
t_j\te_{j}=-\widetilde{p}_{ij} \quad
(i,j=n-\ell+1,\ldots,n)
\label{5.6}
\end{equation}
where $\widetilde{p}_{ij}$ is the $ij$-th entry in $P^{-1}$.
Indeed, if $i\neq j$, then \eqref{5.6} follows from \eqref{5.3},
\eqref{2.4a} and
definition \eqref{0.25a} of $\q_i$. If $i=j$, then \eqref{5.6} follows
directly from \eqref{5.3}. By \eqref{5.2},
$\left[\widetilde{p}_{ij}\right]_{i,j=\ell+1}^{n}=\tP_{22}$, which
allows us to rewrite equalities \eqref{5.6} in the matrix form as
\begin{equation}
{\bf C}^* \, {\mathbb P} \, {\bf C}=-\tP_{22}\quad\mbox{where}\quad
{\bf C}={\rm diag} \, \left(t_{\ell+1}\te_{\ell+1},
\; t_{\ell+2}\te_{\ell+2}, \, \ldots, \, t_{n}\te_{n}\right).
\label{5.7}
\end{equation}
Since the matrix ${\bf C}$ is invertible, all the statements
in \eqref{5.4} follow from \eqref{5.7}. This completes the proof of
Theorem \ref{T:0.15}.\qed

\medskip

To prove Theorem \ref{T:0.21} we shall use the following result
(see \cite[Lemma 2.4]{bol} for the proof).
\begin{La}
Let $P\in\C^{n\times n}$ be an invertible Hermitian matrix and let
\begin{equation}
P=\left[\begin{array}{cc} P_{11} & P_{12} \\ P_{21} &
P_{22}\end{array}\right]\quad\mbox{and}\quad
P^{-1}=\left[\begin{array}{cc} \tP_{11} & \tP_{12} \\ \tP_{21} &
\tP_{22}\end{array}\right]
\label{5.8}
\end{equation}
be two conformal decompositions of $P$ and of $P^{-1}$ with
$P_{22}, \, \tP_{22}\in\C^{\ell\times \ell}$.
Furthermore, let $\tP_{22}$ be negative semidefinite. Then
$$
{\rm sq}_-P_{11}={\rm sq}_-P-\ell.
$$
\label{L:5.1}
\end{La}
{\bf Proof of Theorem \ref{T:0.21}:} We start with several remarks.
We again assume (without loss of generality) that a picked parameter
$\cE\in\cS_0$ satisfies condition ${\bf C}_{1-3}$ at
$t_1,\ldots,t_{n-\ell}$ and conditions \eqref{5.1} at the remaining
$\ell$ interpolation nodes. Under these non-restrictive assumptions we
will show that the function $w={\bf T}_\Theta[\cE]$ belongs to the class
$\cS_{\kappa-\ell}$. Throughout the proof, we shall be using partitionings
\eqref{2.11a}, \eqref{2.11d}. Note that by Theorem \ref{T:0.15},
the block $\tP_{22}$ is necessarily negative semidefinite. Then by Lemma
\ref{L:5.1}, ${\rm sq}_-P_{11}=\kappa-\ell$. Furthermore, since $\cE$
meets condition ${\bf C}_{1-3}$ at $t_1,\ldots,t_{n-\ell}$, the function
$w={\bf T}_\Theta[\cE]$ satisfies interpolation conditions \eqref{0.13} at
each of these points. Then by Remark \ref{R:0.6}, $w$ has at least
${\rm sq}_-P_{11}=\kappa-\ell$ negative squares.

\medskip

It remains to show that it has at most $\kappa-\ell$ negative squares.
This will be done separately for the cases when $\tP_{22}$ is negative
definite and when $\tP_{22}$ is negative semidefinite and singular.

\medskip

Conditions \eqref{5.1} mean that $\cE$ is a solution of a
boundary Nevanlinna--Pick interpolation problem with the data set
consisting of $\ell$ interpolation nodes $t_i$, unimodular numbers
$\q_i$ and nonnegative numbers ${\mathbb
P}_{ii}=-{\displaystyle\frac{\widetilde{p}_{ii}}{|\te_{i}|^2}}$
for $i=n-\ell+1,\ldots,n$. The Pick matrix ${\mathbb P}$ of the problem
is defined in \eqref{5.3}.

\medskip

{\bf Case 1: $\tP_{22}<0$:} In this case ${\mathbb P}>0$ (by \eqref{5.7})
and  by the first statement in Theorem \ref{T:2.6}, $\cE$ admits a
representation
\begin{equation}
\cE={\bf T}_{\widehat{\Theta}}[\widehat{\cE}]
\label{5.9}
\end{equation}
for some $\widehat{\cE}\in\cS_0$ where, according to \eqref{0.17}, the
coefficient matrix $\widehat{\Theta}$ in \eqref{5.9} is of the form
\begin{equation}
\widehat{\Theta}(z)=I_{2}+(z-\mu)\left[\begin{array}{c}M \\
E_2\end{array}\right]
(zI-T_2)^{-1}{\mathbb P}^{-1}(I-\mu T_2^*)^{-1}\left[\begin{array}{cc}
M^* & -E_2^*\end{array}\right]
\label{5.10}
\end{equation}
where the matrices
\begin{equation}
T_2={\rm diag} \, (t_{n-\ell+1}, \, \ldots, \, t_n),\quad
E_2=\begin{bmatrix} 1 & \ldots & 1\end{bmatrix}
\label{5.10a}
\end{equation}
are exactly the same as in \eqref{2.11a}, \eqref{2.11d}) and
\begin{equation}
M=\begin{bmatrix} \q_{n-\ell+1} & \q_{n-\ell+2} & \ldots &
\q_{n} \end{bmatrix}
\label{5.11}
\end{equation}
Self-evident equalities
$$
\begin{bmatrix}\q_i \\ 1\end{bmatrix} \cdot \frac{1}{z-t_i} \cdot
t_i\te_i=-\begin{bmatrix}\tc_i \\ \te_i\end{bmatrix}
\cdot \frac{1}{1-z\bar{t}_i}\quad (i=n-\ell+1,\ldots,n)
$$
can be written in the matrix form as
\begin{equation}
\begin{bmatrix}M \\ E_2\end{bmatrix}(zI-T_2)^{-1}{\bf C}=
-\begin{bmatrix}\tC_2 \\ \tE_2\end{bmatrix}(I-z T_2^*)^{-1}
\label{5.10b}
\end{equation}
where ${\bf C}$ is defined in \eqref{5.7}, whereas
$$
\tE_2=\begin{bmatrix}\te_{n-\ell+1} & \ldots &
\te_n\end{bmatrix}\quad\mbox{and}\quad
\tC_2=\begin{bmatrix}\tc_{n-\ell+1} & \ldots & \tc_n\end{bmatrix}
$$
are the matrices from the two last partitionings in \eqref{2.11d}.
On account of \eqref{5.10b} and \eqref{5.7}, we rewrite the formula
\eqref{5.10} as
$$
\widehat{\Theta}(z)=
I_{2}-(z-\mu)\left[\begin{array}{c}\tC_2 \\ \tE_2\end{array}\right]
(I-z T_2^*)^{-1}\tP_{22}^{-1}(\mu I-T_2)^{-1}
\left[\begin{array}{cc}\tC_2^* & -\tE_2^*\end{array}\right].
$$
Then its inverse can be represented as
$$
\widehat{\Theta}(z)^{-1}=
I_{2}+(z-\mu)\left[\begin{array}{c}\tC_2 \\ \tE_2\end{array}\right]
(I- \mu T_2^*)^{-1}\tP_{22}^{-1}(zI-T_2)^{-1}
\left[\begin{array}{cc}\tC_2^* & -\tE_2^*\end{array}\right]
$$
and coincides with the function $\widetilde{\Theta}^{(2)}$
from \eqref{2.11c}. Therefore, by Lemma \ref{L:2.5},
\begin{equation}
\Theta(z)=\Theta^{(1)}(z)\widehat{\Theta}(z)^{-1}
\label{5.12}
\end{equation}
where $\Theta^{(1)}$ is given in \eqref{2.11b}. Note that
\begin{equation}
\Theta^{(1)}\in{\mathcal W}_{\kappa_1}\quad\mbox{where} \; \;
\kappa_1={\rm sq}_-P_{11}=\kappa-\ell.
\label{5.13}
\end{equation}
Substituting \eqref{5.9} into \eqref{0.24} (that is, representing $w$ as a
result of composition of two linear fractional transformations) and taking
into account \eqref{5.12} we get
$$
w:={\bf T}_{\Theta}[\cE]=
{\bf T}_{\Theta}[{\bf T}_{\widehat{\Theta}}[\widehat{\cE}]]=
{\bf T}_{\Theta\widehat{\Theta}}[\widehat{\cE}]=
{\bf T}_{\Theta^{(1)}}[\widehat{\cE}].
$$
Since $\cE\in\cS_0$ and due to \eqref{5.12}, the last equality
guarantees (by Remark \ref{R:dop}) that $w$ has at most
$\kappa_1=\kappa-\ell$ negative squares
which completes the proof of Case 1.

\medskip

{\bf Case 2: $\tP_{22}\le 0$ is singular:} In this case ${\mathbb P}$ is
positive semidefinite and singular (again, by \eqref{5.7}) and  by the
second statement in Theorem \ref{T:2.6}, $\cE$ admits a representation
\begin{equation}
\cE(z)=\frac{x^*(I-zT_2^*)^{-1}E_2^*}{x^*(I-zT_2^*)^{-1}M^*}
\label{5.14}
\end{equation}
where $x$ is any nonzero vector such that ${\mathbb P}x=0$.
Letting $y:={\bf C}^{-1}x$ we have (due to \eqref{5.7})
\begin{equation}
\tP_{22}y=0
\label{5.15}
\end{equation}
and, on account of \eqref{5.10b}, we can rewrite \eqref{5.14} as
\begin{equation}
\cE(z)=\frac{y^*{\bf C}^*(I-zT_2^*)^{-1}E_2^*}{
y^*{\bf C}^*(I-zT_2^*)^{-1}M^*}
=\frac{y^*(zI-T_2)^{-1}\tE_2^*}{y^*(zI-T_2)^{-1}\tC_2^*}.
\label{5.16a}
\end{equation}
Since $\cE$ is a finite Blaschke product (again by the
second statement in Theorem \ref{T:2.6}) it satisfies the symmetry
relation $\cE(z)=(\overline{\cE(1/\bar{z})})^{-1}$ which together with
\eqref{5.16a} gives another representation for $\cE$:
\begin{equation}
\cE(z)=\frac{\tC_2(I-z T_2^*)^{-1}y}{\tE_2(I-z T_2^*)^{-1}y}.
\label{5.16}
\end{equation}
We will use the latter formula and \eqref{4.7a} to get
an explicit expression for the kernel $K_w(z,w)$. Setting
$$
u(z)=\tC_2(I-z T_2^*)^{-1}y\quad \mbox{and}\quad
v(z)=\tE_2(I-z T_2^*)^{-1}y
$$
for short and making use of the second Stein identity in \eqref{2.3}
we have
\begin{eqnarray*}
v(\zeta)^*v(z)-u(\zeta)^*u(z)&=&y^*(I-\bar{\zeta} T_2)^{-1}
\left[\tE_2^*\tE_2-\tC_2^*\tC_2\right](I-z T_2^*)^{-1}y\\
&=&y^*(I-\bar{\zeta} T_2)^{-1}
\left[\tP_{22}-T_2\tP_{22}T_2^*\right](I-z T_2^*)^{-1}y
\end{eqnarray*}
which reduces, due to \eqref{5.15}, to
$$
v(\zeta)^*v(z)-u(\zeta)^*u(z)=-(1-z\bar{\zeta})
y^*(I-\bar{\zeta} T_2)^{-1}T_2\tP_{22}T_2^*(I-z T_2^*)^{-1}y.
$$
Upon dividing both parts in the latter equality by
$(1-z\bar{\zeta})v(z)v(\zeta)^*$
we arrive at
\begin{equation}
\frac{1-\cE(\zeta)^*\cE(z)}{1-\bar{\zeta}z}=-
\frac{y^*}{v(\zeta)^*}(I-\bar{\zeta}
T_2)^{-1}T_2\tP_{22}T_2^*(I-z T_2^*)^{-1}\frac{y}{v(z)}.
\label{5.18}
\end{equation}
Next, we substitute the explicit formula \eqref{5.16} for $\cE$ into
\eqref{4.3} to get
\begin{eqnarray}
\Psi(z)&=&(zI-T)^{-1}\left(\tE^*-\tC^*\cE(z)\right)\nonumber\\
&=&(zI-T)^{-1}(\tE^*\tE_2-\tC^*\tC_2)(I-zT_2^*)^{-1}\cdot \frac{y}{v(z)}.
\label{5.19}
\end{eqnarray}
Substituting partitionings \eqref{2.11a}, \eqref{2.11d} into the
Stein identity \eqref{2.3} and comparing the right block entries we get
$$
\begin{bmatrix}\tP_{12} \\ \tP_{22}\end{bmatrix}
-T\begin{bmatrix}\tP_{12} \\
\tP_{22}\end{bmatrix}T_2^*=\tE\tE_2^*-\tC\tC_2^*
$$
which implies
\begin{eqnarray*}
&&\left(zI-T\right)^{-1}\left\{\tE\tE_2^*-\tC\tC_2^*\right\}
\left(I-zT_2^*\right)^{-1}\\
&&=\left(zI-T\right)^{-1}\begin{bmatrix}
\tP_{12} \\ \tP_{22}\end{bmatrix}
+\begin{bmatrix}\tP_{12} \\ \tP_{22}\end{bmatrix}T_2^*(I-zT_2^*)^{-1}.
\end{eqnarray*}
Now we substitute the last equality into \eqref{4.3} and take into
account \eqref{5.15} to get
$$
\Psi(z)=\left(zI-T\right)^{-1}\begin{bmatrix}
\tP_{12} \\ 0\end{bmatrix}\cdot \frac{y}{v(z)}
+\begin{bmatrix}\tP_{12} \\ \tP_{22}\end{bmatrix}T_2^*(I-zT_2^*)^{-1}
\cdot \frac{y}{v(z)}
$$
On account of partitionings \eqref{2.11a}, the latter equality leads us
to
\begin{eqnarray}
\Psi(\zeta)^*P\Psi(z)&=&\frac{y^*}{v(\zeta)^*}
\left(\tP_{12}^*(\bar{\zeta}I-T_1^*)^{-1}P_{11}(zI-T_1)^{-1}\tP_{12}\right.
\nonumber\\
&&\left.+(I-\bar{\zeta} T_2)^{-1}T_2\tP_{22}T_2^*(I-z
T_2^*)^{-1}\right)\frac{y}{v(z)}
\label{5.20}
\end{eqnarray}
Upon substituting \eqref{5.18} and \eqref{5.20} into \eqref{4.7a}
we get
$$
\frac{1-w(\zeta)^*w(z)}{1-\bar{\zeta}z}=\frac{y^*}{V_\cE(\zeta)^*v(\zeta)^*}
\cdot \tP_{12}^*(\bar{\zeta}I-T_1^*)^{-1}P_{11}(zI-T_1)^{-1}\tP_{12}\cdot
\frac{y}{V_\cE(z)v(z)}.
$$
Thus, the kernel $K_w(z,\zeta)$ admits a representation
$$
K_w(z,\zeta)=R(\zeta)^*P_{11}R(z)\quad\mbox{where}\quad
R(z)=\frac{y\tP_{21}T_1^*(I-zT_1^*)^{-1}}{v(z)V_\cE(z)}
$$
and thus, $$
{\rm sq}_-K_w\le {\rm sq}_-P_{11}=\kappa-\ell
$$
which completes the proof of the theorem.\qed
\begin{Rk}
{\rm At this point Theorem \ref{T:0.9} is completely proved: the necessity
part follows from Theorem \ref{T:2.3} and from the necessity part in
Theorem \ref{T:2.2}; the sufficiency part follows (as was explained in
introduction) from Corollary \ref{C:0.22} and Theorem \ref{T:0.21} which
have been  already proved.}
\label{R:6:2}
\end{Rk}
\begin{Rk}
{\rm We also proved the sufficiency part in Theorem \ref{T:2.2}
when the Pick matrix $P$ is invertible}.
\label{R:6.3}
\end{Rk}
Indeed, in this case, every solution $w$ to the FMI \eqref{2.2u} is of the
form \eqref{2.9u}, by Theorem \ref{T:2.3}. But every function of this form
solves Problem \ref{P:0.6}, by Theorem \ref{T:0.9}.

\section{The degenerate case}
\setcounter{equation}{0}

In this section we study Problem \ref{P:0.6} in the case when the Pick
matrix $P$ of the problem (defined in \eqref{0.10}) is singular. In the
course of the study we will prove Theorem \ref{T:0.8} and will complete
the proof of Theorem \ref{T:2.2}.
\begin{Tm}
Let the Pick matrix $P$ defined in \eqref{0.10} be singular with
${\rm rank} \, P=\ell<n$. Then there is a unique generalized Schur
function $w$ such that
\begin{equation}
{\rm sq}_-{\bf K}_w(z,\zeta)=\kappa
\label{7.1}
\end{equation}
where ${\bf K}_w(z,\zeta)$ is the kernel defined in \eqref{2.1u}.
Furthermore,
\begin{enumerate}
\item This unique function $w$ is the ratio of two finite Blaschke
products
\begin{equation}
w(z)=\frac{B_1(z)}{B_2(z)}
\label{7.2}
\end{equation}
with no common zeroes and such that
\begin{equation}
{\rm deg} \, B_1+{\rm deg} \, B_2={\rm rank} \, P.
\label{7.3}
\end{equation}
\item This unique function $w$ belongs to the generalized Schur class
$\cS_{\kappa^\prime}$ where $\kappa^\prime={\rm deg} \, B_2\le \kappa$
and satisfies conditions
\begin{equation}
d_w(t_i)\le \gamma_i\quad\mbox{and}\quad w(t_i)=w_i\qquad (i=1,\ldots,n)
\label{7.4}
\end{equation}
at all but $\kappa-\kappa^\prime$ interpolation nodes (that is, $w$
is a solution to Problem \ref{P:0.6}).
\item The function $w$ satisfies conditions
$$
d_w(t_i)=\gamma_i\quad\mbox{and}\quad w(t_i)=w_i
$$
at at least $n-{\rm rank} \, P$ interpolation nodes.
\end{enumerate}
\label{T:7.1}
\end{Tm}
{\bf Proof:} Without  loss of generality we can assume
that the top $\ell\times \ell$ principal submatrix $P_{11}$ of $P$ is
invertible and has $\kappa$ negative eigenvalues.
We consider conformal partitionings
\begin{equation}
T=\begin{bmatrix}T_1 & 0 \\ 0 & T_2
\end{bmatrix},\quad E=\begin{bmatrix}E_1 & E_2\end{bmatrix},\quad
C=\begin{bmatrix}C_1 & C_2\end{bmatrix}
\label{7.5}
\end{equation}
and
\begin{equation}
P=\begin{bmatrix} P_{11} & P_{12} \\ P_{21} & P_{22}\end{bmatrix},\quad
\det \, P_{11}\neq 0,\quad {\rm sq}_-P_{11}=\kappa={\rm sq} \, P.
\label{7.6}
\end{equation}
Since ${\rm rank} \, P_{11}={\rm rank} \, P$, it follows that
$P_{22}-P_{21}P_{11}^{-1}P_{12}$ the Schur complement of $P_{11}$ in $P$,
is the zero matrix, i.e.,
\begin{equation}
P_{22}=P_{21}P_{11}^{-1}P_{12}.
\label{7.7}
\end{equation}
Furthermore, it is readily seen that the $i$-th row of the block $P_{21}$
in \eqref{7.6} can be written in the form
$$
{\bf e}_i^* P_{21}=\begin{bmatrix}{\displaystyle
\frac{1-w_{\ell+i}^*w_1}{1-\bar{t}_{\ell+i}t_1}} & \ldots &
{\displaystyle\frac{1-w_{\ell+i}^*w_{\ell}}
{1-\bar{t}_{\ell+i}t_{\ell}}}
\end{bmatrix}=\left(E_1-w_{\ell+i}^*C_1\right)\left(
I-\bar{t}_{\ell+i}T_1\right)^{-1}
$$
and similarly, the $j$-th column in $P_{12}$ is equal to
\begin{equation}
P_{12}{\bf e}_j=\left(I-t_{\ell+j}T_1^*\right)^{-1}
\left(E_1^*-C_1^*w_{\ell+j}\right)
\label{7.8}
\end{equation}
(recall that ${\bf e}_j$ stands for the $j$-th column of the identity
matrix of an appropriate size). Taking into account that the $ij$-th entry
in $P_{22}$ is equal to ${\displaystyle\frac{1-w_{\ell+i}^*w_{\ell+j}}
{1-\bar{t}_{\ell+i}t_{\ell+j}}}$ (if $i\neq j$) or to $\gamma_{\ell+i}$
(if $i=j$) we write the equality \eqref{7.7} entrywise and get the
equalities
\begin{equation}
\frac{1-w_{i}^*w_j}{1-\bar{t}_it_j}=
\left(E_1-w_{i}^*C_1\right)\left(
I-\bar{t}_{i}T_1\right)^{-1}P_{11}^{-1}\left(
I-t_jT_1^*\right)^{-1}\left(E_1^*-{w}_{j}C_1^*\right)
\label{7.9}
\end{equation}
for $i\neq j\in\{\ell+1,\ldots,n\}$ and the equalities
\begin{equation}
\gamma_{i}=\left(E_1-w_{i}^*C_1\right)\left(
I-\bar{t}_{i}T_1\right)^{-1}P_{11}^{-1}\left(
I-t_iT_1^*\right)^{-1}\left(E_1^*-{w}_{i}C_1^*\right)
\label{7.10}
\end{equation}
for $i=\ell+1,\ldots,n$. The rest of the proof is broken into a number
of steps.

\medskip

{\bf Step 1:} {\em If $w$ is a meromorphic function such that \eqref{7.1}
holds, then it is necessarily of the form
\begin{equation}
w={\bf T}_{\Theta^{(1)}}[\cE]:=
\frac{\Theta^{(1)}_{11}\cE+\Theta^{(1)}_{12}}
{\Theta^{(1)}_{21}\cE+\Theta^{(1)}_{22}}
\label{7.11}
\end{equation}
for some Schur function $\cE\in\cS_0$, where $\Theta^{(1)}$ is given
in \eqref{2.11b}}.

\medskip

{\bf Proof of Step 1:} Write the kernel ${\bf K}_w(z,\zeta)$ in the block
form as
\begin{equation}
{\bf K}_w(z,\zeta)=
\begin{bmatrix}P_{11} & P_{12} & F_1(z) \\
P_{21} & P_{22} & F_2(z)\\
F_1(\zeta)^* & F_2(\zeta)^* & K_w(z,\zeta)\end{bmatrix}
\label{7.13}
\end{equation}
where $F_1$ and $F_2$ are given in \eqref{2.6m}. The kernel
\begin{eqnarray*}
{\bf K}^{1}_w(z,\zeta)&:=&\begin{bmatrix}
P_{11} & F_1(z) \\ F_1(\zeta)^* & K_w(z,\zeta)\end{bmatrix}\nonumber\\
&=&\begin{bmatrix}P_{11} &
(I-zT_1^*)^{-1}(E_1^*-C_1^*w(z)) \\
(E_1-w(\zeta)^*C_1)(I-\bar{\zeta}T_1)^{-1} & K_w(z,\zeta)\end{bmatrix}
\end{eqnarray*}
is contained in ${\bf K}_w(z,\zeta)$ as a principal submatrix and
therefore, ${\rm sq}_-{\bf K}^{1}_w\le \kappa$. On the other hand,
${\bf K}^{1}_w$ contains $P_{11}$ as a principal submatrix and therefore
${\rm sq}_-{\bf K}^{1}_w\ge {\rm sq}_-P_{11}=\kappa$. Thus,
\begin{equation}
{\rm sq}_-{\bf K}^{1}_w=\kappa.
\label{7.14}
\end{equation}
Recall that $P_{11}$ is an invertible Hermitian matrix with
$\kappa$ negative eigenvalues and satisfies the first Stein identity in
\eqref{2.3}. Then we can apply Theorem
\ref{T:2.3} (which is already proved for the case when the Pick
matrix is invertible) to the FMI \eqref{7.14}. Upon this
application we conclude that $w$ is of the form \eqref{7.11}
with some $\cE\in\cS_0$ and $\Theta^{(1)}$ of the form \eqref{2.11b}

\medskip

{\bf Step 2:} {\em Every function of the form \eqref{7.11} solves
the following truncated Problem \ref{P:0.6}: it belongs
to the generalized Schur class
$\cS_{\kappa^\prime}$ for some $\kappa^\prime\le \kappa$
and satisfies conditions
$$
d_w(t_i)\le \gamma_i\quad\mbox{and}\quad w(t_i)=w_i\qquad
(i=1,\ldots,\ell)
$$
at all but $\kappa-\kappa^\prime$ interpolation nodes.}

\medskip

{\bf Proof of Step 2:} The Pick matrix for the indicated truncated
interpolation problem is $P_{11}$ which is invertible and has $\kappa$
negative eigenvalues. Thus, we can apply Theorem \ref{T:0.9} (which is
already proved for the nondegenerate case) to get the desired statement.

\medskip

The rational function $\Theta^{(1)}$ is analytic
and $J$-unitary at $t_i$ for every $i=\ell+1,\ldots,n$. Then we can
consider the numbers $a_i$ and $b_i$ defined by
\begin{equation}
\begin{bmatrix} a_i \\ b_i\end{bmatrix}=\Theta^{(1)}(t_i)^{-1}
\begin{bmatrix} w_i \\ 1\end{bmatrix} \quad\mbox{for} \; \;
i=\ell+1,\ldots,n.
\label{7.15}
\end{equation}
It is clear from \eqref{7.15} that $|a_i|+|b_i|>0$. Furthermore,

\bigskip

{\bf Step 3:} {\em It holds that
\begin{equation}
|a_i|=|b_i|\neq 0\quad\mbox{and}\quad \frac{a_i}{b_i}=\frac{a_j}{b_j}
\quad \mbox{for}\quad i,j=\ell+1,\ldots,n.
\label{7.16}
\end{equation}}
{\bf Proof of Step 3:} Let $i\in\{\ell+1,\ldots,n\}$.
Since the matrix $\Theta^{(1)}(t_i)^{-1}$ is $J$-unitary and since
$|w_i|=1$, we conclude from \eqref{7.15} that
\begin{eqnarray}
|a_i|^2-|b_i|^2=\begin{bmatrix} a_i^* & b_i^*\end{bmatrix}J
\begin{bmatrix} a_i \\ b_i\end{bmatrix}&=&
\begin{bmatrix} w_i^* & 1 \end{bmatrix}\Theta^{(1)}(t_i)^{-*}J
\Theta^{(1)}(t_i)^{-1}\begin{bmatrix} w_i \\ 1\end{bmatrix}\nonumber\\
&=&\begin{bmatrix} w_i^* & 1 \end{bmatrix}J\begin{bmatrix} w_i \\
1\end{bmatrix}=|w_i|^2-1=0.\label{7.19}
\end{eqnarray}
Thus, $|a_i|=|b_i|$ and, since $|a_i|+|b_i|>0$, the first statement in
\eqref{7.16} follows. Similarly to \eqref{7.19}, we have
\begin{equation}
a_i^*a_j-b_i^*b_j=\begin{bmatrix} w_i^* & 1
\end{bmatrix}\Theta^{(1)}(t_i)^{-*}J
\Theta^{(1)}(t_j)^{-1}\begin{bmatrix} w_j \\ 1\end{bmatrix}
\label{7.20}
\end{equation}
for every choice of $i,j\in\{\ell+1,\ldots,n\}$. By a virtue of formula
\eqref{1.299},
\begin{equation}
\frac{\Theta^{(1)}(\zeta)^{-*}J\Theta^{(1)}(z)^{-1}-J}{1-z\bar{\zeta}}=
\left[\begin{array}{r}C_1 \\ -E_1\end{array}\right]
(I-\bar{\zeta}T_1)^{-1}
P_{11}^{-1}(I-zT_1^*)^{-1}\begin{bmatrix}C_1^* & -E_1^*\end{bmatrix}.
\label{7.21}
\end{equation}
Substituting the latter formula (evaluated at $\zeta=t_i$ and $z=t_j$)
into the right hand side expression in \eqref{7.20} and taking into
account that $\begin{bmatrix} w_i^* & 1 \end{bmatrix}J\begin{bmatrix} w_j
\\ 1\end{bmatrix}=w_i^*w_j-1$, we get
\begin{eqnarray*}
a_i^*a_j-b_i^*b_j&=&
w_i^*w_j-1+(1-\bar{t}_it_j)\left(E_1-w_{i}^*C_1\right)\left(
I-\bar{t}_{i}T_1\right)^{-1}P_{11}^{-1}\\
&&\qquad\qquad\qquad\times\left(
I-t_jT_1^*\right)^{-1}\left(E_1^*-{w}_{j}C_1^*\right).
\end{eqnarray*}
The latter expression is equal to zero, by \eqref{7.9}. Therefore,
$a_i^*a_j=b_i^*b_j$ and consequently,
$$
\frac{a_j}{b_j}=\frac{b_i^*}{a_i^*}=\frac{a_i}{b_i}
$$
where the second equality holds since $|a_i|=|b_i|$.

\medskip

{\bf Step 4:} {\em Let $a_i$ and $b_i$ be defined as in \eqref{7.15}.
Then the row vectors
\begin{equation}
A=\begin{bmatrix} a_{\ell+1} & \ldots & a_n\end{bmatrix},\qquad
B=\begin{bmatrix} b_{\ell+1} & \ldots & b_n\end{bmatrix}
\label{7.22}
\end{equation}
can be represented as follows:
\begin{equation}
\begin{bmatrix}A \\ B\end{bmatrix}=\begin{bmatrix}C \\ E\end{bmatrix}
\left(\mu I-T\right)^{-1}\begin{bmatrix} -P_{11}^{-1}P_{12}
\\ I \end{bmatrix}\left(\mu I-T_2\right).
\label{7.23}
\end{equation}}
{\bf Proof of Step 4:} First we substitute the formula \eqref{2.13}
for the inverse of $\Theta^{(1)}$ into \eqref{7.15} to get
$$
\begin{bmatrix}a_{i} \\ b_{i}\end{bmatrix}=\begin{bmatrix} w_i \\
1\end{bmatrix}+(t_{i}-\mu)\left[\begin{array}{c}C_1 \\
E_1\end{array}\right](\mu I-T_1)^{-1}P_{11}^{-1}
(I-t_{i}T_1^*)^{-1}\left(E_1^*-C_1^*w_i\right)
$$
for $i=\ell+1,\ldots,n$ and then we make use of \eqref{7.8} and of the
vector ${\bf e}_i$ to write the
latter equalities in the form
$$
\begin{bmatrix}A \\ B\end{bmatrix}{\bf e}_i=
\begin{bmatrix} w_{\ell+i} \\ 1\end{bmatrix}-
\left[\begin{array}{c}C_1 \\
E_1\end{array}\right](\mu I-T_1)^{-1}P_{11}^{-1}P_{12}
{\bf e}_i(\mu-t_{\ell+i})
$$
for $i=1,\ldots,n-\ell$. Now we transform the right hand side
expression in the latter equality as follows
\begin{eqnarray*}
\begin{bmatrix}A \\ B\end{bmatrix}{\bf e}_i&=&
\begin{bmatrix}C_2 \\ E_2\end{bmatrix}{\bf e}_i-
\begin{bmatrix}C_1 \\ E_1\end{bmatrix}(\mu I-T_1)^{-1}P_{11}^{-1}P_{12}
\left(\mu I-T_2\right){\bf e}_i \\
&=&\left(\begin{bmatrix}C_2 \\ E_2\end{bmatrix}\left(\mu
I-T_2\right)^{-1}-
\begin{bmatrix}C_1 \\ E_1\end{bmatrix}(\mu I-T_1)^{-1}P_{11}^{-1}P_{12}
\right)\left(\mu I-T_2\right){\bf e}_i\\
&=&\begin{bmatrix}C \\ E\end{bmatrix}
\left(\mu I-T\right)^{-1}\begin{bmatrix} -P_{11}^{-1}P_{12}
\\ I \end{bmatrix}\left(\mu I-T_2\right){\bf e}_i
\end{eqnarray*}
and since the latter equality holds for every $i\in\{1,\ldots,n-\ell\}$,
\eqref{7.23} follows.
\begin{Rk}
{\rm Comparing \eqref{7.23} and \eqref{2.16} we conclude that
$$
\begin{bmatrix}A \\ B\end{bmatrix}=
\Theta^{(1)}(z)^{-1}\begin{bmatrix} C \\
E\end{bmatrix} (zI-T)^{-1}\left[\begin{array}{c}-P_{11}^{-1}P_{12} \\
1\end{array}\right](z I -T_2).
$$
By the symmetry principle,
$\Theta^{(1)}(z)^{-1}=J\Theta^{(1)}(1/\bar{z})^{*}J$ and thus, the latter
identity can be written equivalently as
$$
\left[\begin{array}{r}A \\ -B\end{array}\right](z I -T_2)^{-1}
=\Theta^{(1)}(/\bar{z})^*\left[\begin{array}{r}C \\ -E\end{array}\right]
(zI-T)^{-1}\left[\begin{array}{c}-P_{11}^{-1}P_{12} \\
1\end{array}\right]
$$
Taking adjoints and replacing $z$ by $1/\bar{z}$ in the resulting
identity we obtain eventually
\begin{equation}
(I -zT_2^*)^{-1}\begin{bmatrix}A^* & -B^*\end{bmatrix}=
\left[\begin{array}{cc}-P_{21}P_{11}^{-1} &1\end{array}\right]
(I-zT^*)^{-1}\begin{bmatrix}C^* & -E^*\end{bmatrix}\Theta^{(1)}(z).
\label{7.23a}
\end{equation}
}\label{R:dop7}
\end{Rk}
{\bf Step 5:} {\em A function $w$ of the form \eqref{7.11} satisfies
the FMI \eqref{7.1} only if the corresponding parameter $\cE$
is the unimodular constant}
\begin{equation}
\cE(z)\equiv \cE_0:=
\frac{a_{\ell+1}}{b_{\ell+1}}=\ldots=\frac{a_n}{b_{n}}.
\label{7.24}
\end{equation}
{\bf Proof of Step 5:} Let us  consider the Schur
complement ${\bf S}$ of the block $P_{11}$ in \eqref{7.13}:
$$
{\bf S}(z,\zeta)=\begin{bmatrix}P_{22} & F_2(z)\\
F_2(\zeta)^* & K_w(z,\zeta)\end{bmatrix}-
\begin{bmatrix}P_{21}\\ F_1(\zeta)^*\end{bmatrix}
P_{11}^{-1}\begin{bmatrix}P_{12} & F_1(z)\end{bmatrix}
$$
Since
$$
{\rm sq}_-{\bf K}_w={\rm sq}_-P_{11}+{\rm sq}_-{\bf S}=\kappa+{\rm
sq}_-{\bf S},
$$
it follows that the FMI \eqref{7.1} is equivalent to positivity of
${\bf S}$ on $\rho(w)\cap\D$:
\begin{equation}
{\bf S}(z,\zeta)\succeq 0.
\label{7.25}
\end{equation}
Since the $``11''$ block in ${\bf S}(z,\zeta)$ equals $P_{22}-
P_{21}P_{11}^{-1}P_{12}$ which is the zero matrix (by \eqref{7.7}),
the positivity condition \eqref{7.25} guarantees the the
nondiagonal entries in ${\bf S}$ vanish everywhere in $\D$:
\begin{equation}
F_2(z)-P_{21}P_{11}^{-1}F_1(z)\equiv 0.
\label{7.26}
\end{equation}
By \eqref{2.6t}, the latter identity can be written as
\begin{equation}
\begin{bmatrix}-P_{21}P_{11}^{-1} & I\end{bmatrix}
(I-zT^*)^{-1}(E^*-C^*w(z))\equiv 0.
\label{7.27}
\end{equation}
We already know from Step 1, that $w$ is of the form \eqref{7.11} for
some $\cE\in\cS_0$. Now we will show that \eqref{7.27} holds for
$w$ of the form \eqref{7.11} if and only if the corresponding parameter
$\cE$ is subject to
\begin{equation}
A^*\cE(z)\equiv B^*
\label{7.29}
\end{equation}
where $A$ and $B$ are given in \eqref{7.22}. Indeed, it is easily seen
that for $w$ of the form \eqref{7.11}, it holds that
$$
E^*-C^*w=\left(\Theta^{(1)}_{21}\cE+\Theta^{(1)}_{22}\right)^{-1}
\begin{bmatrix}
-C^* & E^*\end{bmatrix}\begin{bmatrix}\Theta^{(1)}_{11} &
\Theta^{(1)}_{12} \\ \Theta^{(1)}_{21} & \Theta^{(1)}_{22}\end{bmatrix}
\begin{bmatrix}\cE \\ 1 \end{bmatrix}
$$
and therefore, identity \eqref{7.27} can be written equivalently in terms
of the parameter $\cE$ as
$$
\begin{bmatrix}-P_{21}P_{11}^{-1} & I\end{bmatrix}(I-zT^*)^{-1}
\begin{bmatrix} C^* & -E^*\end{bmatrix}\Theta^{(1)}(z)
\begin{bmatrix}\cE(z) \\ 1\end{bmatrix}\equiv 0
$$
which is, due to \eqref{7.23a}, the same as
$$
(I-zT_2^*)^{-1}
\begin{bmatrix} A^* & B^*\end{bmatrix}J\begin{bmatrix}
\cE(z) \\ I \end{bmatrix}\equiv 0.
$$
The latter identity is clearly equivalent to \eqref{7.29}. Writing
\eqref{7.29} entrywise we get the system of equalities
$$
a_i^*\cE(z)\equiv b_i^* \qquad (i=\ell+1,\ldots,n).
$$
This system is consistent, by \eqref{7.16}, and it clearly admits a
unique solution $\cE_0$ defined as in \eqref{7.24}. Combining Step 1
and Step 5, we can already conclude that the FMI \eqref{7.1} has at most
one  solution: the only candidate is the function
\begin{equation}
w={\bf T}_{\Theta^{(1)}}[\cE_0]
\label{7.30}
\end{equation}
where $\cE_0$ is the unimodular constant defined in \eqref{7.24}.
The next step will show that this function indeed is a solution
to the FMI \eqref{7.1}.

\medskip

{\bf Step 6:} {\em The function \eqref{7.30} satisfies the FMI
\eqref{7.1} and interpolation conditions
\begin{equation}
d_w(t_i)=\gamma_i\quad\mbox{and}\quad w(t_i)=w_i\quad\mbox{for}\; \;
i=\ell+1,\ldots,n.
\label{7.31}
\end{equation}}
{\bf Proof of Step 6:} First we note that since $\Theta^{(1)}$ is a
rational $J$-inner function of McMillan degree $\ell$ and since $\cE_0$ is
a unimodular constant, the function $w$ of the form \eqref{7.30} is
a rational function of degree $\ell$ which is unimodular on $\T$.
Therefore, $w$ is the ratio of two finite Blachke products satisfying
\eqref{7.3}. Since $w$ belongs to $\cS_{\kappa^\prime}$ (by Step 2), it
has $\kappa^\prime$ poles inside $\D$ and thus, the denominator $B_2$ in
\eqref{7.2} is a finite Blachke product of order $\kappa^\prime$.

\smallskip

It was shown in the proof of Step 5 that
equation \eqref{7.29} is equivalent to \eqref{7.26})
and thus, for the function $w$ of the form \eqref{7.30}, it holds
that
\begin{equation}
F_2(z)\equiv P_{21}P_{11}^{-1}F_1(z)
\label{7.32}
\end{equation}
which is the same, due to definitions \eqref{2.6m}, as
\begin{equation}
(I-zT_2^*)^{-1}(E_2^*-C_2^*w(z))\equiv
P_{21}P_{11}^{-1}(I-zT_1^*)^{-1}(E_1^*-C_1^*w(z)).
\label{7.33}
\end{equation}
Next we show that for $w$ of the form
\eqref{7.30} it holds that
\begin{equation}
K_w(z,\zeta)\equiv F_1(\zeta)^*P_{11}^{-1}F_1(z)
\label{7.34}
\end{equation}
or, which is the same,
\begin{equation}
\frac{1-w(\zeta)^*w(z)}{1-\bar{\zeta}z}\equiv
(E_1-w(\zeta)^*C_1)(I-\bar{\zeta}T_1)^{-1}P_{11}^{-1}
(I-zT_1^*)^{-1}(E_1^*-C_1^*w(z)).
\label{7.35}
\end{equation}
Indeed, on account of \eqref{7.21},
\begin{eqnarray}
&&(E_1-w(\zeta)^*C_1)(I-\bar{\zeta}T_1)^{-1}P_{11}^{-1}
(I-zT_1^*)^{-1}(E_1^*-C_1^*w(z))\nonumber \\
&&=\begin{bmatrix}w(\zeta)^*& 1\end{bmatrix}
\frac{\Theta^{(1)}(\zeta)^{-*}J\Theta^{(1)}(z)^{-1}-J}{1-z\bar{\zeta}}
\begin{bmatrix}w(z) \\ 1\end{bmatrix}\nonumber \\
&&=\frac{1-w(z)w(\zeta)^*}{1-z\bar{\zeta}}+
\begin{bmatrix}w(\zeta)^*& 1\end{bmatrix}
\frac{\Theta^{(1)}(\zeta)^{-*}J\Theta^{(1)}(z)^{-1}}{1-z\bar{\zeta}}
\begin{bmatrix}w(z) \\ 1\end{bmatrix}.\label{7.36}
\end{eqnarray}
Representation \eqref{7.30} is equivalent to
$$
\begin{bmatrix}w(z) \\ 1\end{bmatrix}=\Theta^{(1)}(z)\begin{bmatrix}\cE_0
\\ 1\end{bmatrix}\frac{1}{v(z)},\quad\mbox{where} \; \;
v(z)=\Theta^{(1)}_{21}(z)\cE_0+
\Theta^{(1)}_{22}(z),
$$
and therefore,
$$
\begin{bmatrix}w(\zeta)^*& 1\end{bmatrix}\Theta^{(1)}(\zeta)^{-*}J
\Theta^{(1)}(z)^{-1}\begin{bmatrix}w(z) \\ 1\end{bmatrix}=
\frac{|\cE_0|^2-1}{v(z)v(\zeta)^*}\equiv 0,
$$
since $|\cE_0|=1$. On account of this latter equality, \eqref{7.36}
implies \eqref{7.34}. By \eqref{7.7}, \eqref{7.32} and \eqref{7.34},
the kernel ${\bf K}_w(z,\zeta)$ defined in \eqref{2.1u}
and partitioned as in \eqref{7.13}, can be represented also in the form
$$
{\bf K}_w(z,\zeta)=\left[\begin{array}{c} P_{11} \\ P_{21} \\
F_1(\zeta)^* \end{array}\right]P_{11}^{-1}
\left[\begin{array}{ccc} P_{11} & P_{12} & F_1(z)\end{array}\right]
$$
and the latter representation implies that ${\rm sq}_-{\bf K}_w={\rm
sq}_-P_{11}=\kappa$, i.e., that $w$ of the form \eqref{7.30} satisfies
the FMI \eqref{7.1}. It remains to check that $w$ satisfies interpolation
conditions \eqref{7.31}. Since $w$ is a ratio of two finite Blaschke
products, it is analytic on $\T$. Let $t_i$ ($\ell<i\le n$) be an
interpolation node. Comparing the residues at $z=t_i$ of both parts in
the identity \eqref{7.33} we get
$$
-t_i{\bf e}_i{\bf e}_i^*\left(E_2^*-C_2^*w(t_i)\right)=0
$$
which is equivalent to
$$
1-w_i^*w(t_i)=0
$$
or, since $|w_i|=1$, to the second condition in \eqref{7.31}. On the
other hand, letting $z, \, \zeta\to t_i$ in \eqref{7.35} and taking into
account that $w(t_i)=w_i$, we get
\begin{eqnarray*}
d_w(t_i)&=&(E_1-w(t_i)^*C_1)(I-\bar{t_i}T_1)^{-1}P_{11}^{-1}
(I-t_iT_1^*)^{-1}(E_1^*-C_1^*w(t_i))\\
&=&(E_1-w_i^*C_1)(I-\bar{t_i}T_1)^{-1}P_{11}^{-1}
(I-t_iT_1^*)^{-1}(E_1^*-C_1^*w_i)
\end{eqnarray*}
which together with \eqref{7.10} implies the first condition in
\eqref{7.31}.

\medskip

The first statement of the Theorem is proved. Statement 2 follows by Step
2 and \eqref{7.31}: the function $w$
meets interpolation conditions \eqref{7.4} at all but
$\kappa-\kappa^\prime$ interpolation nodes (and all the exceptional
nodes are in $\{t_1,\ldots,t_\ell\}$). Statement 3 follows from
\eqref{7.31}.\qed
\begin{Rk}
{\rm Statement 2 in Theorem \ref{T:7.1} completes the proof of sufficiency
part in  Theorem \ref{T:2.2}: if $P$ is singular, then a (unique) solution
of the FMI \eqref{2.2u} solves Problem \ref{P:0.6}. }
\label{R:7.2}
\end{Rk}

\section{An example}
\setcounter{equation}{0}

In this section we present a numerical example illustrating the
preceding analysis. The data set of the problem is as follows:
\begin{equation}
t_1=1, \; \; t_2=-1, \; \; w_1=1, \; \; w_2=-1, \; \; \gamma_1=1, \; \;
\gamma_2=0.
\label{11.1}
\end{equation}
Then the matrices \eqref{0.18} take the form
$$
T=\left[\begin{array}{cr} 1 & 0 \\ 0 &-1\end{array}\right]\quad
\mbox{and}\quad
\begin{bmatrix}C\\ E\end{bmatrix}=
\left[\begin{array}{cr}1 & -1\\ 1 & 1\end{array}\right]
$$
and since ${\displaystyle\frac{1-w_1^*w_2}{1-\bar{t}_1t_2}}=1$ we have
also
$$
P=\begin{bmatrix}1 & 1 \\ 1 & 0\end{bmatrix}\quad\mbox{and}\quad
P^{-1}=\left[\begin{array}{cr} 0 & 1 \\ 1 &-1\end{array}\right].
$$
It is readily seen that $P$ is invertible and has one negative
eigenvalue. Thus, Problems \ref{P:0.4}, \ref{P:0.5} and \ref{P:0.6}
take the following form.

\medskip

{\bf Problem \ref{P:0.5}:} {\em Find all functions $w\in\cS_1$ such
that}
\begin{equation}
w(1)=1,\quad d_w(1)\le 1,\quad w(-1)=-1,\quad d_w(-1)\le 0.
\label{11.2}
\end{equation}
{\bf Problem \ref{P:0.4}:} {\em Find all functions $w\in\cS_1$
that satisfy conditions \eqref{11.2} with equalities in the second
and in the fourth conditions}.

\medskip

{\bf Problem \ref{P:0.6}:} {\em Find all functions $w$ such that
either
\begin{enumerate}
\item $w\in\cS_1$ and satisfies all the conditions in \eqref{11.2} or
\item $w\in\cS_0$ and satisfies the two first conditions in \eqref{11.2} or
\item $w\in\cS_0$ and satisfies the two last conditions in \eqref{11.2}.
\end{enumerate}}
Letting $\mu=i$, we get by the formula \eqref{0.17} for $\Theta$:
\begin{eqnarray*}
\Theta(z)&=&I_{2}+(z-i)\left[\begin{array}{cr} 1 & -1 \\ 1 &
1\end{array}\right]\begin{bmatrix}\frac{1}{z-1} & 0 \\ 0 &
\frac{1}{z+1}\end{bmatrix}
\left[\begin{array}{cr} 0 & 1 \\ 1 &-1\end{array}\right]
\begin{bmatrix}\frac{1}{1-i} & 0 \\ 0 &
\frac{1}{1+i}\end{bmatrix}\left[\begin{array}{rr} 1 & -1 \\ -1 &
-1\end{array}\right]\\
&=& \frac{1}{2(z^2-1)}\left[\begin{array}{cc}
(i-1)z^2+2(1+2i)z-1-i &(3i-1)z^2+2z+i-1 \\
(i+1)z^2-2z+1+3i & (1-i)z^2+2(2i-1)z+1+i\end{array}\right]
\end{eqnarray*}
and thus, by Theorem \ref{T:0.9}, all the solutions $w$ to Problem
\ref{P:0.6} are parametrized by the linear fractional formula
\begin{equation}
w(z)=\frac{\left[(i-1)z^2+2(1+2i)z-1-i\right]\cE(z)+(3i-1)z^2+2z+i-1}
{\left[(i+1)z^2-2z+1+3i\right]\cE(z)+ (1-i)z^2+2(2i-1)z+1+i}
\label{11.3}
\end{equation}
when the parameter $\cE$ runs through the Schur class $\cS_0$.
Furthermore, formula \eqref{2.2} in the present setting gives
$$
\left[\begin{array}{cc}\tc_1 & \tc_2\\ \te_1 & \te_2\end{array}\right]
=\left[\begin{array}{cr} 1 & -1 \\ 1 & 1\end{array}\right]
\begin{bmatrix}\frac{1}{i-1} & 0 \\ 0 &
\frac{1}{i+1}\end{bmatrix}
\left[\begin{array}{cr} 0 & 1 \\ 1 &-1\end{array}\right]
\begin{bmatrix}1-i & 0 \\ 0 & 1+i\end{bmatrix}=\left[\begin{array}{rr}
1 & 1-i \\ -1 & -1-i\end{array}\right]
$$
and since the diagonal entries of $P^{-1}$ are $\widetilde{p}_{11}=0$
and $\widetilde{p}_{11}=-1$, we also have
$$
\q_1:=\frac{\tc_1}{\te_1}=-1,\quad \q_2:=\frac{\tc_2}{\te_2}=i, \quad
\frac{\widetilde{p}_{11}}{|\te_1|^2}=0,\quad
\frac{\widetilde{p}_{22}}{|\te_2|^2}=-\frac{1}{2}.
$$
By Theorem \ref{T:0.11}, every function $w$ of the form \eqref{11.3}
also solves Problem \ref{P:0.5}, unless the parameter $\cE$
is subject to
\begin{equation}
\cE(1)=-1\quad\mbox{and}\quad d_\cE(1)=0
\label{11.4}
\end{equation}
or to
\begin{equation}
\cE(-1)=i\quad\mbox{and}\quad d_\cE(-1)\le \frac{1}{2}.
\label{11.5}
\end{equation}
On the other hand, Theorem \ref{T:0.10} tells us that every function $w$
of the form \eqref{11.3} solves Problem \ref{P:0.4}, unless the parameter
$\cE$ is subject to
$$
\cE(1)=-1\quad\mbox{and}\quad d_\cE(1)<\infty
$$
or to
$$
\cE(-1)=i\quad\mbox{and}\quad d_\cE(-1)<\infty.
$$
Thus, every parameter $\cE\in\cS_0$ satisfying conditions
\eqref{11.4} or \eqref{11.5} leads to a solution $w$ of Problem
\ref{P:0.6} which is not a solution to Problem \ref{P:0.5}.
For these special solutions, it looks curious to track
which conditions in \eqref{11.2} are satisfied and which are not. This
will also illustrate propositions 4 and 5 in Theorem \ref{T:0.20}.

\medskip

First we note that there is only one Schur function $\cE\equiv -1$
satisfying  conditions \eqref{11.4} (this is the case indicated in the
fifth part in Theorem \ref{T:0.20}). The corresponding function
$w$ obtained via \eqref{11.3}, equals
$$
w(z)=\frac{2iz^2-4iz+2i}{-2iz^2+4iz-2i}\equiv -1.
$$
It satisfies all the conditions in \eqref{11.2} but the first one.

\medskip

All other ``special'' solutions of Problem \ref{P:0.6} are exactly
all Schur functions satisfying the two first conditions in \eqref{11.2}.
Every such function does not satisfy at least one of the two last
conditions in \eqref{11.2}. We present several examples omitting
straightforward computations:

\medskip

{\bf Example 1:} The function
$$
\cE(z)=\frac{2iz+2}{(1-i)z-1-3i}
$$
belongs to $\cS_0$ and satisfies $\cE(-1)=i$ and
$d_{\cE}(-1)=\frac{1}{2}$ (i.e., it meets condition \eqref{0.29} at
$t_2$). Substituting this parameter into \eqref{11.3} we get the function
$$
w(z)=\frac{z-i}{iz+1-2i}
$$
which belongs to $\cS_0$ and satisfies (compare with \eqref{11.2})
$$
w(1)=1,\quad d_w(1)=1,\quad w(-1)=\frac{1+i}{3i-1},\quad d_w(-1)=\infty.
$$
{\bf Example 2:} The function
$$
\cE(z)=\frac{(3-i)z-(1+i)}{-(1+i)z+3i-1}
$$
belongs to $\cS_0$ and satisfies (as in Example 1) $\cE(-1)=i$ and
$d_{\cE}(-1)=\frac{1}{2}$.  Substituting this parameter into \eqref{11.3}
we get the function $w(z)\equiv 1$ which belongs to $\cS_0$ and satisfies
(compare with \eqref{11.2})
$$
w(1)=1,\quad d_w(1)=0,\quad w(-1)=1,\quad d_w(-1)=0.
$$
{\bf Example 3:} The function
$$
\cE(z)=\frac{\left[(3+i)z+1-i\right]e^{\frac{z-1}{z+1}}-2iz-2}
{-2(1+iz)e^{\frac{z-1}{z+1}}+(i-1)z+3i+1}
$$
belongs to $\cS_0$ and satisfies $\cE(-1)=i$ and
$d_{\cE}(-1)=\frac{1}{2}$. Substituting this parameter into \eqref{11.3}
we get the function
$$
w(z)=\frac{\left[(2-i)z-1\right]e^{\frac{z-1}{z+1}}-z+i}
{(z-i)e^{\frac{z-1}{z+1}}-iz+2i-1}
$$
which belongs to $\cS_0$ and fails to have a boundary nontangential limit
at $t_2=-1$.

\bibliographystyle{amsplain}
\providecommand{\bysame}{\leavevmode\hbox to3em{\hrulefill}\thinspace}

\end{document}